\documentclass[12pt,a4paper]{amsart}
\usepackage[T1]{fontenc}
\usepackage{lmodern}
\usepackage[a4paper,margin=28mm]{geometry}
\usepackage{amsmath,amssymb,amsthm,mathtools,mathrsfs}
\usepackage{enumitem}
\usepackage{microtype}
\usepackage{aliascnt}
\usepackage{tikz}
\usepackage{xcolor}
\usepackage[colorlinks,linkcolor=blue,anchorcolor=green,citecolor=red]{hyperref}
\usepackage[nameinlink,noabbrev]{cleveref}
\numberwithin{equation}{section}
\newtheorem{theorem}{Theorem}[section]
\newaliascnt{proposition}{theorem}
\newtheorem{proposition}[proposition]{Proposition}
\aliascntresetthe{proposition}
\newaliascnt{lemma}{theorem}
\newtheorem{lemma}[lemma]{Lemma}
\aliascntresetthe{lemma}
\newaliascnt{corollary}{theorem}
\newtheorem{corollary}[corollary]{Corollary}
\aliascntresetthe{corollary}
\theoremstyle{definition}
\newaliascnt{definition}{theorem}
\newtheorem{definition}[definition]{Definition}
\aliascntresetthe{definition}
\theoremstyle{remark}
\newaliascnt{remark}{theorem}
\newtheorem{remark}[remark]{Remark}
\aliascntresetthe{remark}
\crefname{theorem}{Theorem}{Theorems}
\crefname{proposition}{Proposition}{Propositions}
\crefname{lemma}{Lemma}{Lemmas}
\crefname{corollary}{Corollary}{Corollaries}
\crefname{definition}{Definition}{Definitions}
\crefname{remark}{Remark}{Remarks}
\newcommand{\X}{{\overline{X}}}
\newcommand{\A}{\mathcal A}
\newcommand{\B}{\mathcal B}
\newcommand{\G}{\mathcal G}
\newcommand{\M}{\mathcal M}
\newcommand{\E}{\mathcal E}

\newcommand{\ii}{\sqrt{-1}}
\newcommand{\RR}{\mathbb R}
\newcommand{\CC}{\mathbb C}
\newcommand{\ZZ}{\mathbb Z}
\newcommand{\dd}{\partial}
\newcommand{\db}{\bar\partial}
\newcommand{\Id}{\operatorname{Id}}
\newcommand{\tr}{\operatorname{tr}}
\newcommand{\Hom}{\operatorname{Hom}}
\newcommand{\End}{\operatorname{End}}
\newcommand{\ad}{\operatorname{ad}}
\newcommand{\diag}{\operatorname{diag}}
\newcommand{\rk}{\operatorname{rk}}
\newcommand{\supp}{\operatorname{supp}}

\newcommand{\pardeg}{\operatorname{pardeg}}
\newcommand{\parmu}{\operatorname{par\mu}}
\newcommand{\vol}{\operatorname{Vol}}

\newcommand{\KM}{\mathrm{KM}}

\DeclareMathOperator{\pa}{par}
\DeclareMathOperator{\Ad}{Ad}
\DeclareMathOperator{\re}{Re}
\DeclareMathOperator{\im}{Im}
\newcommand{\SU}{SU}
\newcommand{\SO}{SO}
\newcommand{\PU}{PU}
\newcommand{\U}{U}
\setlist[enumerate]{label=\textup{(\roman*)},leftmargin=2.2em,itemsep=3pt}
\hypersetup{pdftitle={Compactification of parabolic Hermitian--Einstein moduli spaces on Gauduchon surfaces},pdfauthor={Xilun Li and Gang Tian},pdfsubject={Parabolic bundles, singular instantons, and analytic compactification}}

\title[Compactification of parabolic Hermitian--Einstein moduli spaces]{Compactification of parabolic Hermitian--Einstein moduli spaces on complex surfaces}

\author{Xilun Li}
\address{The Abdus Salam International Centre for Theoretical Physics\\Trieste\\ Italy}
\email{xli@ictp.it}
\author{Gang Tian}
\address{Beijing International Center for Mathematical Research, Peking University, Beijing 100871, China}
\email{gtian@math.pku.edu.cn}
\date{}
\begin{document}
\maketitle
\begin{abstract}
Let $\X$ be a compact complex surface with a Gauduchon metric, and let $D\subset\X$ be a smooth connected divisor. We construct an Uhlenbeck compactification of the moduli space of rank-two parabolic Hermitian--Einstein connections with prescribed holomorphic determinant.
\end{abstract}
\tableofcontents

\section{Introduction}

\;\;\; Compactification plays a fundamental role in understanding the global geometry of moduli spaces. Uhlenbeck's theory provides a compactification for the moduli space of four-dimensional ASD connections \cite{Uhlenbeck82a,Uhlenbeck82b,DK}, and hence also gives a compactification for the moduli space of Hermitian--Einstein connections \cite{LubkeTeleman}. A natural question is to develop an analogous compactification theory for parabolic Hermitian--Einstein connections.

In the projective setting, Maruyama--Yokogawa constructed a compact projective moduli space of semistable parabolic sheaves using Gieseker-type stability \cite{MY92}. Balaji--Dey--Parthasarathi constructed an algebro-geometrical compactification of parabolic bundle moduli spaces on projective surfaces with rational weights \cite{BDP08}, where the compactified space has a similar structure as Uhlenbeck's compactness. The algebraic arguments cannot work for general Gauduchon surfaces, and in fact we cannot compactify the nonK\"{a}hler case in a pure complex-geometric way for some cases.

Let $(\X,J,\omega)$ be a compact complex surface equipped with a Gauduchon metric, and let $D\subset \X$ be a smooth divisor. We consider the parabolic bundles $\mathcal E_*=(\mathcal E,\mathcal L,\alpha_1,\alpha_2)$ of rank $2$ on $(\X,D)$, which admits a Hermitian--Einstein metric $H$ on $X:=\X\setminus D$ compatible with the parabolic structure. We fix the smooth structures $(E,L)$ of the holomorphic bundles $(\mathcal E,\mathcal L)$ and fix the determinant $\det\mathcal E\simeq\mathcal{P}$, where $\mathcal{P}$ is a fixed holomorphic line bundle. Then the isomorphism classes of such $(\mathcal E_*,H)$ form a moduli space $\mathcal M_h^{\mathrm{par,HE}}(E,L;\mathcal P)$ (\Cref{def:parHEmoduli}). Our goal is to give a compactification of such moduli space using the gauge theory.

There are two other related moduli spaces. The first one is the moduli space of polystable parabolic bundles $\M_{\mathrm{pst}}^{\mathrm{par}}(E,L;\mathcal P)$ (\Cref{def:polystable-moduli}). By the parabolic Kobayashi--Hitchin correspondence in the nonK\"{a}hler setting \cite{LiTian}, we can show these moduli spaces are bijective. In the classical non-parabolic case, the irreducible part of the moduli space of Hermitian--Einstein connections, which corresponds to the strictly stable part of the moduli space of polystable bundles, has a complex space structure, and the local structures have Kuranishi description \cite{Kobayashi87,LubkeTeleman}. However, in the nonK\"{a}hler setting, the reducible part may be of odd dimension \cite{Teleman15,BTT18}. Therefore the total moduli space is \emph{not} a complex-geometric object, compared with the fact the moduli space has a projective variety structure in the projective setting.

The second one is the moduli space of twisted ASD $PU(2)$-connections $\M^{\KM}_\beta(E,L;a_{P,\mathrm{par}})$ (\Cref{def:gaugemoduli}), where the $SU(2)$-case was studied by Kronheimer--Mrowka \cite{KM93}. They considered ASD connections on compact Riemannian four-manifold, which is singular along a smooth embedded surface. Their description used a model connection with prescribed holonomy and weighted Sobolev spaces, then they gave a compactification of the moduli space, which can be regarded as a singular version of Uhlenbeck--Donaldson compactification. A natural idea is to modify the Kronheimer--Mrowka's theory to the $PU(2)$-setting,  then identify the moduli spaces of parabolic Hermitian--Einstein connections and twisted ASD connections to get the compactification, since the parabolic setting shows a similar asymptotic behaviour as the twisted connections. Biquard proved a similar bijection in the $SU(2)$ case for projective surfaces \cite{Biquard96}.

Thus our first result is
\begin{theorem}[\Cref{thm:inverse-map}, \Cref{thm:parabolic-KH}]
    There exist natural well-defined bijections:
    \begin{align*}
        \mathcal M_h^{\mathrm{par,HE}}(E,L;\mathcal P)\simeq \M_{\mathrm{pst}}^{\mathrm{par}}(E,L;\mathcal P)\simeq \M^{\KM}_\beta(E,L;a_{P,\mathrm{par}}).
    \end{align*}
\end{theorem}

The moduli space $\M^{\KM}_\beta(E,L;a_{P,\mathrm{par}})\subset \mathcal A^{1,p}_\beta /\mathcal G^{2,p}_\beta$ is a quotient of a weighted Sobolev affine space, hence has a quotient topology induced by the weighted Sobolev norms. Under the bijections above, this also gives a topology to two other moduli spaces. It is plausible that the restriction of the pull-back topology on the strictly stable part of $\M_{\mathrm{pst}}^{\mathrm{par}}(E,L;\mathcal P)$ and the topology of complex space structure coincide, while we always use the pull-back topology throughout this article. Since we use the pull-back topology, the bijections are homeomorphisms.

Now we can apply the $PU(2)$ version of Kronheimer--Mrowka's compactness theorem, and construct suitable ideal parabolic connections and ideal topology. Our second result is
\begin{theorem}[\Cref{cor:compactification}, \Cref{prop:open-topology}, \Cref{def:idealparHE}]
    There exists a compact, Hausdorff, second-countable, metrizable topological space 
    \begin{align*}
        \overline{\mathcal M^{\mathrm{par,HE}}_h}(E,L;\mathcal P)\supset \mathcal M_h^{\mathrm{par,HE}}(E,L;\mathcal P)
    \end{align*}
    such that $M_h^{\mathrm{par,HE}}$ is an open subset and the restriction topology is same as the induced topology by the weighted Sobolev norms.

    Moreover, any boundary element has the form $(\mathcal F_*,H,\phi_F,Z_0,Z_1)$, where
\begin{enumerate}
\item $\mathcal F_*=(\mathcal F,M,\alpha_1,\alpha_2)$
is a rank-two parabolic holomorphic bundle on $(\X,D)$. 
\item $H$ is a compatible Hermitian--Einstein metric on $\mathcal F|_X$.
\item $\phi_F:\det\mathcal F\xrightarrow{\sim}\mathcal P$ is a holomorphic line-bundle isomorphism.
\item $Z_0$ and $Z_1$ are effective zero-cycles on
$\X$ satisfying
\begin{align*}
    Z_0|_X=Z_1|_X,\quad c_2(E)=c_2(F)+\deg Z_0,\quad \deg L&=\deg M+\deg Z_0-\deg Z_1.
\end{align*}
\end{enumerate}
\end{theorem}

We discuss the technical part of the paper. 
The correspondence between parabolic Hermitian--Einstein connections and twisted ASD connections needs arguments of two directions. One direction, from twisted ASD connections to parabolic HE, mainly follows from Biquard \cite{Biquard96}, where he studied the parabolic extension of holomorphic bundles on $\X\setminus D$ equipped Hermitian metrics with $L^p$ curvature for some $p>2$. Another direction, from parabolic HE to twisted ASD connections, needs to show weighted Sobolev estimates on the compatible HE connections. The main step of the estimate is to show any compatible HE connection has $L^2$ curvature. A similar result was proved by Biquard \cite{Biquard96} for projective surfaces, and by Li--Narasimhan \cite{LN99} for K\"{a}hler surfaces. Their arguments both depend on the specific construction of the HE connections by the heat flow method which cannot be naturally extended to the non-K\"{a}hler setting, however our arguments only depend on a weak compatible condition and Hermitian--Einstein condition.

We also need to deal with topological issues to modify the compactness of twisted ASD $SU(2)$ connections to prescribed determinant twisted ASD $U(2)$ connections. 
Since the projective unitary frame bundle is a $PU(2)$ bundle, we first extend to $PU(2)\cong SU(2)/\{\pm 1\}$ case. Since the local arguments carry over without essential modification, this extension is straightforward. In fact, Kronheimer--Mrowka briefly discussed the $SO(3)$ setting in \cite[\S~2(iv)]{KM93}, thus it suffices to deal with the details.
A further topological issue arises when lifting the projective limit back to a rank-two bundle with the prescribed determinant. Such a lift always exists by considering Steifel--Whitney class $w_2$. However, the projective gauge may not admit a lift, whose obstruction lies in $H^1(\overline X\setminus S;\mathbb Z_2)$. By tensoring a $2$-torsion line bundle, we can eliminate the obstruction.

The organization of this paper is as follows:

In \Cref{sec:pre}, we first review some concepts on parabolic bundles. Then we give several preparatory propositions.

In \Cref{sec:pu2-km-compactness}, we give the $PU(2)$ version of the Kronheimer--Mrowka's compactness theorem.

In \Cref{sec:detlift}, we discuss the topology of $PU(2)$ bundles and bundle pairs $(E,L)$. Then we lift the compactness of the $PU(2)$ case to  prescribed determinant $U(2)$ case.

In \Cref{sec:HE-to-KM}, we consider the parabolic Hermitian--Einstein connections. We first show $F_A\in L^2$, and then show $F_A\in L^p$ for some $p>2$ by techinques of Sibner--Sibner and Kronheimer--Mrowka. Therefore we construct a well-defined map $\mathfrak R:\mathcal M_h^{\mathrm{par,HE}}(E,L;\mathcal P)\to\mathcal M_\beta^{\mathrm{KM}}(E,L;a_{P,\mathrm{par}})$.

In \Cref{sec:KM-to-HE}, we consider the twisted ASD connections. Applying Biquard's extension theorem and constructing a compatible Hermitian--Einstein metric give the inverse map $\mathfrak S=\mathfrak R^{-1}$, hence we obtain the bijection.

In \Cref{sec:parabolic-polystability}, we prove the bijection between parabolic Hermitian--Einstein moduli space and parabolic polystable moduli space by the parabolic Kobayashi--Hitchin correspondence.

In \Cref{sec:ideal}, we construct the ideal parabolic Hermitian--Einstein pairs and topological bases inducing ideal convergence. By showing the compactness and metrizability of the topology, we get the desired compactification.

This is the first of our series on compactification of moduli spaces of stable parabolic bundles. In subsequent papers, we will consider more general cases, including those cases of normal crossings as well as higher dimensions. 
\newline

\textbf{Declaration of AI use} The main idea and the outline of the proof are due to the authors. The estimates of \Cref{thm:energy} and \Cref{prop:mixedlinear} are suggested by AI. The topological arguments in \Cref{lem:pairclassification} and \Cref{thm:ideal-topology} are also based on suggestions from AI. In the early version of the draft, AI pointed out the subtlety in \Cref{rem:projective-gauge-exception} and the examples in \Cref{rem:example-not-W2p} and \Cref{rem:example-not-canonical}.

All AI-generated suggestions are verified and revised by authors. The authors are fully responsible for all assertions in this paper.

\section{Preliminaries}\label{sec:pre}
We fix throughout the following data:
\begin{itemize}[leftmargin=2.4em]
\item $(\overline{X},J,\omega)$ is a connected compact complex surface equipped with a Gauduchon metric, i.e. $\partial\bar\partial\omega=0$;
\item $D\subset \overline{X}$ is a smooth connected divisor and $X:=\overline{X}\setminus D$;
\item $\mathcal{P}=(P,\bar\partial_{P})$ is a fixed holomorphic line bundle on $\overline{X}$;
\item $(E,h)$ is a smooth Hermitian bundle of rank $2$ on $\overline{X}$, satisfying $\det E\simeq P$;
\item $L\subset E|_D$ is a smooth complex line subbundle;
\item $0<\alpha_1<\alpha_2< 1$ are fixed weights, and $\beta:=\frac{1}{2}(\alpha_2-\alpha_1)\in(0,\frac{1}{2})$, $\bar\alpha:=\frac12(\alpha_1+\alpha_2)$.
\item Fix an exponent $p>2$ such that
\begin{align}\label{eq:p-assumption}
    2<p<\min\{p(\beta),2+\eta(\beta)\},\qquad\epsilon_p:=1-\frac2p<\min\{2\beta,1-2\beta\},
\end{align}
where $p(\beta)$ is defined in \cite[Proposition~2.4 (iii)]{KM93}, $\eta(\beta)$ is defined in \cite[Proposition~4.11]{KM93}.
\end{itemize}

A \textbf{divisor chart} is a local chart $U\subset\X$ such that $(U,U\cap D)$ diffeomorphic to $(B^2\times B^2,B^2\times\{0\})$. An \textbf{adapted chart} is a local chart $U\subset X$ or a divisor chart $U$.

\subsection{Parabolic bundles and compatible metrics}

\begin{definition}\label{def:gauduchondegree}
Let $\mathcal F$ be a torsion-free coherent sheaf on $\X$. Define the \textbf{degree} of $\mathcal{F}$
\begin{align*}
\deg_\omega(\mathcal F):=\frac{\ii}{2\pi}\int_{\overline{X}} F_h\wedge \omega,
\end{align*}
where $h$ is any smooth Hermitian metric on the determinant line bundle $\det\mathcal F:=(\wedge^{\rk\mathcal F}\mathcal F)^{**}$ and $F_h$ is the Chern curvature of $(\det\mathcal F,h)$. This is well-defined by the Gauduchon condition.
\end{definition}

\begin{definition}\label{def:parabolic}
A \textbf{rank-$2$ parabolic bundle} $\mathcal E_*=(\mathcal E,\mathcal L,\alpha_1,\alpha_2)$ is a holomorphic vector bundle $\mathcal E$ on $\overline{X}$ together with a holomorphic line subbundle $\mathcal L\subset \mathcal E|_D$
and the weights $0<\alpha_1<\alpha_2<1$.

Its \textbf{parabolic degree} is defined by
\begin{align*}
\pardeg_\omega(\mathcal E_*):=\deg_\omega(\mathcal E)+(\alpha_1+\alpha_2)\int_D \omega,
\end{align*}
and its \textbf{parabolic slope} is 
\begin{align*}
    \parmu_\omega(\mathcal E_*):=\frac12\pardeg_\omega(\mathcal E_*).
\end{align*}
\end{definition}

\begin{definition}\label{def:inducedweight}
Let $\mathcal S\subset\mathcal E$ be a saturated rank-one coherent subsheaf, i.e. the quotient sheaf $\mathcal{E}/\mathcal{S}$ is torsion-free. Define the \textbf{induced subsheaf weight}
$$
\epsilon(\mathcal S)=
\begin{cases}
\alpha_2,&\text{if $\mathcal S|_D=\mathcal L$ generically,}\\
\alpha_1,&\text{otherwise.}
\end{cases}
$$
Then
$$
\pardeg_\omega(\mathcal S_*)
 :=\deg_\omega\mathcal S+\epsilon(\mathcal S)\int_D\omega.
$$
\end{definition}

\begin{definition}\label{def:stability}
The parabolic bundle $\mathcal E_*$ is \textbf{stable}, respectively \textbf{semistable}, if
$$
\pardeg_\omega(\mathcal S_*)
 <\parmu_\omega(\mathcal E_*),
\qquad\text{respectively}\qquad
\pardeg_\omega(\mathcal S_*)
 \leq\parmu_\omega(\mathcal E_*),
$$
for every saturated rank-one coherent subsheaf $\mathcal S\subset\mathcal E$.
It is \textbf{polystable} if it is stable or a direct sum of stable parabolic line bundles of the same
parabolic slope.
\end{definition}

Choose a tubular neighbourhood $N(D)$ of $D$ and extend $L$ to a smooth line subbundle of $E|_{N(D)}$, still denoted by $L$. Let
\begin{align*}
E|_{N(D)}=L\oplus Q,\qquad Q:=L^{h\perp}.
\end{align*}
Let $\sigma$ denote the canonical section of $\mathcal O_{\overline{X}}(D)$, and choose a smooth Hermitian metric $h_D$ on $\mathcal O_{\X}(D)$. The model parabolic metric on $E|_X$ is
\begin{align}\label{eq:parmodel}
    h_{\pa}:=\|\sigma\|^{2\alpha_2} h|_L\oplus \|\sigma\|^{2\alpha_1} h|_Q
\end{align}
near $D$, extended smoothly over the rest of $X$.

\begin{definition}[Compatibility]\label{def:compatibility}
A smooth Hermitian metric $H$ on $\mathcal E|_X$ is \textbf{compatible with $\mathcal E_*$} if,
writing
$$b_H:=h_{\pa}^{-1}H,$$
one has
$$
 b_H,b_H^{-1}\in L^\infty(\End\mathcal E,h_{\pa}),\qquad
 \bar\partial_{\End\mathcal E}b_H\in L^2(h_{\pa},\omega),\qquad
 |\Lambda_\omega F_H|_H\in L^1(X).
$$
These are the conditions of \cite[Definition~2.11]{LiTian}.  
\end{definition}

\begin{lemma}[Independence of the parabolic model metric]
\label{lem:parmodelindependence}
Let $\mathcal E_*$ be a fixed rank-two parabolic bundle with line subbundle
$L\subset \mathcal E|_D$ and weights
$0<\alpha_1<\alpha_2<1$.  Let $h_{\pa}^{(0)}$ and $h_{\pa}^{(1)}$ be two
parabolic model metrics of the form \eqref{eq:parmodel}, constructed using
possibly different smooth Hermitian metrics on $\mathcal O_{\X}(D)$,
different smooth extensions of $L$, different smooth complementary line
subbundles, and different smooth positive Hermitian metrics on the
resulting splittings.  Set
\begin{align*}
 q:=(h_{\pa}^{(0)})^{-1}h_{\pa}^{(1)}
 \in \Gamma\bigl(X,\End\mathcal E\bigr).
\end{align*}
Then
\begin{align*}
 q,q^{-1}\in L^\infty(\End\mathcal E),
 \qquad
 \bar\partial_{\End\mathcal E}q,\,
 \bar\partial_{\End\mathcal E}q^{-1}\in L^2.
\end{align*}
Consequently, for a smooth Hermitian metric $H$ on $\mathcal E|_X$,
compatibility with respect to $h_{\pa}^{(0)}$ is equivalent to compatibility
with respect to $h_{\pa}^{(1)}$.  In particular, Definition
\ref{def:compatibility} is independent of the choice of the parabolic model metric.
\end{lemma}

\begin{proof}
It suffices to consider the behavior near $D$. Fix a divisor chart $U=\Delta_z\times\Delta_w$, $D\cap U=\{z=0\}$, and choose an adapted holomorphic frame $(e_1,e_2)$ of $\mathcal E$ on
$U$, with $e_1|_D$ spanning $L$.  Put $r=|z|$, $\gamma=\alpha_2-\alpha_1\in(0,1)$, and introduce the diagonal model
\begin{align*}
 G=\diag(r^{\alpha_2},r^{\alpha_1}),\qquad K=G^2.
\end{align*}

Consider one of the two model metrics, denoted by $h_{\pa}$.
Choose a smooth frame $(f_1,f_2)$ adapted to the smooth splitting used
in the construction of $h_{\pa}$ and orthonormal for its smooth factor
$h_L\oplus h_Q$.  If
\begin{align*}
 (f_1,f_2)=(e_1,e_2)C
\end{align*}
for a smooth map $C:U\to GL(2,\mathbb C)$, then the fact that both
$e_1|_D$ and $f_1|_D$ span $L$ gives $C_{21}|_D=0$.
Hence $C_{21}=O(r)$. The same property holds for $C^{-1}$.

In a local holomorphic trivialization of $\mathcal O_{\X}(D)$ in which $\sigma=z$, the norm of the defining section has the form $\|\sigma\|=r\rho$
with $\rho$ a smooth strictly positive function.  The corresponding smooth diagonal factors $\rho^{\alpha_2}$ and $\rho^{\alpha_1}$ do not
affect any of the estimates below.  Thus it is enough to consider the weighted conjugation $GCG^{-1}$.
Its off-diagonal entries satisfy
\begin{align*}
 (GCG^{-1})_{12}=O(r^\gamma),
 \qquad
 (GCG^{-1})_{21}=O(r^{1-\gamma}),
\end{align*}
while its diagonal entries remain smooth and bounded.  The same
statements hold for its inverse.  Differentiation therefore gives,
with respect to the Chern connection $D_K$ of the diagonal model,
\begin{align*}
 \bigl|D_K(GCG^{-1})\bigr|
 +
 \bigl|D_K(GC^{-1}G^{-1})\bigr|
 \leq
 C\bigl(1+r^{\gamma-1}+r^{-\gamma}\bigr).
\end{align*}
Since the normal volume element is uniformly comparable to
$r\,dr\,d\theta$, one has
\begin{align*}
 \int_0^\varepsilon r^{2\gamma-2}r\,dr<\infty,
 \qquad
 \int_0^\varepsilon r^{-2\gamma}r\,dr<\infty,
\end{align*}
because $0<\gamma<1$.  It follows that the matrix $u$ of $h_{\pa}$ relative to $K$, defined by $u:=K^{-1}h_{\pa}$, satisfies
\begin{align*}
 u,u^{-1}\in L^\infty,
 \qquad
 D_Ku\in L^2.
\end{align*}
Indeed, in the corresponding singular unitary frames $u$ is obtained
from the preceding weighted frame-change matrix and its adjoint by
matrix multiplication, and the assertion follows from the product
rule.  The same conclusion holds for $D_Ku^{-1}$.

Apply this conclusion to $h_{\pa}^{(0)}$ and $h_{\pa}^{(1)}$, and write
\begin{align*}
 u_a:=K^{-1}h_{\pa}^{(a)},\qquad a=0,1.
\end{align*}
Then $q=u_0^{-1}u_1$.
Hence $q$ and $q^{-1}$ are bounded, and the product rule gives
\begin{align*}
 D_Kq
 =
 -u_0^{-1}(D_Ku_0)u_0^{-1}u_1
 +u_0^{-1}D_Ku_1
 \in L^2.
\end{align*}
The same argument applied to $q^{-1}$ gives $D_Kq^{-1}\in L^2$.
Since the $(0,1)$ part of the Chern connection on
$\End\mathcal E$ is the fixed Dolbeault operator
$\bar\partial_{\End\mathcal E}$, this proves
\begin{align*}
 \bar\partial_{\End\mathcal E}q,\,
 \bar\partial_{\End\mathcal E}q^{-1}\in L^2.
\end{align*}
The boundedness of $q$ and $q^{-1}$ also shows that the norms induced
by $h_{\pa}^{(0)}$ and $h_{\pa}^{(1)}$ on $\End\mathcal E$ are uniformly
equivalent.

It remains to prove independence of compatibility.  Define
\begin{align*}
 b_a:=(h_{\pa}^{(a)})^{-1}H,\qquad a=0,1.
\end{align*}
Since $q=(h_{\pa}^{(0)})^{-1}h_{\pa}^{(1)}$, one has $b_1=q^{-1}b_0$.
Suppose first that $H$ is compatible with respect to $h_{\pa}^{(0)}$.
Then $b_0,b_0^{-1}\in L^\infty$, and therefore
$b_1,b_1^{-1}\in L^\infty$.  Moreover,
\begin{align*}
 \bar\partial_{\End\mathcal E}b_1
 =
 (\bar\partial_{\End\mathcal E}q^{-1})\,b_0
 +q^{-1}\bar\partial_{\End\mathcal E}b_0,
\end{align*}
and both terms on the right belong to $L^2$ by the estimates already
proved.  The remaining compatibility condition,
\begin{align*}
 |\Lambda_\omega F_H|_H\in L^1(X),
\end{align*}
depends only on $H$ and not on the choice of model metric.  Thus $H$
is compatible with respect to $h_{\pa}^{(1)}$.  Interchanging the two
models proves the converse.
\end{proof}

\subsection{Determinant normalization}\label{subsec:hunitary}
The gauge theory we will apply later works for the connections which are unitary for a smooth Hermitian metric, while the parabolic connections are unitary for the parabolic metrics, which are singular along the divisor $D$. Hence we need to use gauge transformation (which are singular along $D$) to connect these two settings. The following elementary construction will be used repeatedly.

\begin{definition}[$h$-unitary representative]\label{def:hunitaryrep}
Let $F\to U$ be a smooth complex vector bundle with a fixed smooth Hermitian metric $h$, and let $H$ be another Hermitian metric on $F|_U$. Assume $\nabla^H$ is an $H$-unitary connection on $F|_U$.

An $h$-unitary connection $\nabla$ is called an
\textbf{$h$-unitary representative} of $\nabla^H$ if
\begin{align*}
\nabla=(S_H)_*\nabla^H:=S_H\circ\nabla^H\circ S_H^{-1}
\end{align*}
for an isometry $S_H:(F,H)\to(F,h)$ satisfying $\det S_H=\det\widehat S_H$, where $\widehat S_H:(F,H)\to(F,h)$ is the unique positive
$h$-self-adjoint isometry.

The representative induced by $\widehat S_H$ is called the \textbf{canonical $h$-unitary representative}.
\end{definition}

\begin{proposition}\label{prop:hrepcurvature}

The integrability $F^{0,2}=0$ and the trace-free Hermitian--Einstein equation are preserved under $h$-unitary representative. Pointwise curvature norms agree when the two sides are computed with the corresponding metrics.
\end{proposition}
\begin{proof}
The operator $S_H$ is a bundle isometry $(F,H)\to(F,h)$.  Conjugation of a connection
by a bundle isomorphism conjugates its curvature.
Taking $(0,2)$ parts and applying $\ii\Lambda_\omega$ commute with conjugation on the
bundle factor.
\end{proof}

We modify the Hermitian metric $h$ such that the induced determinant metric $\det h$ satisfies the Hermitian--Einstein equation.

\begin{proposition}[Normalized determinant metric]\label{prop:detnormalization}
There is a smooth Hermitian metric $h_P$ on the fixed holomorphic determinant line bundle
$\mathcal{P}$ such that the singular metric
\begin{align}\label{eq:hppar}
    h_{P,\pa}:=\|\sigma\|^{2(\alpha_1+\alpha_2)}h_P\qquad\text{on }X
\end{align}
has the curvature form
\begin{equation}\label{eq:detcurvature}
F_{h_{P,\pa}}=\Big(F_{h_P}+(\alpha_1+\alpha_2)F_{h_D}\Big)\Big|_X,
\end{equation}
and satisfies
\begin{equation}\label{eq:detconstant}
\ii\Lambda_\omega F_{h_{P,\pa}}
 =\frac{2\pi}{\vol_\omega(\X)}
   \left(\deg_\omega \mathcal{P}+(\alpha_1+\alpha_2)\int_D\omega\right)=:2\lambda\in\RR.
\end{equation}
The $h_P$-unitary representative of the Chern connection of $(\mathcal{P}|_X,h_{P,\pa})$ is
\begin{equation}\label{eq:detconnection}
a_{P,\mathrm{par}}=a_{h_P}+(\alpha_1+\alpha_2)(\partial-\bar\partial)\log \|\sigma\|.
\end{equation}

\end{proposition}
\begin{proof}
By the standard computation, the conformal change $h_{P,\pa}=\|\sigma\|^{2(\alpha_1+\alpha_2)}h_P$ satisfies
\begin{align*}
    \nabla_{h_{P,\pa}}=\nabla_{h_P}+(\alpha_1+\alpha_2)\partial\log\|\sigma\|^2,\qquad F_{h_{P,\pa}}=F_{h_P}+(\alpha_1+\alpha_2)\bar\partial\partial\log\|\sigma\|^2.
\end{align*}

In a holomorphic frame $e_D$ of $\mathcal O(D)$ with $\sigma=ze_D$,
$$
\log \|\sigma\|^2=\log|z|^2+\log h_D(e_D,e_D).
$$
Since $\partial\bar\partial\log|z|^2=0$ on $X$, the Chern curvature of
$h_{P,\pa}$ on $X$ is exactly \eqref{eq:detcurvature}.  

The isometry $S:(P,h_{P,\pa})\to (P,h_P)$ in this case is obviously $S(v)=\|\sigma\|^{\alpha_1+\alpha_2}v$. Then \eqref{eq:detconnection} follows from
\begin{align*}
    a_{P,\pa}=S\nabla_{h_{P,\pa}}S^{-1}=\nabla_{h_{P,\pa}}-d\log\|\sigma\|^{\alpha_1+\alpha_2}=a_{h_P}+(\alpha_1+\alpha_2)(\partial-\bar\partial)\log\|\sigma\|.
\end{align*}

Start from an arbitrary smooth metric on $P$ and modify it by $e^fh_P$. The curvature changes as
\begin{align*}
    \ii\Lambda_\omega F_{e^fh_P}=\ii\Lambda_\omega F_{h_P}-\ii\Lambda_\omega\partial\bar\partial f.
\end{align*}
Thus it suffices to solve the elliptic equation
\begin{align*}
    \ii\Lambda_\omega\dd\db f=\ii\Lambda_\omega F_{h_{P,\pa}}-2\lambda=\ii\Lambda_\omega\Big(F_{h_P}+(\alpha_1+\alpha_2)F_{h_D}\Big)-2\lambda \qquad \text{on } X.
\end{align*}
The right hand side can be extended to $\X$ smoothly and has the zero average, so the elliptic theory gives a smooth solution $f\in C^\infty(\X)$ to the equation.
\end{proof}

From now on, we can assume the smooth Hermitian metric $h$ on $E$ satisfies
\begin{align}\label{eq:dethnormalized}
    \det h=h_P,\qquad \det h_{\pa}=h_{P,\pa}
\end{align}
under a suitable conformal change of \Cref{prop:detnormalization}.

\subsection{Hermitian--Einstein connections}

Let $P_E\to\X$ be the principal $\U(2)$ bundle of $h$-unitary frames, let
$P_{\det E}\to\X$ be the unitary frame bundle of $\det E$, and put
$$
 \widehat P_E:=P_E/\U(1)_\Delta,
 \qquad
 \U(1)_\Delta:=\{z\Id_{\CC^2}:z\in\U(1)\}.
$$
Thus $\widehat P_E$ is a principal $\PU(2):=\U(2)/\U(1)_\Delta$ bundle.  The homomorphism
\begin{equation}\label{eq:groupdetproj}
 \vartheta:\U(2)\longrightarrow \U(1)\times\PU(2),
 \qquad g\longmapsto(\det g,[g])
\end{equation}
is a surjective two-fold covering with kernel $\{\pm\Id\}$.  There is a corresponding
fiberwise two-fold covering
\begin{equation}\label{eq:bundledetprojcover}
 \Theta:P_E\longrightarrow P_{\det E}\times\widehat P_E,
 \qquad
 (e_1,e_2)\longmapsto(e_1\wedge e_2,[(e_1,e_2)]),
\end{equation}
which is equivariant with respect to \eqref{eq:groupdetproj}.

\begin{proposition}[Determinant--projective decomposition of a connection]
\label{prop:detproj}
The map \eqref{eq:bundledetprojcover} give a natural bijection
between
\begin{enumerate}
\item unitary connections $A$ on $E$, and
\item pairs $(a,B)$ consisting of a unitary connection $a$ on $\det E$ and a
$\PU(2)$ connection $B$ on $\widehat P_E$.
\end{enumerate}
We denote $a=\det A$ and $B=A^{\rm proj}$ under this bijection.  
Moreover, the curvatures satisfy 
\begin{equation}\label{eq:curvdetproj}
 F_{\det A}=\tr F_A,
 \qquad
 F_{A^{\rm proj}}=
 F_A^\perp:=F_A-\frac12(\tr F_A)\Id,
\end{equation}
under the Lie-algebra identifications $\mathfrak u(2)\simeq\mathfrak{u}(1)\oplus\mathfrak{pu}(2)$ induced by \eqref{eq:groupdetproj}.

\end{proposition}

\begin{proof}
  Since the homomorphism $\vartheta$ is a two-fold covering, its differential at the identity is an isomorphism of Lie algebras
  \begin{align*}
      d\vartheta:\mathfrak u(2)\longrightarrow
\mathfrak u(1)\oplus\mathfrak{pu}(2),
\qquad
\xi\longmapsto(\tr\xi,[\xi]).
  \end{align*}
  The Lie algebra $\mathfrak{pu}(2)=\mathfrak u(2)/\sqrt{-1}\mathbb R\cdot\Id$
is canonically identified with $\mathfrak{su}(2)$ by
\begin{align}\label{eq:pu2su2}
    [\xi]\longmapsto\xi^\perp:=\xi-\frac12(\tr\xi)\Id.
\end{align}

Under this identification,
\begin{align*}
d\vartheta(\xi)
=
\bigl(\tr\xi,\xi^\perp\bigr),
\end{align*}
and its inverse is
\begin{align*}
(a,b)\longmapsto b+\frac12a\cdot\Id ,
\qquad
a\in\mathfrak u(1),\quad b\in\mathfrak{su}(2).
\end{align*}

We next prove that $\Theta$ induces a bijection between connections on $P_E$ and connections on $P_{\det E}\times\widehat P_E$. Let $A$ be a unitary connection on $E$, and let $\omega_A\in\Omega^1(P_E;\mathfrak u(2))$ be its principal connection form.  Define
\begin{align*}
    \widetilde\omega_A:=d\vartheta\circ\omega_A
\in
\Omega^1(P_E;\mathfrak u(1)\oplus\mathfrak{pu}(2)).
\end{align*}

For $k\in\{\pm\Id\}$, the connection-equivariance identity gives
\begin{align*}
    R_k^*\omega_A=\Ad_{k^{-1}}\omega_A=\omega_A,
\end{align*}
since $\{\pm\Id\}$ is central in $\U(2)$.  Hence $R_k^*\widetilde\omega_A=\widetilde\omega_A.$

Since $\Theta:P_E\to P_{\det E}\times\widehat P_E$ is the quotient by the free action of the finite group $\{\pm\Id\}$, there is therefore a unique one-form
\begin{align*}
\omega_Q\in
\Omega^1(Q;\mathfrak u(1)\oplus\mathfrak{pu}(2))
\end{align*}
such that
\begin{align}\label{eq:omegaQA}
    \Theta^*\omega_Q=\widetilde\omega_A=d\vartheta\circ\omega_A.
\end{align}

We show that $\omega_Q$ is a principal connection form. For any $\eta\in\mathfrak u(1)\oplus\mathfrak{pu}(2)$, put $\xi=(d\vartheta)^{-1}\eta\in\mathfrak u(2)$. By $\vartheta$-equivariance of $\Theta$, the corresponding fundamental vector fields satisfy $d\Theta(\xi^\#)=\eta^\#$. Consequently, we have
\begin{align*}
    \omega_Q(\eta^\#)
    =\omega_Q\circ d\Theta(\xi^\#)
    =\Theta^*\omega_Q(\xi^\#)
    =d\vartheta\circ\omega_A(\xi^\#)
    =d\vartheta(\xi)
    =\eta.
\end{align*}

For equivariance, let $h\in U(1)\times PU(2)$.  Choose $g\in U(2)$ with $\vartheta(g)=h$, since $\vartheta$ is surjective. Then the equivariance implies $\Theta\circ R_g=R_h\circ\Theta$.

We have
\begin{align*}
    \Theta^*\bigl(R_h^*\omega_Q\bigr)
=R_g^*\bigl(\Theta^*\omega_Q\bigr)
=R_g^*\bigl( d\vartheta\circ\omega_A\bigr)
=d\vartheta\bigl(\Ad_{g^{-1}}\omega_A\bigr)
=\Ad_{h^{-1}}
(d\vartheta\circ\omega_A)=\Theta^*\bigl(\Ad_{h^{-1}}\omega_Q\bigr).
\end{align*}
Since $\Theta$ is a covering map, pullback by $\Theta$ is injective on differential forms.  Hence
\begin{align*}
R_h^*\omega_Q=\Ad_{h^{-1}}\omega_Q.
\end{align*}
Thus $\omega_Q$ is a principal $U(1)\times PU(2)$-connection on $P_{\det E}\times\widehat P_E$.

Conversely, if $\omega_Q$ is a principal $U(1)\times PU(2)$-connection on $P_{\det E}\times\widehat P_E$, the same calculation shows that
\begin{align*}
    \omega_A:=(d\vartheta)^{-1}\Theta^*\omega_Q
\end{align*}
is a principal $\U(2)$-connection on $P_E$. Therefore connections on $P_E$ are naturally in bijection with connections on $P_{\det E}\times\widehat P_E$.

It remains to identify a connection on $P_{\det E}\times\widehat P_E$ with a pair of connections on its two factors.  Let
\begin{align*}
    \pi_1:P_{\det E}\times\widehat P_E\to P_{\det E},\qquad \pi_2:P_{\det E}\times\widehat P_E\to \widehat P_E
\end{align*}
be the projections.  Given connection forms
\begin{align*}
    \omega_a\in\Omega^1(P_{\det E};\mathfrak u(1)),
\qquad
\omega_B\in\Omega^1(\widehat P_E;\mathfrak{pu}(2)),
\end{align*}
the form
\begin{align}\label{eq:QaB}
    \omega_Q:=\bigl(\pi_1^*\omega_a,\pi_2^*\omega_B\bigr)
\end{align}
is a connection form on the principal
$\U(1)\times\PU(2)$ bundle $P_{\det E}\times\widehat P_E$.

Conversely, write an arbitrary connection form on $Q$ as $\omega_Q=(\omega_1,\omega_2)$. The first component $\omega_1$ vanishes on vertical vectors generated by
$\mathfrak{pu}(2)$, because the reproduction property gives
\begin{align*}
\omega_Q(0,\eta)^\#=(0,\eta).
\end{align*}
Moreover, $\omega_1$ is invariant under the $\PU(2)$ action, since the adjoint action of the second factor of
$\U(1)\times\PU(2)$ acts trivially on the first Lie-algebra factor.  Thus $\omega_1$ is basic for the projection $\pi_1$ and descends uniquely to a connection form $\omega_a$ on $P_{\det E}$.  The same argument shows that $\omega_2$ is basic for $\pi_2$ and descends uniquely to a connection form $\omega_B$ on $\widehat P_E$.  Hence \eqref{eq:QaB} gives a bijection between connections on $P_{\det E}\times\widehat P_E$ and pairs $(a,B)$.

Combining the preceding two bijections gives the asserted natural correspondence
\begin{align*}
    A\longleftrightarrow(\det A,A^{\rm proj}).
\end{align*}

We finally compute the curvatures.  Let $F_A
=
d\omega_A+\frac12[\omega_A\wedge\omega_A]$ be the curvature form of $A$.  Since $d\vartheta$ is a Lie-algebra homomorphism, \eqref{eq:omegaQA} implies $$\Theta^*F_Q=d\vartheta(F_A).$$

Under $d\vartheta(\xi)=(\tr\xi,[\xi])$, the two components of $\Theta^*F_Q$ are therefore
\begin{align*}
    F_{\det A}=\tr F_A,
\qquad
F_{A^{\rm proj}}=[F_A].
\end{align*}

Under the canonical identification
$\mathfrak{pu}(2)\simeq\mathfrak{su}(2)$ in \eqref{eq:pu2su2}, the class $[F_A]$ is represented by its unique trace-free representative $F_A^\perp$. Hence, we have
\begin{align*}
    F_{A^{\rm proj}}=F_A^\perp:=F_A-\frac12(\tr F_A)\Id.
\end{align*}
\end{proof}

\begin{remark}\label{rem:detprojgauge}
The determinant-projective decomposition is a bijection for connections. However, the
analogous statement for gauge transformations is different.  A unitary automorphism
$u$ induces the pair $(\det u,[u])$.  If $\det u=1$, its projectivization is a
projective gauge transformation, but a projective gauge transformation need not lift
to a determinant-one gauge.  The obstruction is the pullback of the central double
cover $\SU(2)\to\PU(2)$ and will be kept explicitly in \Cref{sec:detlift}.
\end{remark}

\begin{lemma}\label{lem:scalarHE}
Let $A$ be a unitary connection on $E|_{X}$ satisfying $\det A=a_{P,\mathrm{par}}$. Then the following statements of Hermitian--Einstein condition are equivalent:
\begin{enumerate}[label=\textup{(\roman*)}]
\item $F_A^{0,2}=0$ and there exists a constant $\mu\in\mathbb R$ such that $\sqrt{-1}\Lambda_\omega F_A=\mu\,\Id$;
\item $F_A^{0,2}=0$ and $\sqrt{-1}\Lambda_\omega F_A^\perp=0$;
\item $(F_A^\perp)^+=0$, i.e. $F_A^\perp$ satisfies ASD equation.
\end{enumerate}
Moreover, in {\rm(i)} one necessarily has $\mu=\lambda$, where $\lambda$ is defined by
\eqref{eq:detconstant}.
\end{lemma}

\begin{proof}
Since
\begin{align*}
F_A=F_A^\perp+\frac12F_{\det A}\,\Id
\end{align*}
and $\det A=a_{P,\mathrm{par}}$, \eqref{eq:detconstant} gives
\begin{align*}
\sqrt{-1}\Lambda_\omega F_{\det A}=2\lambda.
\end{align*}
Hence
\begin{align*}
\sqrt{-1}\Lambda_\omega F_A
=
\sqrt{-1}\Lambda_\omega F_A^\perp+\lambda\,\Id.
\end{align*}
Therefore (i) and (ii) are equivalent, and whenever they hold one has
$\mu=\lambda$.

We next compare (ii) and (iii).  Since $A$ is unitary,
\begin{align*}
(F_A^\perp)^{2,0}
=
-\bigl((F_A^\perp)^{0,2}\bigr)^*.
\end{align*}
Thus $F_A^{0,2}=0$ is equivalent to
\begin{align*}
(F_A^\perp)^{2,0}=(F_A^\perp)^{0,2}=0,
\end{align*}
because $F_{\det A}$ is of type $(1,1)$.

On a Hermitian surface, the real two-forms decompose orthogonally as
\begin{align*}
\Lambda^2_{\mathbb R}=
\mathbb R\omega\oplus\bigl(\Lambda^{2,0}\oplus\Lambda^{0,2}\bigr)_{\mathbb R}\oplus\Lambda^{1,1}_0,
\end{align*}
where
\begin{align*}
\Lambda^+=\mathbb R\omega\oplus
\bigl(\Lambda^{2,0}\oplus\Lambda^{0,2}\bigr)_{\mathbb R},
\qquad
\Lambda^-=\Lambda^{1,1}_0.
\end{align*}
Consequently, a unitary trace-free curvature form is anti-self-dual if and only if it is of type
$(1,1)$ and primitive.  Hence
\begin{align*}
(F_A^\perp)^+=0
\end{align*}
is equivalent to
\begin{align*}
(F_A^\perp)^{0,2}=0,
\qquad
\Lambda_\omega F_A^\perp=0.
\end{align*}
Since $F_{\det A}^{0,2}=0$, this is precisely condition {\rm(ii)}.
\end{proof}

\begin{lemma}\label{lem:scalarremoval}
Let $f\in C^{\infty}(X)\cap L^{\infty}(X)$ with $df\in L^2(X)$. Suppose
$$\ii\Lambda_\omega\partial\bar\partial f=\lambda.$$
Then $\lambda=0$ and $f$ is constant.
\end{lemma} 
\begin{proof}
    Let $r$ denote the radial function of the tubular neighborhood $N(D)$. Choose $\eta\in C^\infty(\RR,[0,1])$ with $\eta|_{(-\infty,0]}\equiv 0$, $\eta|_{[1,+\infty)}\equiv 1$, and $|\eta'|\le C$. For $0<\varepsilon<1$, we construct the cut-off function
\begin{align*}
\chi_\varepsilon(r):=\eta\left(\frac{\log(r/\varepsilon^2)}{|\log\varepsilon|}\right).
\end{align*}
Then $\chi_\varepsilon\in C^\infty(X)$ satisfies
\begin{align}\label{eq:logcutoff}
    \chi_\varepsilon|_{r\le \varepsilon^2}\equiv 0,\quad \chi_\varepsilon|_{r\ge \varepsilon}\equiv 1,\quad |d\chi_\varepsilon|\le\frac{C}{r|\log\varepsilon|}. 
\end{align}
Consequently,
\begin{align*}
    \|d\chi_\varepsilon\|_{L^2(X)}^2\leq
C\int_{\varepsilon^2}^{\varepsilon}\frac{rdr}{r^2|\log\varepsilon|^2}
\leq
\frac{C}{|\log\varepsilon|}=o(1).
\end{align*}

Let $\varphi\in C^\infty(\X,\mathbb R)$. The equation and integration by parts gives
\begin{align*}
\lambda\int_X\chi_\varepsilon\varphi\,dV_\omega&=\ii\int_X\chi_\varepsilon\varphi\,\Lambda_\omega\dd\db f\,dV_\omega\\
&=\ii\int_X
\chi_\varepsilon\varphi\,
\partial\bar\partial f\wedge\omega\\
&=\ii\int_X\db f\wedge\dd(\chi_\varepsilon\varphi\omega)\\
&=\ii\int_X\bar\partial f\wedge\chi_\varepsilon\partial(\varphi\omega)+\ii\int_X\varphi\,\bar\partial f\wedge\partial\chi_\varepsilon\wedge\omega\\
&=-\ii\int_X f\db\chi_\varepsilon\wedge\partial(\varphi\omega)+\ii\int_X\chi_\varepsilon f\wedge\dd\db(\varphi\omega)+\ii\int_X\varphi\,\bar\partial f\wedge\partial\chi_\varepsilon\wedge\omega.
\end{align*}

By Cauchy--Schwarz, as $\varepsilon\to 0$,
\begin{align*}
\left|\int_X f\db\chi_\varepsilon\wedge\dd(\varphi\omega)\right|
&\le
C\|\varphi\|_{C^1}\|f\|_{L^2}\|d\chi_\varepsilon\|_{L^2}\to0,\\
    \left|\int_X \varphi\,\db f\wedge\partial\chi_\varepsilon\wedge\omega\right|
&\le
C\|\varphi\|_{L^\infty}\|df\|_{L^2}\|d\chi_\varepsilon\|_{L^2}\to0.
\end{align*}

Hence, letting $\varepsilon\to0$, we obtain
\begin{align}\label{eq:ellipticdistribution}
    \lambda\int_\X\varphi\,dV_\omega=\int_X\ii\,f\dd\db(\varphi\omega),
\end{align}
for every $\varphi\in C^\infty(\X)$. Thus
\begin{align*}
\sqrt{-1}\Lambda_\omega\partial\bar\partial f=\lambda
\end{align*}
holds on $\X$ in the distribution sense.

Taking $\varphi\equiv1$ in \eqref{eq:ellipticdistribution} gives
\begin{align*}
\lambda\,\operatorname{Vol}_\omega(\X)=\ii\int_X f\dd\db\omega=0,
\end{align*}
hence $\lambda=0$.

Now we have
\begin{align*}
    \ii\Lambda_\omega\dd\db f=0
\end{align*}
distributionally on the compact manifold $\X$. Since $f\in L^\infty(X)\subset L^2(\X)$, standard interior elliptic regularity gives $f\in C^\infty(\X)$. Finally, the strong maximum principle implies that $f$ is constant.
\end{proof}

\begin{proposition}[Determinant normalization of a compatible HE metric]\label{prop:metricdet}
Let $H$ be a Hermitian--Einstein metric on $E|_X$ compatible with the parabolic structure $\mathcal{E}_*$ and $h_{\pa}$ be the singular model metric normalized in \eqref{eq:dethnormalized}.  

Then $\det H/\det h_{\pa}$ is a positive constant. So we can always assume $\det H=\det h_{\pa}$ without changing its Chern connection or compatibility.
\end{proposition}
\begin{proof}
$H$ is a compatible Hermitian--Einstein metric means
\begin{align*}
    \ii\Lambda_\omega F_H=\lambda_H \Id,\quad b_H,b_H^{-1}\in L^{\infty}(h_{\pa}),\quad \db_{\End \E}b_H\in L^2(h_{\pa}),
\end{align*}
where $b_H:=h_{\pa}^{-1}H\in\Gamma(X,\End\E)$.

Take
$$
f:=\log\frac{\det H}{\det h_{\pa}}=\log\det b_H\in L^{\infty}(X,\RR).
$$
We have
\begin{align*}
    \db f=\db\log\det b_H=\tr\big(b_H^{-1}\db b_H\big)\in L^2.
\end{align*}

 Taking determinant in the Hermitian--Einstein equation and normalization \eqref{eq:dethnormalized} gives
$$
\ii\Lambda_\omega F_{\det H}=2\lambda_H,\qquad \ii\Lambda_\omega F_{\det h_{\pa}}=2\lambda.
$$
Thus
\begin{align*}
    \ii\Lambda_\omega\dd\db f=2(\lambda-\lambda_H).
\end{align*}
\cref{lem:scalarremoval} gives
$\lambda_H=\lambda$ and $f$ is constant.
\end{proof}

\section{The \texorpdfstring{$PU(2)$}{PU(2)} version of Kronheimer--Mrowka compactness}
\label{sec:weighted}\label{sec:pu2-km-compactness}
 
 Kronheimer--Mrowka develop a singular version of Uhlenbeck's compactness theory for $\SU(2)$ instanton \cite{KM93}, and they discuss the $SO(3)$ case very briefly \cite[\S~2(iv)]{KM93}. The local arguments are similar, however the global topology is different due to the Stiefel--Whitney class $w_2$. The purpose of this section is to formulate the compactness theorem for singular
 $PU(2)$ instanton in a form that can be applied to the projective unitary frame bundle. We use the fixed identification
$PU(2)\cong SO(3)$ and work with principal $SO(3)$ bundles.  

This section doesn't involve the complex structure. Throughout the section, $(\X,g)$ is a closed oriented smooth Riemannian
four-manifold, $D\subset\X$ is a closed connected oriented smoothly embedded
surface, and let $X:=\X\setminus D$. Fix a closed tubular neighborhood $N$ of $D$, identified with the unit disk
bundle of the normal bundle of $D$.  Let $r$ be the radial coordinate,
let $Y=\partial N$, and let $\eta\in\Omega^1(N\setminus D)$ be a connection one-form on the circle bundle $Y\to D$, extended to $N\setminus D$ by radial projection and
normalized by
\begin{align*}
        \int_{S^1_{\mathrm{fiber}}}\eta=2\pi.
\end{align*}
We also fix a smooth cut-off function $\chi(r):N\to[0,1]$ that is identically
one near $D$ and is supported in the interior of $N$.

\subsection{Projective bundles and singular connections}
Let
\begin{align*}
    U(1)_\Delta:=\{z\Id_{\mathbb C^2}:z\in U(1)\}\subset U(2),\qquad PU(2):=U(2)/U(1)_\Delta\cong SU(2)/\{\pm1\}.
\end{align*}
 Conjugation on $\mathfrak{su}(2)$ gives the fixed isomorphism
\begin{equation}\label{eqn:pu2=so3}
    \begin{aligned}
        \rho:PU(2)&\xrightarrow{\;\cong\;}
        SO\bigl(\mathfrak{su}(2),-\tr_{\mathbb C^2}\bigr)\cong SO(3),\\
        [A]&\longmapsto(X\mapsto AXA^{-1}),
    \end{aligned}
\end{equation}
where the latter isomorphism is induced by the orthonormal basis of $\big(\mathfrak{su}(2),-\tr_{\mathbb C^2}\big)$:
\begin{align}\label{eqn:su2-basis}
    e_1:=\frac{1}{\sqrt{2}}\begin{pmatrix}
        \ii&0\\
        0&-\ii
        \end{pmatrix},\
    e_2:=\frac{1}{\sqrt{2}}\begin{pmatrix}
        0&1\\
        -1&0
        \end{pmatrix},\
    e_3:=\frac{1}{\sqrt{2}}\begin{pmatrix}
        0&\ii\\
        \ii&0
        \end{pmatrix}.
\end{align}

On $\mathfrak{so}(3)$ we use the invariant trace $\tr_3$ characterized by
\begin{align*}
        \tr_3\bigl(d\rho(\xi)d\rho(\zeta)\bigr)
        =\tr_{\mathbb C^2}(\xi\zeta),
        \qquad \xi,\zeta\in\mathfrak{pu}(2)\cong\mathfrak{su}(2),
\end{align*}
and the inner product
\begin{align*}
        \langle a,b\rangle_{\mathfrak{so}(3)}:=-\tr_3(ab).
\end{align*}
Equivalently, in the standard three-dimensional real representation,
$\tr_3=\frac14\tr_{\mathbb R^3}$.  

\begin{lemma}
\label{lem:diagonal-torus-under-fixed-identification}
Let $T_P\subset PU(2)$ be the diagonal torus defined by
\begin{align*}
    T_P:=\left\{ \left[ \begin{pmatrix} z_1 & 0 \\ 0 & z_2 \end{pmatrix} \right] \in PU(2):z_1,z_2\in U(1) \right\}.
\end{align*}
Under the identification (\ref{eqn:pu2=so3}), $T_P$ is equivalent to the standard subgroup $T_S\cong SO(2)\subset SO(3)$:
\begin{align}\label{eqn:Ts}
    T_S:=\left\{\begin{pmatrix} 1 & 0 & 0\\ 0 & \cos\theta &-\sin\theta\\ 0 &\sin\theta & \cos\theta\end{pmatrix}:\theta\in\mathbb R\right\}\subset SO(3).
\end{align}
Moreover, for $e_1\in\mathfrak{su}(2)\cong\mathfrak{pu}(2)$ defined in (\ref{eqn:su2-basis}), we have
\begin{align*}
    d\rho(\sqrt{2}e_1):=d\rho\begin{pmatrix}
        \ii & \\ & -\ii    
    \end{pmatrix}
    =\begin{pmatrix}
        0 & 0 &0\\0 & 0 & -2\\0 & 2 & 0
    \end{pmatrix}=:2J_0\in\mathfrak{so}(3).
\end{align*}
\end{lemma}

\begin{proof}
We first regard $PU(2)$ as $SU(2)/\{\pm1\}$, where
\begin{align*}
    SU(2)=\left\{\begin{pmatrix}a&b\\-\bar{b}&\bar{a}\end{pmatrix}:|a|^2+|b|^2=1\right\}.
\end{align*}
A direct calculation gives
\begin{align*}
    \rho\left[\begin{pmatrix}a&b\\-\bar{b}&\bar{a}\end{pmatrix}\right]=\begin{pmatrix}
        |a|^2-|b|^2 & 2\im(a\bar{b}) & 2\re(a\bar{b})\\
        2\im(ab) & \re(a^2+b^2) & -\im(a^2-b^2)\\
        -2\re(ab) & \im(a^2+b^2) & \re(a^2-b^2)
    \end{pmatrix}.
\end{align*}
Thus for the diagonal torus $T_P\subset PU(2)$,
\begin{align*}
    \rho\left[\begin{pmatrix}z_1&\\&z_2\end{pmatrix}\right]=\rho\left[\begin{pmatrix}e^{\ii\theta/2}&\\&e^{-\ii\theta/2}\end{pmatrix}\right]=\begin{pmatrix}
        1 & 0 & 0\\ 0 & \cos\theta &-\sin\theta\\ 0 &\sin\theta & \cos\theta
    \end{pmatrix},
\end{align*}
where $z_1z_2^{-1}=e^{\ii\theta}$.

Finally,
\begin{align*}
    d\rho\begin{pmatrix}
        \ii & \\ & -\ii    
    \end{pmatrix}
    =\left.\frac{d}{dt}\right|_{t=0}\rho\left[\begin{pmatrix}e^{\ii t} &\\& e^{-\ii t}\end{pmatrix}\right]
    =\left.\frac{d}{dt}\right|_{t=0}\begin{pmatrix}
        1 & 0 & 0\\ 0 & \cos 2t &-\sin 2t\\ 0 &\sin 2t & \cos 2t
    \end{pmatrix}=2J_0.
\end{align*}
\end{proof}
\begin{remark}\label{rmk:orientedSO(2)}
    From now on, we fix an oriented subgroup $(T_S,J_0)\subset SO(3)$, which means a subgroup together with a fixed positive generator $J_0$ rather than $-J_0$. This orientation distinguishes the line subbundle $L$ in the splitting $E|_N=L\oplus Q$.
\end{remark}

\begin{definition}[Reduction of structure group]
\label{def:principal-reduction}
Let $H\subset G$ be a closed subgroup and let $P\to B$ be a principal
$G$ bundle.  A \textbf{reduction} of $P$ to $H$ is a principal $H$ subbundle
$R\subset P$ such that the canonical map
\begin{align*}
    R\times_H G\longrightarrow P,
        \qquad [(q,g)]\longmapsto qg,
\end{align*}
is an isomorphism.  Equivalently, it is a smooth section of the associated
bundle $P/H\to B$.
\end{definition}

\begin{definition}[Surface reduction]
\label{def:surface-reduction}
Let $\widehat P\to\X$ be a principal $SO(3)$ bundle.  A \textbf{surface reduction} of $\widehat P$ along
$D$ is a reduction
\begin{align*}
        R\subset \widehat P|_N
\end{align*}
to the subgroup $T_S\cong SO(2)\subset SO(3)$ defined in (\ref{eqn:Ts}). A \textbf{projective bundle pair} $(\widehat P,R)$ is an $SO(3)$ bundle $\widehat P$ together with a surface reduction $R$.
\end{definition}

Let
\begin{align*}
        \mathfrak g_{\widehat P}:=\widehat P\times_{\operatorname{Ad}}\mathfrak{so}(3)
\end{align*}
be the adjoint bundle of $\widehat P$.  The fixed subgroup $T_S\cong SO(2)\subset SO(3)$ naturally induces an orthogonal
splitting
\begin{align*}
        \mathfrak{so}(3)=\mathbb RJ_0\oplus(\mathbb RJ_0)^\perp.
\end{align*}
Next we show the surface reduction induces a splitting of the adjoint bundle.

\begin{proposition}
\label{prop:surface-reduction-adjoint-splitting}
A surface reduction $R\subset \widehat P|_N$ determines a canonical orthogonal
splitting
\begin{align*}
        \mathfrak g_{\widehat P}|_N=\lambda_R\oplus K_R,
\end{align*}
where
\begin{align*}
        \lambda_R:=R\times_{SO(2)}\mathbb RJ_0,
        \qquad
        K_R:=R\times_{SO(2)}(\mathbb RJ_0)^\perp.
\end{align*}
Moreover, $\lambda_R$ is a trivial line bundle with a nowhere vanishing global section $J_R$, and $K_R$ is a complex line bundle.
\end{proposition}

\begin{proof}
Since the adjoint action of $SO(2)$ fixes $J$ and preserves
$(\mathbb RJ)^\perp$, the associated-bundle construction gives the
splitting.

We define a section of $\lambda_R$ by
\begin{align*}
        J_R(x):=[(q,J_0)],\qquad q\in R_x.
\end{align*}
For any $x\in N$, $q_1,q_2\in R_x$, since
\begin{align*}
    [(q_1,J_0)]=[(q_2\cdot q_2^{-1}q_1,J_0)]=[(q_2,\operatorname{Ad}_{q_2^{-1}q_1}J_0)]=[(q_2,J_0)],
\end{align*}
$J_R$ is well-defined and trivializes $\lambda_R$.  Moreover, we claim $K_R$ has a complex structure to be a complex line bundle:
\begin{align*}
        \mathcal J_R:=\ad(J_R)|_{K_R},\quad [(q,\xi)]\mapsto [(q,[J_0,\xi])].
\end{align*}
Take
\begin{align*}
    A:=\begin{pmatrix}
        0 & 0 & 1\\
        0 & 0 & 0\\
        -1& 0 & 0
    \end{pmatrix},\qquad
    B:=\begin{pmatrix}
        0 & -1 & 0\\
        1 & 0 & 0\\
        0 & 0 & 0
    \end{pmatrix},
\end{align*}
then $(\RR J_0)^\perp=\mathrm{span}_{\RR}\{A,B\}$. Direct computation gives $[J_0,A]=B$, $[J_0,B]=-A$. Hence, $\mathcal{J}_R:K_R\to K_R$ and $\mathcal{J}_R^2=-\Id_{K_R}$.
\end{proof}

\begin{proposition}[Reduction induced by an ordered line]\label{prop:line-reduction}
Let $E$ be a Hermitian rank-two bundle and let
$E|_N=L\oplus Q$ be an ordered orthogonal splitting.  The adapted projective unitary frames give
a reduction $R_L$ of $\widehat P_E=\PU(E)$ along $D$. Under
$\mathfrak g_{\widehat P_E}\cong\mathfrak{su}(E)$ one has
\begin{equation}\label{eq:Kline}
 K_{R_L}\cong\Hom(Q,L)
\end{equation}
as complex line bundles.  Conversely, for a fixed unitary lift $E$, the surface reduction of $\widehat P_E$ determines the line subbundle $L$.
\end{proposition}

\begin{proof}
Denote by $P_E$ the principal $\U(2)$-bundle of unitary frames of $E$. A unitary frame $(e_1,e_2)$ is said to be adapted to the splitting if $e_1\in L$, $e_2\in Q$. Then the set of adapted unitary frames forms a principal $\U(1)\times\U(1)$ subbundle $P_{L,Q}\subset P_E|_N$. Since the scalar subgroup $U(1)_\Delta$ is contained in $\U(1)\times\U(1)$, passing to the projective quotient gives a principal subbundle
\begin{align*}
R_L:=P_{L,Q}/\U(1)_\Delta
\subset
P_E/\U(1)_\Delta
=
\widehat P_E.
\end{align*}
Its structure group is
\begin{align*}
(\U(1)\times\U(1))/\U(1)_\Delta
\cong \U(1)\cong SO(2).
\end{align*}
Thus $R_L$ is a surface reduction of $\widehat P_E$.

We now identify the associated rank-two real bundle $K_{R_L}$. Fix $x\in N$ and an adapted unitary frame $(e_1,e_2)$ of $E_x=L_x\oplus Q_x$.  With respect to this frame, every element of $\mathfrak{su}(E_x)$ has the form
\begin{align*}
\xi=
\begin{pmatrix}
\sqrt{-1}t & \nu\\
-\bar{\nu} & -\sqrt{-1}t
\end{pmatrix},
\qquad
t\in\RR, \nu\in\CC.
\end{align*}
The diagonal part is precisely the fiber of $\lambda_{R_L}$, while the off-diagonal part is the fiber of $K_{R_L}$.

The upper-right entry $\nu$ represents a complex-linear map $Q_x\to L_x$. Thus, fiberwise, we have an isomorphism
\begin{align}\label{eq:KR-Hom(Q,L)}
    K_{R_L,x}\longrightarrow \Hom(Q_x,L_x),
\qquad
\begin{pmatrix}
0&\nu\\
-\overline{\nu}&0
\end{pmatrix}
\longmapsto
\bigl(e_2\mapsto \nu e_1\bigr).
\end{align}
It's easy to verify that this map is independent of the chosen adapted unitary frame. This map gives an isomorphism of real rank-two vector bundles $K_{R_L}\to \Hom(Q,L)$.

By definition, the complex structure on $K_{R_L}$ is $\mathcal J_{R_L}=\ad(J_{R_L})|_{K_{R_L}}$. Under the identification $d\rho(T)=J_0$ by \Cref{lem:diagonal-torus-under-fixed-identification}, where $T:=\frac{1}{\sqrt{2}}e_1$, this corresponds to the commutator action $\ad(T)$ on the off-diagonal part of $\mathfrak{su}(2)$.

For
\begin{align*}
\xi_\nu:=
\begin{pmatrix}
0&\nu\\
-\overline{\nu}&0
\end{pmatrix},
\end{align*}
a direct computation gives
\begin{align*}
\ad(T)(\xi_\nu)=\xi_{\sqrt{-1}\nu}.
\end{align*}
Under the isomorphism \eqref{eq:KR-Hom(Q,L)}, this is exactly multiplication by $\sqrt{-1}$ on $\Hom(Q_x,L_x)$. Hence \eqref{eq:KR-Hom(Q,L)} is an isomorphism of complex line bundles.

We finally prove the converse.  Assume that the unitary rank-two bundle $E$ is fixed and that surface reduction $R\subset\widehat P_E|_N$ is given.

Let
\begin{align*}
\pi:P_E|_N\longrightarrow\widehat P_E|_N
\end{align*}
be the quotient map.  The inverse image
\begin{align*}
\widetilde R:=\pi^{-1}(R)
\end{align*}
is a reduction of $P_E|_N$ to the inverse image of the $\SO(2)$-subgroup under $\U(2)\to\PU(2)$, which is exactly is the diagonal torus $T^2:=\U(1)\times\U(1)\subset\U(2)$ by \Cref{lem:diagonal-torus-under-fixed-identification}.

Therefore $\widetilde R$ is a $T^2$-reduction of the unitary frame bundle $P_E|_N$. Let $\chi_i:T^2\to U(1)_i$ be the projection to each factor for $i=1,2$, then
\begin{align*}
L_R:=\widetilde R\times_{\chi_1}\mathbb C,
\qquad
Q_R:=\widetilde R\times_{\chi_2}\mathbb C,
\end{align*}
and the standard representation of $T^2$ on $\mathbb C^2$ gives an orthogonal splitting into two complex line bundles
\begin{align*}
E|_N\cong L_R\oplus Q_R.
\end{align*}

This procedure recovers the original line bundle $L$.  This proves the converse and completes the proof.
\end{proof}
\begin{remark}
    There is a subtlety on the Weyl element, which is an element of $U(2)$ exchanging the two coordinate lines. Passing to the projective group, this Weyl element conjugates $T_S$ by inversion and hence reverses its orientation.  Since the orientation of $T_S$ is already fixed in \Cref{rmk:orientedSO(2)}, this Weyl conjugation is not allowed. Consequently $\chi_1$ and $\chi_2$ are not identified, and the first line bundle $L$ is intrinsically distinguished.
\end{remark}

\begin{definition}
    Let $B$ be a smooth connection on $\widehat P$, and $\nabla_B$ be the induced connection on $\mathfrak{g}_{\widehat P}|_N=\lambda_R\oplus K_R$. $B$ is called \textbf{preserving the reduction} $R$ near $D$ if one of the following equivalent conditions holds on $N$:
    \begin{itemize}
        \item The horizontal distribution restricts to a connection on $R$;
        \item $\nabla_B$ preserves both summands $\lambda_R$ and $K_R$;
        \item $\nabla_B J_R=0$.
    \end{itemize}
\end{definition}
The connection preserving the surface reduction always exists by partition of unity.

\begin{definition}\label{def:model-so3-connection}
For $0<\beta<\frac12$, choose a smooth connection $B$ preserving $R$ near $D$.  The
\textbf{model $\beta$-twisted connection} on $\widehat P|_X$ is
\begin{equation}\label{eq:model}
 B_\beta:=B+2\beta\chi(r)\eta\otimes J_R.
\end{equation}
On $K_R$ the corresponding covariant derivative is
$$
\nabla^{K_R}_{B_\beta}=\nabla^{K_R}_{B}+2\beta\chi\eta\otimes\mathcal J_R.
$$
\end{definition}
The connection one-form is singular because $|\eta|\sim r^{-1}$ in a smooth metric,
but its curvature on $X$ is
$$F_{B_\beta}=F_B+2\beta\,d(\chi\eta)\otimes J_R,$$
which extends smoothly across $D$. For any positively oriented meridian $\gamma_r$ in a normal disk to $D$, the limiting holonomy satisfies
\begin{align*}
\lim_{r\to0}\operatorname{Hol}_{\gamma_r}(B_\beta)
=\exp(-4\pi\beta J_0)=
\rho\left(\diag(e^{-2\pi i\beta},e^{2\pi i\beta})\right).
\end{align*}

\begin{lemma}\label{lem:smooth-model-change}
Given a projective bundle pair $(\widehat P,R)$.
Let $(B,\eta,\chi)$ and $(B',\eta',\chi')$ be two choices of auxiliary data in the construction of the model singular connection, where $B$ and $B'$ are smooth connections on $\widehat P$ preserving
$R$ near $D$, $\eta$ and $\eta'$ are connection one-forms on the oriented
normal circle bundle, normalized by
\begin{align*}
\int_{S^1_{\rm fib}}\eta=\int_{S^1_{\rm fib}}\eta'=2\pi,
\end{align*}
and $\chi,\chi'$ are smooth cutoffs which are equal to $1$ near $D$ and are
supported in the interior of $N$.  Let
\begin{align*}
B_\beta:=B+2\beta\chi\eta\otimes J_R,
\qquad
B'_\beta:=B'+2\beta\chi'\eta'\otimes J_R.
\end{align*}
Then $B'_\beta-B_\beta\in \Gamma(X,T^*X\otimes\mathfrak g_{\widehat P})$
extends across $D$ smoothly. In fact, we have $(B'_\beta-B_\beta)|_N\in \Gamma(N,T^*X\otimes\mathfrak \lambda_R)$.

In particular, the singular part of the model connection is determined by $(R,\beta)$ up
to addition of a smooth adjoint-valued one-form.
\end{lemma}

\begin{proof}
Since $B$ and $B'$ are smooth connections on the same principal bundle,
their difference is a globally defined smooth adjoint-valued one-form $B'-B\in
\Omega^1(\X;\mathfrak g_{\widehat P})$. Thus it remains to consider the singular terms.

Let $\pi_Y:Y\to D$ be the oriented normal circle bundle.  The difference of two connection
one-forms on the same principal $S^1$-bundle is basic.  Hence there exists
a unique one-form $\alpha\in\Omega^1(D)$ such that
\begin{align*}
\eta'-\eta=\pi_Y^*\alpha
\end{align*}
on $Y$.  After radial extension from $Y$ to $N\setminus D$, this difference
is the pullback of $\alpha$ under the tubular projection $\pi_N:N\to D$. Consequently $\eta'-\eta=\pi^*_N\alpha$ extends smoothly across the zero section $D$.

Now we have
\begin{align*}
    B'_\beta-B_\beta
&=(B'-B)+2\beta\bigl(\chi'\eta'-\chi\eta\bigr)\otimes J_R\\
&=(B'-B)+2\beta\Bigl[\chi'(\eta'-\eta)+(\chi'-\chi)\eta\Bigr]\otimes J_R,
\end{align*}
and each term on the right is smooth across $D$. Hence $B'_\beta-B_\beta$ is a global smooth
adjoint-valued one-form on $\X$.

Since both $B$ and $B'$ preserves $R$ near $D$, we have $\nabla_B J_R=\nabla_{B'}J_R=0$. Thus
\begin{align*}
    [B'-B,J_R]=0,
\end{align*}
which gives $(B'-B)|_N\in \Gamma(N,T^*X\otimes\lambda_R)$. It's obvious that $2\beta\chi\eta\otimes J_R\in\Gamma(N,T^*X\otimes\lambda_R)$, therefore we have $(B'_\beta-B_\beta)|_N\in \Gamma(N,T^*X\otimes\lambda_R)$.
\end{proof}

\begin{definition}[Topological labels]\label{def:so3-topological-labels}
For a projective bundle pair $(\widehat P,R)$, we define
\begin{align*}
    w:=w_2(\widehat P)\in H^2(\X,\ZZ_2),\quad
 k:=-\frac14\langle p_1(\mathfrak g_{\widehat P}),[\X]\rangle\in\frac14\ZZ,
 \quad
 l:=-\frac12\langle e(K_R),[D]\rangle\in\frac12\ZZ.
\end{align*}
 The numbers $k$ and $l$ are called the
\textbf{instanton number} and \textbf{monopole number}.
\end{definition}

\begin{remark}
    These three labels are not independent on each other. Actually they satisfy a congruence relation:
    \begin{equation}\label{eq:labelcongruences}
 -4k\equiv\langle\mathfrak P(w),[\X]\rangle\pmod4,
 \qquad
 2l\equiv\langle w,[D]\rangle\pmod2,
\end{equation}
where $\mathfrak P$ is the Pontryagin square.
\end{remark}

\begin{lemma}\label{lem:labels}
If $\widehat P=\PU(E)$ and $R=R_L$ is induced by an ordered splitting
$E|_N=L\oplus Q$ in \Cref{prop:line-reduction}, then
\begin{align*}
 w_2(\widehat P)&=c_1(E)\pmod2,\\
 p_1(\mathfrak{su}(E))&=c_1(E)^2-4c_2(E),\\
 e(K_R)[D]&=2\deg L-c_1(E)\cdot D.
\end{align*}
Consequently
\begin{equation}\label{eq:unitarylabels}
 k=c_2(E)-\frac14c_1(E)^2=:k(E),
 \qquad
 l=\frac12c_1(E)\cdot D-\deg L=:l(E,L).
\end{equation}
\end{lemma}
\begin{proof}
For $\widehat P=\PU(E)$, the adjoint bundle is canonically $\mathfrak g_{\widehat P}\cong\mathfrak{su}(E)$. Thus it suffices to compute the characteristic classes of $\mathfrak{su}(E)$. By the splitting principle for complex vector bundles, we can assume $E=L_1\oplus L_2$ splits as a sum of complex line bundles. Write $x_i=c_1(L_i)$ for $i=1,2$. Then
\begin{align*}
c_1(E)=x_1+x_2,
\qquad
c_2(E)=x_1x_2.
\end{align*}

Equip the splitting with an orthogonal direct-sum Hermitian metric, then an element of $\mathfrak{su}(E)$ has the form
\begin{align*}
\begin{pmatrix}
\sqrt{-1}t & \phi\\
-\phi^* & -\sqrt{-1}t
\end{pmatrix},
\qquad
t\in\mathbb R,\quad
\phi\in\Hom(L_2,L_1).
\end{align*}
Consequently, as an oriented real vector bundle,
\begin{align*}
\mathfrak{su}(E)\cong\underline{\mathbb R}\oplus M_{\mathbb R},
\end{align*}
where the complex line $M:=\Hom(L_2,L_1)=L_1\otimes L_2^{-1}$. Its first Chern class is $c_1(M)=x_1-x_2$.

For any complex line bundle $M$, one has
\begin{align*}
w_2(M_{\mathbb R})=c_1(M)\pmod 2.
\end{align*}
Since the real line summand is trivial, we have
\begin{align*}
w_2(\mathfrak{su}(E))=w_2(M_{\mathbb R})=x_1-x_2=x_1+x_2=c_1(E) \pmod2.
\end{align*}

We next compute the first Pontryagin class.  For a complex line bundle $M$, its underlying oriented real two-plane bundle satisfies
\begin{align*}
p_1(M_{\mathbb R})=c_1(M)^2.
\end{align*}

Since the trivial real line contributes no Pontryagin class, and therefore
\begin{align*}
p_1(\mathfrak{su}(E))=p_1(M_{\RR})=(x_1-x_2)^2=(x_1+x_2)^2-4x_1x_2=c_1(E)^2-4c_2(E).
\end{align*}

It remains to compute the Euler class along $D$.  By Proposition~\ref{prop:line-reduction}, there is an isomorphism of complex line bundles
\begin{align*}
K_R\cong\Hom(Q,L)
=
L\otimes Q^{-1}.
\end{align*}
Thus
\begin{align*}
e(K_R)=c_1(K_R)=c_1\bigl(\Hom(Q,L)\bigr)=
c_1(L)-c_1(Q)\in H^2(D;\mathbb Z).
\end{align*}

Since $E|_D=L\oplus Q$, we have
\begin{align*}
c_1(E)|_D=c_1(L)+c_1(Q),
\end{align*}
therefore
\begin{align*}
e(K_R)=2c_1(L)-c_1(E)|_D.
\end{align*}
Since
\begin{align*}
\deg L:=\langle c_1(L),[D]\rangle,
\qquad
c_1(E)\cdot D:=\langle c_1(E)|_D,[D]\rangle,
\end{align*}
we obtain
\begin{align*}
e(K_R)[D]=2\deg L-c_1(E)\cdot D.
\end{align*}
\end{proof}

\subsection{Weighted Sobolev spaces and local analysis}
For $1\le p<\infty$, define the weighted Sobolev space
\begin{align*}
    W^{m,p}_{B_\beta}(X,\mathfrak{g}_{\widehat P}):=\left\{u\in W^{m,p}_{\mathrm{loc}}(X,\mathfrak{g}_{\widehat P}):\nabla^j_{B_\beta}u\in L^p(X),\; 0\le j\le m\right\},
\end{align*}
equipped with the weighted Sobolev norm
\begin{equation}\label{eq:covweighted}
 \|u\|_{W^{m,p}_{B_\beta}}:=\sum_{j=0}^m\|\nabla_{B_\beta}^{\,j}u\|_{L^p(X)}.
\end{equation}

\begin{proposition}[Independence of the weighted spaces from the auxiliary model]\label{prop:weighted-norm-indep-model}
Let $B_\beta$ and $B'_\beta$ be two model singular connections associated
with the same $(\widehat P,R,\beta)$ as in
Lemma~\ref{lem:smooth-model-change}.  For every $1<p<\infty$ and
$m=0,1,2$, there exists a constant $C\geq1$ such that
\begin{align*}
C^{-1}\|u\|_{W^{m,p}_{B_\beta}}
\leq
\|u\|_{W^{m,p}_{B'_\beta}}
\leq
C\|u\|_{W^{m,p}_{B_\beta}}.
\end{align*}
Consequently the corresponding Sobolev spaces coincide with equivalent norms.   
\end{proposition}

\begin{proof}
By Lemma~\ref{lem:smooth-model-change}, $a:=B'_\beta-B_\beta$ extends to a smooth section of
$T^*\X\otimes\mathfrak g_{\widehat P}$, and satisfies $a|_N\in \Gamma(N,T^*X\otimes\lambda_R)$. The two induced covariant derivatives satisfy
\begin{align*}
\nabla_{B'_\beta}u=\nabla_{B_\beta}u+[a,u].
\end{align*}
Hence
\begin{align*}
\|\nabla_{B'_\beta}u\|_{L^p}
\leq
\|\nabla_{B_\beta}u\|_{L^p}+C\|u\|_{L^p}.
\end{align*}

Applying $\nabla_{B'_\beta}$ once more gives
\begin{align*}
\nabla_{B'_\beta}^2u
=\nabla_{B_\beta}^2u+(\nabla_{B_\beta}a)* u
+a*\nabla_{B_\beta}u+a*a*u.
\end{align*}
Since $a|_N\in \Gamma(N,T^*X\otimes\lambda_R)$ and $[J_R,J_R]=0$,
\begin{align*}
    \nabla_{B_\beta}a=\nabla_B a+2\beta\chi\eta\otimes [J_R,a]=\nabla_B a\in C^\infty.
\end{align*}
Thus
\begin{align*}
\|\nabla_{B'_\beta}^2u\|_{L^p}
\leq
C\left(
\|\nabla_{B_\beta}^2u\|_{L^p}
+
\|\nabla_{B_\beta}u\|_{L^p}
+
\|u\|_{L^p}
\right).
\end{align*}
Now we have
\begin{align*}
\|u\|_{W^{m,p}_{B'_\beta}}
\leq
C\|u\|_{W^{m,p}_{B_\beta}},
\qquad m=0,1,2.
\end{align*}
Interchanging $B_\beta$ and $B'_\beta$ gives the reverse inequality.

\end{proof}

Next we give a local characterization of the weighted Sobolev norm $W^{m,p}_{B_\beta}$ near the divisor with respect to the splitting $\mathfrak g_{\widehat P}|_N=\lambda_R\oplus K_R$, which essentially established in \cite[Proposition~3.7]{KM93}. Indeed, since the reference connection $B$ preserves the reduction, 
we have $\nabla_{B_\beta}J_R=0$. Thus $\nabla_{B_\beta}(f J_R)=df\otimes J_R$ acts trivially on the trivial line bundle summand $\lambda_R=\mathbb R J_R$. On the rotating summand $K_R$, the weighted norms are more complicated.

Let $\nabla$ be any fixed smooth covariant derivative on $K_R$. Define, for $m=1,2$,
\begin{equation}\label{eq:scalarweighted}
\|v\|_{W_D^{m,p}}:=\sum_{j=0}^m
\bigl\|r^{-(m-j)}\nabla^jv\bigr\|_{L^p},\qquad v\in\Gamma(K_R).
\end{equation}
Thus
\begin{align*}
    \|v\|_{W_D^{1,p}}&=\|r^{-1}v\|_{L^p}+\|\nabla v\|_{L^p},\\
    \|v\|_{W_D^{2,p}}&=\|r^{-2}v\|_{L^p}+\|r^{-1}\nabla v\|_{L^p}+\|\nabla^2v\|_{L^p}.
\end{align*}

\begin{lemma}[\cite{KM93}, Lemma 3.5]
\label{lem:angular-poincare}
Let $U=\Delta_z\times\Delta_w$, $D\cap U=\{z=0\}$, where $z=re^{\sqrt{-1}\theta}$. Let $(V,\nabla_0)\to U$ be a
smooth Hermitian vector bundle with a smooth unitary connection.
For $\lambda\in\RR\setminus\ZZ$, define on $U\setminus D$
\begin{align*}
\nabla_\lambda:=\nabla_0+\sqrt{-1}\lambda\,d\theta\,\Id_V.
\end{align*}
Then for every $1<p<\infty$, there is a constant $C$ such that
\begin{align*}
\|r^{-1}s\|_{L^p}+\|\nabla_0s\|_{L^p}
\le
C\bigl(\|\nabla_\lambda s\|_{L^p}+\|s\|_{L^p}\bigr)
\end{align*}
for every $s\in C^\infty_c(U\setminus D;V)$.
\end{lemma}
\begin{proof}
For the model connection $\widehat\nabla_\lambda:=d+\sqrt{-1}\lambda\,d\theta$ on the trivialization of the Hermitian bundle, this is precisely \cite[Lemma~3.5]{KM93}, which gives
\begin{align}\label{eq:KM-Poincare}
\|r^{-1}s\|_{L^p}\le C \|\widehat \nabla_\lambda s\|_{L^p}.
\end{align}
In general write the connection as $\nabla_\lambda=\widehat \nabla_\lambda+A$, where $A$ is a smooth matrix-valued one-form.  Applying \eqref{eq:KM-Poincare} gives
\begin{align*}
\|r^{-1}s\|_{L^p}\le C\|\widehat\nabla_\lambda s\|_{L^p}
\le
C\bigl(\|\nabla_\lambda s\|_{L^p}+\|As\|_{L^p}\bigr)
\le C\bigl(\|\nabla_\lambda s\|_{L^p}+\|A\|_{L^\infty}\|s\|_{L^p}\bigr).
\end{align*}
Since $d\theta=O(r^{-1})$ and $\nabla_0=\nabla_\lambda-\ii\lambda d\theta \Id_V$,
\begin{align*}
    \|\nabla_0 s\|_{L^p}\le \|\nabla_\lambda s\|_{L^p}+C\|r^{-1}s\|_{L^p}\le
C\bigl(\|\nabla_\lambda s\|_{L^p}+\|s\|_{L^p}\bigr).
\end{align*}
\end{proof}

\begin{proposition}[Local form of the weighted norm]
\label{prop:localweighted}
Let $U\subset N$ and write $u=u_\lambda+u_K$ with respect to the orthogonal splitting $\mathfrak g_{\widehat P}|_U=\lambda_R\oplus K_R$.

For $m=1,2$ and $1<p<\infty$, the following are equivalent:
\begin{enumerate}
\item $u\in W_{B_\beta}^{m,p}(U)$;
\item $u_\lambda\in W^{m,p}(U;\lambda_R)$ and $u_K\in W_D^{m,p}(U;K_R)$.

\end{enumerate}
Actually, the two norms are equivalent:
\begin{align*}
   \|u\|_{W_{B_\beta}^{m,p}(K)}\quad \text{  and  }\quad  \|u_\lambda\|_{W^{m,p}(K)}+\|u_K\|_{W_D^{m,p}(K)},
\end{align*}
for any compact subset $K\subset U$.

Moreover, if $p>2$, then
\begin{align*}
W_D^{2,p}\hookrightarrow C^0,
\end{align*}
and every $v\in W_D^{2,p}$ satisfies $v|_D=0$.
\end{proposition}
\begin{proof}
The statement is the projective version of the local norm characterization proved by Kronheimer--Mrowka in \cite[\S3(ii)--(iii), Lemmas~3.5--3.6 and Proposition~3.7]{KM93}.  We explain precisely why their argument applies in the present setting.

As mentioned before, $B$ preserves the reduction implies $\nabla_{B_\beta}u_\lambda=\nabla_Bu_\lambda$. Since $B$ is smooth, the restriction of $B_\beta$ to $\lambda_R$ is equivalent to the ordinary $W^{m,p}$-norm.

We next consider the complex line bundle $(K_R,\mathcal J_R)$. In a local trivialization of $K_R$ and the normal circle bundle $Y$, one may write $\mathcal{J}_R=\ii$ and $\eta=d\theta+\pi^*\zeta$, where $\zeta$ is a smooth one-form on the divisor chart. Then the connection induced by $B_\beta$ on $K_R$ has the form
\begin{align*}
\nabla_{B_\beta}^{K}
=d+2\sqrt{-1}\beta\,d\theta+A_{\rm sm},
\end{align*}
where $A_{\rm sm}$ is a smooth one-form with values in $\operatorname{End}_{\mathbb C}(K_R)$. This is exactly the scalar local model considered in \cite[\S3(ii)]{KM93}.  There they consider the trivial
complex line with connection
\begin{align*}
\nabla_\alpha=\nabla+\sqrt{-1}\alpha\,d\theta,
\end{align*}
for $\alpha\notin\mathbb Z$. \cite[Lemma~3.5]{KM93} (see also \Cref{lem:angular-poincare}) proves the fundamental estimate using Green function
\begin{align*}
\|r^{-1}s\|_{L^p}\leq C_\alpha\|\nabla_\alpha s\|_{L^p}.
\end{align*}
Moreover, since $\nabla_\alpha-\nabla=O(r^{-1})$, \cite[Lemma~3.6]{KM93} shows that for every $k\geq1$, the following two norms are equivalent:
\begin{align*}
    \sum_{j=0}^k\|\nabla_\alpha^j s\|_{L^p},\qquad \sum_{j=0}^k
\|r^{-(k-j)}\nabla^js\|_{L^p}.
\end{align*}

Applying this result with $\alpha=2\beta\in(0,1)$, we have
\begin{align*}
u_K\in W_{B_\beta}^{m,p}
\quad\Longleftrightarrow\quad
u_K\in W_D^{m,p},
\qquad m=1,2,
\end{align*}
with equivalent local norms.

Suppose now that $p>2$ and $v\in W_D^{2,p}$.  On a fixed divisor chart, $r$ is bounded above. Hence we have $v\in W^{2,p}$. Since $p>2$, the Sobolev--Morrey embedding in real dimension four gives
\begin{align*}
W^{2,p}
\hookrightarrow
C^{0,\alpha},
\qquad
\alpha=2-\frac4p>0.
\end{align*}
Thus $v$ has a continuous representative up to $D$.

If $v(x_0)\neq0$ at some $x_0\in D$, continuity gives a neighborhood $V$ of $x_0$ and a constant $c>0$ such that $|v|\geq c>0$ on $V$. Then we have
\begin{align*}
\int_V |r^{-2}v|^p\,dV \geq c^p C^{-1}
\int_0^\varepsilon r^{1-2p}\,dr=\infty,
\end{align*}
because $p>1$.  This contradicts $r^{-2}v\in L^p$. Hence we have $v|_D=0$.
\end{proof}

The multiplication lemma
\cite[Lemma~3.8]{KM93} gives the continuous maps
\begin{equation}\label{eq:weightedproducts}
\begin{split}
 W^{2,p}_{B_\beta}\cdot W^{2,p}_{B_\beta}&\subset W^{2,p}_{B_\beta},\\
 W^{2,p}_{B_\beta}\cdot W^{1,p}_{B_\beta}&\subset W^{1,p}_{B_\beta},\\
 W^{1,p}_{B_\beta}\cdot W^{1,p}_{B_\beta}&\subset L^p.
\end{split}
\end{equation}

\begin{definition}\label{def:projectiveweighted}
Given a projective bundle pair $(\widehat P,R)$, define the \textbf{space of $\beta$-twisted connections} by
\begin{align*}
    \widehat\A_\beta^{1,p}(\widehat P,R)
 :=B_\beta+W^{1,p}_{B_\beta}(T^*X\otimes\mathfrak g_{\widehat P}).
\end{align*}

For the gauge group, we first define the \textbf{smooth projective gauge group} by
\begin{align*}
    \widehat{\mathcal G}^{\infty}(\widehat P,R):=\left\{u\in\Gamma^\infty\!\left(\X,\SO(\mathfrak g_{\widehat P})\right): u(J_R)=J_R
\text{ on a neighborhood of }D
\right\},
\end{align*}
then we define the \textbf{weighted gauge group} to be the completion of $\widehat{\mathcal G}^{\infty}(\widehat P,R)$ with respect to the weighted Sobolev norm $W^{2,p}_{B_{\beta}}(\End(\mathfrak g_{\widehat P}))$:
\begin{align*}
    \widehat{\mathcal G}^{\,2,p}_\beta(\widehat P,R):=\overline{\widehat{\mathcal G}^{\infty}(\widehat P,R)}^{W^{2,p}_{B_\beta}}.
\end{align*}

Finally, we define the \textbf{twisted instanton moduli space} by
\begin{equation}\label{eq:projectiveKMmoduli}
 \widehat\M_\beta(\widehat P,R)
 =\{B\in\widehat\A_\beta^{1,p}(\widehat P,R):F_B^{+}=0\}/
   \widehat\G_\beta^{2,p}(\widehat P,R).
\end{equation}
\end{definition}

\begin{remark}[Subtlety on the $SO(3)$ weighted gauge group]\label{rem:projective-gauge-exception}
There is a subtlely between $SU(2)$ connection and $SO(3)$ connection for the weighted gauge group.

Our definition of \(\widehat{\mathcal G}^{\,2,p}_\beta(\widehat P,R)\) uses the
completion of the \emph{oriented} smooth gauge group. However, the definition in $SU(2)$ setting is more simple. Analogous to \cite{KM93}, it seems that we could define it as
\begin{align}\label{eq:KMweightedgauge}
    \widetilde{\mathcal G}^{\,2,p}_\beta(\widehat P,R):=\left\{
u\in
\Gamma\!\left(X,\SO(\mathfrak g_{\widehat P})\right): u\in W^{2,p}_{B_\beta}\left(\End(\mathfrak g_{\widehat P})\right)
\right\}.
\end{align}
However, at the exceptional value $\beta=\frac{1}{4}$, two definitions are different:
\begin{align*}
    \widehat{\mathcal G}^{\,2,p}_{\frac{1}{4}}(\widehat P,R)\neq\widetilde{\mathcal G}^{\,2,p}_{\frac{1}{4}}(\widehat P,R).
\end{align*}

Indeed, on a local trivialization near $D$, the singular connection on \(\mathfrak g_{\widehat P}\) is $\nabla_{B_{\beta}}=d+2\beta J_0 d\theta$. The connection induced on the endomorphism bundle is therefore
\begin{align*}
\nabla^{\End}_{B_{\beta}}U
=dU+2\beta[J_0,U]d\theta.
\end{align*}
Consider the two endomorphisms
\begin{align*}
S_1:=
\begin{pmatrix}
0&0&0\\
0&1&0\\
0&0&-1
\end{pmatrix},
\qquad
S_2:=
\begin{pmatrix}
0&0&0\\
0&0&1\\
0&1&0
\end{pmatrix}.
\end{align*}
A direct calculation gives $[J_0,S_1]=2S_2$, $[J_0,S_2]=-2S_1$. Thus, on the real two-plane
\begin{align*}
\mathcal W:=\operatorname{span}_{\mathbb R}\{S_1,S_2\}
\subset
\End(\mathfrak g_{\widehat P}),
\end{align*}
if \(\mathcal J_{\mathcal W}\) is defined by
\begin{align*}
\mathcal J_{\mathcal W}S_1=S_2,
\qquad
\mathcal J_{\mathcal W}S_2=-S_1,
\end{align*}
then
\begin{align*}
[J_0,\cdot]\big|_{\mathcal W}
=
2\mathcal J_{\mathcal W},
\end{align*}
and consequently
\begin{align*}
\nabla^{\End}_{B_{\beta}}\big|_{\mathcal W}
=d+4\beta\,\mathcal J_{\mathcal W}d\theta.
\end{align*}
Hence the induced endomorphism singular connection contains an angular weight \(4\beta\), compared with the weight \(2\beta\) occurring on
\(K_R\subset\mathfrak g_{\widehat P}\). The subtle point is $2\beta\notin\ZZ$ but $4\beta$ may be integral for $0<\beta<1/2$.

For $\beta=\frac14$, the weight $4\beta$ is integral.  In fact
\begin{align*}
U(\theta)
=
\cos\theta\,S_1-\sin\theta\,S_2
\end{align*}
satisfies
\begin{align*}
\left(
\partial_\theta+\mathcal J_{\mathcal W}
\right)U=0.
\end{align*}
The eigenvalue $0$ breaks the arguments in \Cref{prop:localweighted}.

More precisely, the matrix
\begin{align*}
u(\theta)
:=
\begin{pmatrix}
-1&0&0\\
0&\cos\theta&-\sin\theta\\
0&-\sin\theta&-\cos\theta
\end{pmatrix}
\in\SO(3)
\end{align*}
is \(2\pi\)-periodic and satisfies $\nabla^{\End}_{B_{1/4}}u=0$. Thus \(u\) satisfies the covariant \(W^{2,p}_{B_\beta}\) condition in \eqref{eq:KMweightedgauge} for every \(p\). However, \(u(\theta)\) does not extend continuously across $D$ and it reverses the oriented reduction
\begin{align*}
u(\theta)J_0=-J_0.
\end{align*}

This phenomenon is absent in the $SU(2)$ setting. The only nonzero angular weights are $\pm2\beta$.
Since $0<2\beta<1$, these weights are never integral.

Our definition excludes this extra Weyl component by fixing the oriented reduction before taking the Banach completion.
\end{remark}

\begin{remark}
As \cite[Proposition 2.4]{KM93}, we can prove the followings properties. The affine space \(\widehat{\mathcal A}^{1,p}_\beta(\widehat P,R)\) is separable, \(\widehat{\mathcal G}^{2,p}_\beta(\widehat P,R)\) is a Banach Lie group, and its action on \(\widehat{\mathcal A}^{1,p}_\beta(\widehat P,R)\) is smooth. The curvature map $B\longmapsto F_B$ is continuous from \(\widehat{\mathcal A}^{1,p}_\beta\) to \(L^p\). These structures are independent on the choice of auxiliary model data, up to a canonical Banach isomorphism. Since we won't use these facts, we omit the proof.
\end{remark}

\begin{lemma}\label{lem:projectiveboundary}
Every $u\in\widehat\G_\beta^{2,p}(\widehat P,R)$ has a continuous representative on $\X$, and $u|_D$ preserves the surface reduction $R$, i.e. $u|_D(J_R)=J_R$.
\end{lemma}
\begin{proof}
By definition, there exists a sequence $u_n\in\widehat{\mathcal G}^{\infty}(\widehat P,R)$ such that
\begin{align*}
u_n\longrightarrow u
\quad\text{in}\quad
W^{2,p}_{B_\beta}
\bigl(\End(\mathfrak g_{\widehat P})\bigr).
\end{align*}

Consider a \(B_\beta\)-parallel orthogonal decomposition of \(\End(\mathfrak g_{\widehat P})\):
\begin{align*}
\End(\mathfrak g_{\widehat P})=\mathcal E_0\oplus\mathcal E_2\oplus\mathcal E_4,
\end{align*}
where
\begin{align*}
    \mathcal E_0&:=
\End(\lambda_R)\oplus\operatorname{span}_{\mathbb R}\{\Id_{K_R},\mathcal J_R\},\\
    \mathcal E_2&:=
\Hom(\lambda_R,K_R)\oplus\Hom(K_R,\lambda_R),\\
    \mathcal E_4&:=
\left\{
T\in\End(K_R):T\mathcal J_R=-\mathcal J_RT
\right\}.
\end{align*}
The decomposition are according to the corresponding singular weights.

Write
\begin{align*}
u_n=u_{n,0}+u_{n,2}+u_{n,4}
\end{align*}
according to this decomposition.  Since the projections are parallel for the model connection, the convergence holds componentwise.

On \(\mathcal E_0\) the singular connection acts trivially.  Hence the \(W^{2,p}_{B_\beta}\)-norm is equivalent to the ordinary \(W^{2,p}\)-norm.  Thus
\begin{align*}
u_{n,0}\longrightarrow u_0
\quad\text{in }W^{2,p}(U;\mathcal E_0).
\end{align*}
Since \(p>2\) and \(\dim_{\mathbb R}X=4\),
\begin{align*}
W^{2,p}(U)\hookrightarrow C^0(U),
\end{align*}
so \(u_0\) has a continuous representative. Moreover, since $u_{n,0}|_D(J_R)=J_R$, we have $u_0|_D(J_R)=J_R$.

On each summand of \(\mathcal E_2\), the singular connection has locally the form $d\pm 2\ii\beta d\theta$. Since $0<2\beta<1$, this parameter is nonintegral. The arguments in \Cref{prop:localweighted} identify the
\(W^{2,p}_{B_\beta}\)-norm with the weighted norm \(W_D^{2,p}\). Consequently
\begin{align*}
u_{n,2}\longrightarrow u_2
\quad\text{in }W_D^{2,p}(U;\mathcal E_2).
\end{align*}
By \Cref{prop:localweighted}, \(u_2\) has a continuous representative, and $u_2|_D=0$.

We next consider \(\mathcal E_4\). The singular connection has the form \(d+4\ii\beta d\theta\).  If $\beta\neq \frac14$, the same argument gives
\begin{align*}
u_{n,4}\longrightarrow u_4
\quad\text{in }W_D^{2,p}(U;\mathcal E_4),
\end{align*}
and hence
\begin{align*}
u_4\in C^0(U;\mathcal E_4),
\qquad
u_4|_D=0.
\end{align*}

It remains to treat the only exceptional value $\beta=\frac14$. Now locally $\nabla_{B_{1/4}}=d+\ii d\theta$. Since $\exp(2\pi\ii)=\Id$, the map $Q(\theta):=\exp(\ii\theta)$ is single-valued on $U\setminus D$. Note that
\begin{align*}
    \nabla(e^{\ii\theta}f)&=e^{\ii\theta}(df+\ii f d\theta)=e^{\ii\theta}\nabla_{B_{1/4}}f,\\
    \nabla^2(e^{\ii\theta}f)&=\nabla\left(e^{\ii\theta}\nabla_{B_{1/4}}f\right)=e^{\ii\theta}\nabla_{B_{1/4}}^2f,
\end{align*}
thus we have
\begin{align*}
    \|f\|_{W^{2,p}_{B_{1/4}}}=\|e^{\ii\theta}f\|_{W^{2,p}}.
\end{align*}

Recall $u_n\in \widehat{\mathcal G}^{\infty}(\widehat P,R)$, which means $u_n(J_R)=J_R$ in a neighborhood of $D$. Moreover, $u_n\in \Gamma(SO(\mathfrak g_{\widehat P}))$ implies there exists function $t_n$ such that
\begin{align*}
    u_n=\begin{pmatrix}
        1 & 0 & 0\\
        0 & \cos t_n & -\sin t_n\\
        0 & \sin t_n & \cos t_n
    \end{pmatrix}=\Id_{\lambda_R}+\cos t_n \Id_{K_R}+\sin t_n \mathcal{J}_R\in\mathcal{E}_0.
\end{align*}
Therefore $u_{n,2}=u_{n,4}=0$ on some neighborhood of \(D\). The convergence
\begin{align*}
    u_{n,4}\longrightarrow u_4\quad\text{in }W_{B_{1/4}}^{2,p}(U;\mathcal E_4)
\end{align*}
implies
\begin{align*}
    e^{\ii\theta}u_{n,4}\longrightarrow e^{\ii\theta}u_4\quad\text{in }W^{2,p}(U;\mathcal E_4)\hookrightarrow C^0(U;\mathcal{E}_4).
\end{align*}
Thus $e^{\ii\theta}u_4$ extends continuously across \(D\) by zero. Since $|u_4|=|e^{\ii\theta} u_4|$, \(u_4\) also extends continuously across \(D\) by $u_4|_D=0$. This proves the required conclusion also at \(\beta=\frac14\).

Combining the three components, \(u\) has a continuous representative on \(U\) with $u|_D=u_0|_D$.
\end{proof}

The following proposition is the twisted $SO(3)$ version of Uhlenbeck's compactness theorem. Since there's no obstruction on the local lift to $SU(2)$, we can apply Kronheimer--Mrowka's results directly.
\begin{proposition}[Local compactness and removable singularity]
\label{prop:pu2local}
Let $x\in D$, and $U$ be a divisor chart around $x$.

\begin{enumerate}
\item
Let $(\widehat P,R)$ be a projective bundle pair over $U$. Then there exists $\varepsilon>0$ and $U'\subset U$ such that every sequence of ASD connections
\begin{align*}
B_n=B_\beta+a_n,\qquad
a_n\in W^{1,p}_{B_\beta}
(U\setminus D;T^*U\otimes\mathfrak g_{\widehat P}),
\end{align*}
satisfying
\begin{align*}
\sup_n\int_{U\setminus D}|F_{B_n}|^2\,dV_g
<\varepsilon
\end{align*}
admits a subsequence $B_{n_j}$ and gauge transformations $u_j$ over $U'\setminus D$, induced by weighted $SU(2)$ gauge transformations on a local lift, such that
\begin{align*}
u_j^*B_{n_j}=B_\beta+b_j,
\qquad
b_j\to b_\infty
\quad\text{in }W^{1,p}_{B_\beta}(U'\setminus D).
\end{align*}
The connection $B_\beta+b_\infty$ is ASD.

\item
Let $(\widehat P,R)$ be a projective bundle pair over $U\setminus\{x\}$.
Let $B$ be an ASD connection on $\widehat P|_{U\setminus D}$. Suppose $B-B_{\beta}\in W^{1,p}_{B_\beta,\mathrm{loc}}(U\setminus \{x\})$ and
\begin{align*}
\int_{U\setminus D}|F_B|^2\,dV_g<\infty.
\end{align*}

Then there exist a projective bundle pair $(\widehat P',R')$ over $U$, a twisted ASD connection $B'\in\widehat{\mathcal A}^{1,p}_\beta(\widehat P',R')$, and a projective bundle-pair isomorphism
\begin{align*}
\Phi:(\widehat P',R')|_{U\setminus\{x\}}
\longrightarrow(\widehat P,R)
\end{align*}
such that $\Phi^*B=B'$ on $U\setminus D$.
\end{enumerate}
\end{proposition}
\begin{proof}
Since $U$ is contractible and $U\setminus\{x\}$ is homotopy equivalent to $S^3$, every $SO(3)$ bundle over either space admits an $SU(2)$ lift. The surface reductions also admit local lifts.
Applying \cite[Proposition~7.4, 7.7]{KM93} and descending the resulting $SU(2)$ connections and gauge transformations give (i) and (ii).
\end{proof}

\subsection{Projective compactness}
\begin{proposition}[Singular Chern--Weil formula]\label{prop:CW}
Given a projective bundle pair $(\widehat P,R)$. Let $(w,k,l)$ be its topological labels and let $B\in\widehat\A_\beta^{1,p}(\widehat P,R)$. Then we have
\begin{equation}\label{eq:CWprojective}
 \frac1{8\pi^2}\int_X\tr_3(F_B\wedge F_B)
 =k+2\beta l-\beta^2D^2.
\end{equation}
If $B$ is anti-self-dual, then
\begin{equation}\label{eq:CWenergy}
 \frac1{8\pi^2}\int_X|F_B|^2\,dV_g
 =k+2\beta l-\beta^2D^2.
\end{equation}
\end{proposition}

\begin{proof}
Recall $B_\beta=B_0+2\beta\chi\eta\otimes J_R$. Set
$\Omega:=d(\chi\eta)$ on $X$. Since $\eta=d\theta+\pi^*\zeta$ for a smooth $1$-form $\zeta$ on $D$, $\Omega$ extends to a smooth closed two-form on $\X$. For any closed $2$-form $\phi$,
\begin{align*}
    \int_\X\Omega\wedge\phi=\lim_{\varepsilon\to0}\int_{\X\setminus N_\varepsilon}d(\chi\eta)\wedge \phi=\lim_{\varepsilon\to0}\int_{\dd\left(\X\setminus N_\varepsilon\right)}\chi\eta\wedge\phi=\lim_{\varepsilon\to0}\int_{-Y_\varepsilon}\eta\wedge\phi=-2\pi\int_D\phi.
\end{align*}
Therefore $\Omega$ represents $-2\pi c_1(D)$.
In particular, $$\int_{\X}\Omega\wedge\Omega=4\pi^2D^2.$$

Since $B_0$ preserves $R$, there exists a closed real
two-form $f$ on $N$ such that $F_{B_0}|_N=f\otimes J_R$.
The induced connection on the complex line bundle
$(K_R,\mathcal J_R)$ has curvature $f\mathcal J_R$.
Hence
\begin{align*}
e(K_R)=-\left[\frac{f}{2\pi}\right].
\end{align*}
Thus
\begin{align*}
\int_N f\wedge\Omega=\langle-2\pi[D], [f]\rangle
=4\pi^2\langle e(K_R),[D]\rangle
=-8\pi^2l.
\end{align*}

The curvature of the model connection satisfies $F_{B_\beta}=F_{B_0}+2\beta\Omega\otimes J_R$.
Since $\tr_3(J_R^2)=-\frac12$, it follows that
\begin{align*}
\begin{aligned}
\tr_3(F_{B_\beta}\wedge F_{B_\beta})
=\tr_3(F_{B_0}\wedge F_{B_0})-2\beta f\wedge\Omega
-2\beta^2\Omega\wedge\Omega.
\end{aligned}
\end{align*}
The ordinary Chern--Weil formula gives
\begin{align*}
\frac1{8\pi^2}\int_{\X}
\tr_3(F_{B_0}\wedge F_{B_0})
=-\frac14
\langle p_1(\mathfrak g_{\widehat P}),[\X]\rangle
=k.
\end{align*}
Combining these identities yields
\begin{align*}
\frac1{8\pi^2}\int_X
\tr_3(F_{B_\beta}\wedge F_{B_\beta})
=k+2\beta l-\beta^2D^2.
\end{align*}

The same identity holds for every smooth perturbation
$B_\beta+a$ with $a$ taking values in $\lambda_R$
near $D$. Indeed, along the path
$B_t=B_\beta+ta$, the transgression formula gives
\begin{align*}
\frac{d}{dt}\tr_3(F_{B_t}\wedge F_{B_t})
=2d\,\tr_3(a\wedge F_{B_t}).
\end{align*}
The three-form on the right extends smoothly across
$D$, so its integral over $X$ vanishes.

Smooth perturbations are dense in $W^{1,p}_{B_\beta}$ and the curvature map is continuous from $W^{1,p}_{B_\beta}$ to $L^p$.
Since $p>2$ and $X$ has finite volume, the curvature
integral is continuous in this topology.
Thus \eqref{eq:CWprojective} holds for every
$B\in\widehat{\mathcal A}^{1,p}_\beta(\widehat P,R)$.

Finally, if $B$ is ASD, since the curvature is always skew-Hermitian
\begin{align*}
\tr_3(F_B\wedge F_B)=-\tr_3(F_B\wedge*F_B)=\tr_3(F_B\wedge*F_B^*)=|F_B|^2\,dV_g,
\end{align*}
which proves \eqref{eq:CWenergy}.
\end{proof}

The extension of a limiting bundle across a concentration point can differ from the
original extension. We need define the relative topological numbers to record the defect.

\begin{definition}\label{def:relative}

Let $B$ be an adapted ball,
and let $(\widehat P_-,R_-)$ and
$(\widehat P_+,R_+)$ be projective bundle pairs
defined on a neighborhood of $\overline B$.
Let $\Psi:(\widehat P_-,R_-)|_C\to (\widehat P_+,R_+)|_C$ be a projective bundle-pair isomorphism on a collar $C$ of $\partial B$.

The \textbf{relative local labels} of $\Psi$,
denoted by $(k_B(\Psi),l_B(\Psi))$,
are defined by gluing the two bundle pairs
over oppositely oriented copies of $B$
using $\Psi$, with the target pair on the
positively oriented copy.

If $B\cap D=\varnothing$, the resulting
projective bundle over $S^4$ defines
\begin{align*}
k_B(\Psi)=-\frac14\langle p_1,[S^4]\rangle,
\qquad l_B(\Psi)=0.
\end{align*}
If $B$ is a divisor chart, the resulting
projective bundle pair over $(S^4,S^2)$ defines
\begin{align*}
k_B(\Psi)=-\frac14\langle p_1,[S^4]\rangle,
\qquad
l_B(\Psi)=-\frac12\langle e(K),[S^2]\rangle.
\end{align*}

For determinant-fixed rank-two bundle-pair
isomorphisms, the relative labels are those
of their projectivizations.
\end{definition}

\begin{lemma}[Integrality and excision]
\label{lem:localintegral}
The relative labels defined in \Cref{def:relative} satisfy
\begin{align*}
(k_B(\Psi),l_B(\Psi))\in\ZZ^2,
\end{align*}
with $l_B(\Psi)=0$ if $B\cap D=\varnothing$.

Let $(\widehat P_i,R_i)$, $i=0,1$, be two
projective bundle pairs over $(\X,D)$, with
topological labels $(w_i,k_i,l_i)$.
Let $B_1,\ldots,B_m\subset\X$ be pairwise
disjoint adapted balls, and let $U$ be an open
neighborhood of
$\X\setminus\bigcup_{j=1}^m\operatorname{int}B_j$.
Suppose there is a projective bundle-pair
isomorphism
\begin{align*}
\Psi:(\widehat P_1,R_1)|_U
\longrightarrow
(\widehat P_0,R_0)|_U.
\end{align*}
Set $(k_j,l_j):=(k_{B_j}(\Psi),l_{B_j}(\Psi))$.
Then
\begin{equation}\label{eq:projectiveexcision}
k_0-k_1=\sum_{j=1}^m k_j,
\qquad
l_0-l_1=\sum_{j=1}^m l_j,
\qquad
w_0=w_1.
\end{equation}
\end{lemma}

\begin{proof}
For each adapted ball $B$, the gluing construction
in \Cref{def:relative} gives a principal $\SO(3)$ bundle $\widehat Q$ over $S^4$, together with an oriented $\SO(2)$ reduction along $S^2$ when $B\cap D\neq\varnothing$.
Since $H^2(S^4;\ZZ_2)=0$, the bundle $\widehat Q$ admits an $\SU(2)$ lift $E$. In the divisor chart case, the reduction lifts to a splitting $E|_{S^2}=L\oplus L^{-1}$.
By \Cref{lem:labels},
\begin{align*}
k_B(\Psi)=\langle c_2(E),[S^4]\rangle,
\qquad
l_B(\Psi)=-\langle c_1(L),[S^2]\rangle.
\end{align*}
Thus both relative labels are integers.

We now prove the excision identities. Choose smooth connections $A_i$ on $\widehat P_i$ preserving $R_i$ near $D$, such that $A_1=\Psi^*A_0$ on $U$.
Let $p_i$ and $e_i$ be the Chern--Weil forms
representing $p_1(\mathfrak g_{\widehat P_i})$ and
$e(K_{R_i})$, respectively. Then $A_1=\Psi^*A_0$ on $U$ implies
\begin{align*}
p_0|_U=p_1|_U,
\qquad
e_0|_{U\cap D}=e_1|_{U\cap D}.
\end{align*}

For each $B_j$, the restrictions of these
connections glue to connections on the
doubled bundle pair
$(\widehat Q_j,R_j')$ over $(S^4,S^2)$.
By the orientation convention in \Cref{def:relative}, integration gives
\begin{align*}
    \langle p_1(\mathfrak g_{\widehat P_0}),
[\X]\rangle
-\langle p_1(\mathfrak g_{\widehat P_1}),
[\X]\rangle=\int_\X (p_0-p_1)=\sum_{j=1}^m\int_{B_j}(p_0-p_1)=\sum_{j=1}^m
\langle p_1(\mathfrak g_{\widehat Q_j}),
[S^4]\rangle,
\end{align*}
and similarly,
\begin{align*}
    \langle e(K_{R_0}),[D]\rangle
-\langle e(K_{R_1}),[D]\rangle
=\int_D(e_0-e_1)=\sum_{j=1}^m\int_{D\cap B_j}(e_0-e_1)=\sum_{j=1}^m
\langle e(K_{R_j'}),[S^2]\rangle.
\end{align*}l
The definitions of $k$ and $l$ now give
the first two identities in
\eqref{eq:projectiveexcision}.

Finally, excision gives $H^2(\X,U;\ZZ_2)=0$. Hence the restriction map
\begin{align*}
H^2(\X;\ZZ_2)\longrightarrow H^2(U;\ZZ_2)
\end{align*}
is injective. Since
$\widehat P_0|_U\simeq\widehat P_1|_U$,
we obtain $w_0=w_1$.
\end{proof}

\begin{lemma}\label{lem:localenergydefect}
Let $V\subset\X$ be an adapted closed ball, and let $C\subset V$ be a collar of $\partial V$. Let $(\widehat P_i,R_i)$, $i=0,1$, be smooth projective bundle pairs over a neighborhood of $V$.

Let $B_n$ be $\beta$-twisted ASD connections on $(\widehat P_0,R_0)$, and let $B_\infty$ be a $\beta$-twisted ASD connection on $(\widehat P_1,R_1)$.

Suppose there exist projective bundle-pair isomorphisms
\begin{align*}
\Psi_n:(\widehat P_1,R_1)|_C
\longrightarrow(\widehat P_0,R_0)|_C,
\end{align*}
such that
\begin{align*}
\Psi_n^*B_n\to B_\infty
\quad\text{in }W^{1,p}_\beta(C\setminus D).
\end{align*}
Assume that, for sufficiently large $n$,
\begin{align*}
(k_V(\Psi_n),l_V(\Psi_n))=(\kappa,\lambda).
\end{align*}

Then
\begin{align}\label{eq:localenergydefect}
    \lim_{n\to\infty}\frac1{8\pi^2}\left(\int_{V\setminus D}|F_{B_n}|^2\,dV_g-\int_{V\setminus D}|F_{B_\infty}|^2\,dV_g\right)=\kappa+2\beta\lambda.
\end{align}
\end{lemma}

\begin{proof}
By the density of smooth gauges in the local
weighted gauge spaces of \cite[\S3]{KM93},
together with the continuity of the gauge action,
we may approximate the comparisons $\Psi_n$
by smooth bundle-pair isomorphisms without
changing their relative labels or the stated
convergence. We therefore assume that
$\Psi_n$ are smooth.

Set $a_n:=\Psi_n^*B_n-B_\infty$ on $C\setminus D$. By assumption,
\begin{align*}
\|a_n\|_{W^{1,p}_\beta(C\setminus D)}
\to0.
\end{align*}

Choose a smooth cutoff function $\chi$ on $V$
which equals $1$ away from a smaller collar
of $\partial V$ and equals $0$ near $\partial V$.
On the positive copy of $V$, replace $B_n$ in the collar by the connection $\Psi_{n*}(B_\infty+\chi a_n)$.
This agrees with $B_n$ where $\chi=1$ and with
$\Psi_{n*}B_\infty$ near $\partial V$.
Consequently, it glues to $B_\infty$ on the
negative copy of $V$, giving a twisted
connection $\widetilde B_n$ on the doubled
projective bundle pair over $(S^4,S^2)$,
or over $S^4$ if $V\cap D=\varnothing$.

The doubled surface $S^2$ has self-intersection zero.
Thus \Cref{prop:CW} gives
\begin{align*}
\frac1{8\pi^2}\int_{S^4\setminus S^2}
\tr_3(F_{\widetilde B_n}\wedge
F_{\widetilde B_n})
=\kappa+2\beta\lambda.
\end{align*}

Since the curvature map is continuous from $W^{1,p}_\beta$ to $L^p$ \cite[Lemma~3.8]{KM93}, and
\begin{align*}
(B_\infty+\chi a_n)-(B_\infty+a_n)=(\chi-1)a_n\to 0\quad\text{in }W^{1,p}_\beta(C\setminus D),
\end{align*}
the corresponding curvatures differ by a
quantity tending to zero in $L^p$.
As $p>2$, their Chern--Weil four-forms
converge in $L^1$.

It follows that
\begin{align*}
    \int_{S^4\setminus S^2}\tr_3(F_{\widetilde B_n}\wedge F_{\widetilde B_n})=\int_{V\setminus D}\tr_3(F_{B_n}\wedge F_{B_n})-\int_{V\setminus D}\tr_3(F_{B_\infty}\wedge F_{B_\infty})+o(1).
\end{align*}

Since $B_n$ and $B_\infty$ are anti-self-dual,
\begin{align*}
\tr_3(F_{B_n}\wedge F_{B_n})
=|F_{B_n}|^2\,dV_g,
\qquad
\tr_3(F_{B_\infty}\wedge F_{B_\infty})
=|F_{B_\infty}|^2\,dV_g.
\end{align*}
Substitution and passing to the limit prove
\eqref{eq:localenergydefect}.
\end{proof}

\begin{theorem}[$PU(2)$ form of Kronheimer--Mrowka compactness]\label{thm:puKM}
Given an $SO(3)$ bundle with a surface reduction $(\widehat P,R)$ and fixed
$0<\beta<\frac12$.  Let $B_n\in\widehat\A_\beta^{1,p}(\widehat P,R)$ be
anti-self-dual connections.  Then, after passing to a subsequence, there exist
\begin{enumerate}
\item a finite set $S\subset\X$;
\item an $SO(3)$ bundle with a surface reduction $(\widehat P_\infty,R_\infty)$ over $(\X,D)$;
\item a $\beta$-twisted ASD connection
$B_\infty\in\widehat\A_\beta^{1,p}(\widehat P_\infty,R_\infty)$;
\item projective pair isomorphisms
$$
 \Psi_n:(\widehat P_\infty,R_\infty)|_{\X\setminus S}
 \longrightarrow(\widehat P,R)|_{\X\setminus S};
$$
\item integral defect labels $(k_z,l_z)\in\ZZ^2$ for any $z\in S$, with $l_z=0$ for $z\notin D$,
\end{enumerate}
such that $\Psi_n^*B_n\to B_\infty$ in weighted $W^{1,p}$ on every relatively compact subsets of $\X\setminus S$.  If $(k,l,w)$ and $(k_\infty,l_\infty,w_\infty)$ are the topological labels of $(\widehat P,R)$ and $(\widehat P_\infty,R_\infty)$,
then
\begin{equation}\label{eq:PUconservation}
 k=k_\infty+\sum_{z\in S}k_z,
 \qquad
 l=l_\infty+\sum_{z\in S}l_z,
 \qquad
 w_\infty=w.
\end{equation}
Moreover,
\begin{equation}\label{eq:PUmeasure}
 |F_{B_n}|^2dV_g\rightharpoonup
 |F_{B_\infty}|^2dV_g
 +8\pi^2\sum_{z\in S}(k_z+2\beta l_z)\delta_z.
\end{equation}
\end{theorem}

\begin{theorem}
\label{thm:puKM}
Let $(\widehat P,R)$ be a projective bundle
pair over $(\X,D)$.
Let $B_n\in\widehat{\mathcal A}^{1,p}_\beta(\widehat P,R)$ be a sequence of ASD connections.

Then, after passing to a subsequence, there exist
\begin{enumerate}
\item a finite set $S\subset\X$;
\item a projective bundle pair $(\widehat P_\infty,R_\infty)$ over $(\X,D)$;
\item a $\beta$-twisted ASD connection $B_\infty\in\widehat{\mathcal A}^{1,p}_\beta(\widehat P_\infty,R_\infty)$;
\item projective bundle-pair isomorphisms
\begin{align*}
\Psi_n:
(\widehat P_\infty,R_\infty)|_{\X\setminus S}
\longrightarrow
(\widehat P,R)|_{\X\setminus S},
\end{align*}
locally of weighted $W^{2,p}$ regularity;
\item integers $(k_z,l_z)\in\ZZ^2$ for $z\in S$,
with $l_z=0$ if $z\notin D$;
\end{enumerate}
such that, for every compact subset $V\subset\X\setminus S$,
\begin{align*}
\Psi_n^*B_n\to B_\infty
\quad\text{in }W^{1,p}_{B_{\beta,\infty}}
(V\setminus D),
\end{align*}
where $B_{\beta,\infty}$ is any model connection
for $(\widehat P_\infty,R_\infty)$.
The convergence is smooth on compact subsets
of $X\setminus S$, after choosing smooth gauges.

For each $z\in S$, the pair $(k_z,l_z)$ is the
eventually constant relative label of $\Psi_n$
on a collar of any sufficiently small adapted
ball centered at $z$.
Let $(w,k,l)$ and $(w_\infty,k_\infty,l_\infty)$
denote the topological labels of the two bundle
pairs. Then
\begin{align}\label{eq:PUconservation}
    k=k_\infty+\sum_{z\in S}k_z,\quad l=l_\infty+\sum_{z\in S}l_z,\quad w_\infty=w.
\end{align}
Moreover, as finite Radon measures on $\X$,
\begin{equation}\label{eq:PUmeasure}
|F_{B_n}|^2\,dV_g
\rightharpoonup
|F_{B_\infty}|^2\,dV_g+8\pi^2\sum_{z\in S}(k_z+2\beta l_z)\delta_z.
\end{equation}
\end{theorem}

\begin{proof}
By \Cref{prop:CW},
\begin{align*}
\frac1{8\pi^2}\int_X|F_{B_n}|^2\,dV_g
=k+2\beta l-\beta^2D^2
\end{align*}
is independent of $n$.
After passing to a subsequence, the measures
$\mu_n:=|F_{B_n}|^2dV_g$ converge weakly
to a finite Radon measure $\mu$ on $\X$.

Fix a sufficiently small constant
$\varepsilon_0>0$ for which the local
small-energy compactness theorem
\Cref{prop:pu2local} applies on the coordinate
balls of a finite cover of $\X$.
Set
\begin{align*}
S:=\{x\in\X:\mu(\{x\})\geq\varepsilon_0\}.
\end{align*}
Since $\mu(\X)<\infty$, the set $S$ is finite.

Put $\Omega:=\X\setminus S$. For every $x\in\Omega$, there exists a sufficiently small coordinate ball $U_x$ such that $\mu(\overline{U_x})<\varepsilon_0$.
Thus
\begin{align*}
\int_{U_x}|F_{B_n}|^2\,dV_g
<\varepsilon_0
\end{align*}
for all sufficiently large $n$.

By \Cref{prop:pu2local}(i), there exist a subsequence, and local projective gauge transformations $u_{n,x}$ such that $u_{n,x}^*B_n-B_\beta$ converges strongly in
$W^{1,p}_{B_\beta}(U_x\setminus D)$.
The patching argument of \cite[Proposition~7.5]{KM93} applies to these projective gauge transformations. Indeed, on a good cover the local $\SU(2)$
lifts of the transition functions may differ
on triple intersections by central signs $\{\pm1\}$.
These signs disappear under the projection
$\SU(2)\to\PU(2)$, so the projective cocycle
identities are preserved by the patching
construction.

Consequently, after passing to a subsequence, there exist projective gauge transformations $u_n$ on $\widehat P|_{\Omega\setminus D}$,
with continuous boundary values preserving
the oriented reduction along $D\cap\Omega$,
and
\begin{align*}
a_\infty\in
W^{1,p}_{B_\beta,\mathrm{loc}}
(\Omega\setminus D;
T^*\Omega\otimes\mathfrak g_{\widehat P})
\end{align*}
such that, for every $V\subset\subset\Omega$,
\begin{align*}
\lim_{n\to\infty}
\left\|
u_n^*B_n-(B_\beta+a_\infty)
\right\|_{W^{1,p}_{B_\beta}(V\setminus D)}
=0.
\end{align*}

Set $B_\infty^\circ:=B_\beta+a_\infty$. Since the curvature map is continuous from weighted $W^{1,p}$ to $L^p$, we have $F_{B_\infty^\circ}^+=0$. Standard elliptic regularity for the ASD equation gives smooth convergence away from $D$.

Since the convergence is strong in weighted
$W^{1,p}$ and $p>2$, the corresponding
curvatures converge strongly in $L^2$ on
compact subsets of $\X\setminus S$. Hence $\mu|_{\X\setminus S}=|F_{B_\infty}|^2dV_g$.

We next extend the limiting pair and connection across $S$. For each $z\in S$, choose a sufficiently small adapted ball $V_z$ centered at $z$. The limiting connection $B^\circ_\infty$ has finite curvature energy on $V_z\setminus D$ and belongs to the local weighted $W^{1,p}$ space on $V_z\setminus\{z\}$.
Then \Cref{prop:pu2local}(ii) gives an extension of the projective bundle pair and the ASD connection across $z$.
Since the balls $V_z$ can be chosen pairwise disjoint, these extensions yield a smooth projective bundle pair $(\widehat P_\infty,R_\infty)$ over $(\X,D)$ and an ASD connection $B_\infty\in\widehat{\mathcal A}^{1,p}_\beta(\widehat P_\infty,R_\infty)$.

Choose pairwise disjoint adapted balls $V_z$
centered at the points of $S$, with collars
$C_z\subset\X\setminus S$.
The comparison maps $\Psi_n|_{C_z}$ determine
relative labels as in \Cref{def:relative}.
The annular gauge-comparison argument in
\cite[\S7(iv)]{KM93} applies to local
$\SU(2)$ lifts on each $C_z$.
It gives, after a further subsequence,
a fixed homotopy class of the collar
comparison for every $z\in S$.
Consequently the relative labels are
independent of $n$ for all sufficiently
large $n$.
Their integrality and the identities
\eqref{eq:PUconservation} follow from
\Cref{lem:localintegral}.

It remains to identify the curvature measure. Since $\mu$ agrees with $|F_{B_\infty}|^2dV_g$ outside $S$,
there exist $m_z\geq0$ such that
\begin{align*}
\mu=|F_{B_\infty}|^2dV_g+\sum_{z\in S}m_z\delta_z.
\end{align*}
Choose each $V_z$ so that
$\mu(\partial V_z)=0$.
The convergence on $C_z$ and the stabilization
of the relative labels allow us to apply
\Cref{lem:localenergydefect}, obtaining
\begin{align*}
    m_z=\lim_{n\to\infty}
\left(
\int_{V_z\setminus D}|F_{B_n}|^2\,dV_g
-\int_{V_z\setminus D}
|F_{B_\infty}|^2\,dV_g
\right)=8\pi^2(k_z+2\beta l_z).
\end{align*}
This proves \eqref{eq:PUmeasure}.
\end{proof}

\begin{proposition}\label{prop:puadmissible}
Every defect labels in \cref{thm:puKM} satisfy
\begin{equation}\label{eq:admissibility}
 k_z\geq0,
 \qquad
 k_z+l_z\geq0.
\end{equation}
\end{proposition}

\begin{proof}
Since we always have an $SU(2)$ lift locally, \cite[Proposition~8.1]{KM93} applies and gives finitely many admissible sphere-pair $(S^4,S^2)$ labels $(\kappa_j,\lambda_j)\in\ZZ^2$ and ordinary
instanton charges $m_r\in\ZZ_{>0}$ satisfying
\begin{align*}
k_z=\sum_j\kappa_j+\sum_r m_r,
\qquad
l_z=\sum_j\lambda_j.
\end{align*}

It suffices to show that $\kappa_j\geq0$ and
$\kappa_j+\lambda_j\geq0$ for every singular bubble.
Since $H^2(S^4;\mathbb Z_2)=0$, its $SO(3)$ bundle admits an $SU(2)$ lift $E_j\to S^4$. The oriented surface reduction determines an ordered splitting $E_j|_{S^2}=L_j\oplus L_j^{-1}$.
By \Cref{lem:labels},
\begin{align*}
    \kappa_j=\langle c_2(E_j),[S^4]\rangle,\qquad \lambda_j=-\langle c_1(L_j),[S^2]\rangle.
\end{align*}
Since $0<2\beta<1$, \cite[Proposition~2.10 and equation~(45)]{KM11} gives the desired inequality.

\end{proof}

\begin{corollary}\label{cor:quantum}

Every nonzero divisor defect $(k_z+2\beta l_z)\delta_z$ has mass at least 
\begin{align*}
    \varepsilon_\beta:=\min\{2\beta,1-2\beta\}>0.
\end{align*}
\end{corollary}
\begin{proof}
Note that
\begin{align*}
    k_z+2\beta l_z=(1-2\beta)k_z+2\beta(k_z+l_z).
\end{align*}
Then \eqref{eq:admissibility} gives the claimed inequality.
\end{proof}

\section{Determinant-fixed lifts and bundle-pair topology}\label{sec:detlift}
We now return to our parabolic setting: a rank-two Hermitian bundle $E$, determinant line bundle $P$, and line subbundle $L\subset E|_D$. 
\subsection{Prescribed determinant lifts and twisted instantons}

\begin{proposition}[Prescribed determinant lifts]\label{prop:lift}
Let $\widehat P$ be a $\PU(2)$-bundle on $\X$, and let $P$ be a Hermitian line bundle. There exists a Hermitian rank-two bundle $E$ together with unitary bundle isomorphisms
\begin{align*}
\det E\simeq P,
\qquad
PU(E)\simeq \widehat P
\end{align*}
if and only if
\begin{equation}\label{eq:w2lift}
w_2(\widehat P)=c_1(P)\pmod 2.
\end{equation}

If \eqref{eq:w2lift} holds, then the isomorphism classes of triples
\begin{align*}
(E,\varphi,\Phi),
\qquad
\varphi:\det E\xrightarrow{\sim}P,
\qquad
\Phi:PU(E)\xrightarrow{\sim}\widehat P,
\end{align*}
form a torsor for $H^1(\X;\ZZ_2)$, where $(E,\varphi,\Phi)\simeq(E',\varphi',\Phi')$ means there exists a bundle isomorphism $f:E\xrightarrow{\sim}E'$ such that
\begin{align*}
    \varphi'\circ \det f=\varphi,\qquad \Phi'\circ PU(f)=\Phi.
\end{align*}
\end{proposition}

\begin{proof}
Consider the double covering
\begin{align*}
1\longrightarrow \{\pm I\}
\longrightarrow U(2)
\xrightarrow{\ \ q\ \ } PU(2)\times U(1)
\longrightarrow 1,
\qquad
u\longmapsto([u],\det u).
\end{align*}
The pair $(\widehat P,P)$ determines a principal
$PU(2)\times U(1)$-bundle over $X$.  A lift through the above double covering is exactly a Hermitian rank-two bundle $E$ equipped with identifications $PU(E)\simeq\widehat P$, $\det E\simeq P$.

By the standard obstruction theory, there is a universal class
\begin{align*}
e_q\in H^2\bigl(B(PU(2)\times U(1));\ZZ_2\bigr)
\end{align*}
such that a $PU(2)\times U(1)$-bundle lifts to a $U(2)$-bundle if
and only if the pullback of $e_q$ vanishes. In degree two,
\begin{align*}
H^2\bigl(B(PU(2)\times U(1));\ZZ_2\bigr)
\simeq
H^2(BPU(2);\ZZ_2)\oplus H^2(BU(1);\ZZ_2).
\end{align*}

The restriction of the extension to $PU(2)\times\{1\}$ is
\begin{align*}
1\longrightarrow\ZZ_2
\longrightarrow SU(2)
\longrightarrow PU(2)
\longrightarrow1,
\end{align*}
whose obstruction class is the second Stiefel--Whitney class $w_2$, which is well-known in spin geometry. Its restriction to $\{1\}\times U(1)$ is
\begin{align*}
1\longrightarrow\ZZ_2
\longrightarrow U(1)
\xrightarrow{\,z\mapsto z^2\,}
U(1)
\longrightarrow1.
\end{align*}
The obstruction to the square-root of a line bundle is just the mod-$2$ reduction of the first Chern class. Hence $e_q=(w_2,[c_1])$. Thus the lifting obstruction of $(\widehat P,P)$ is therefore the pullback
\begin{align*}
    w_2(\widehat P)+c_1(P)\pmod2.
\end{align*}

Once one lift exists, given a good cover $\{U_i\}$, the transition functions of $(\widehat P,P)$ are $g_{ij}:U_{ij}\to PU(2)\times U(1)$. A lift means determine $\tilde g_{ij}:U_{ij}\to U(2)$ satisfies $q(\tilde g_{ij})=g_{ij}$. Any two lifts satisfy $q(\tilde g_{ij})=q(\tilde g_{ij}')$, thus gives a cocycle
\begin{align*}
    \epsilon_{ij}:=\tilde g_{ij}'\tilde g_{ij}^{-1}:U_{ij}\longrightarrow \ker q=\ZZ_2.
\end{align*}
It's easy to check two lifts coincide if and only if the cohomology class $[\epsilon_{ij}]\in H^1(\X,\ZZ_2)$. On the other hand, any class $\epsilon\in H^1(\X,\ZZ_2)$ gives a flat line bundle $L_\epsilon$ satisfying $L_\epsilon^{\otimes 2}\simeq \underline{\CC}$. If $(E,\varphi,\Phi)$ is a lift, it naturally induces a new lift $(E\otimes L_\epsilon,\widetilde \varphi,\widetilde \Phi)$. Thus the set of isomorphism classes of triples $(E,\varphi,\Phi)$ is a torsor
for $H^1(\X;\ZZ_2)$.

\end{proof}

Let $R_L$ be the oriented projective reduction from \cref{prop:line-reduction}, and
let $B_\beta$ be the projective model \eqref{eq:model} on $\PU(E)$.  By
\cref{prop:detproj}, there is a unique unitary connection $A_0$ on $E|_X$ whose
determinant is $a_{P,\mathrm{par}}$ (see \eqref{eq:detconnection},\eqref{eq:dethnormalized}) and whose projective part is $B_\beta$. 

\begin{definition}\label{def:gaugemoduli}
Define the space of \textbf{determinant-fixed twisted connections} 
\begin{align*}
 \A_\beta^{1,p}(E,L;a_{P,\mathrm{par}}):=A_0+W^{1,p}_{B_\beta}(T^*X\otimes\mathfrak{su}(E)),
\end{align*}
and \textbf{determinant-fixed weighted gauge group}
\begin{align*}
    \G_\beta^{2,p}(E,L):=\{u\in\Gamma(X,\SU(E)):u\in W^{2,p}_{A_0}(\End E)\}.
\end{align*}
The determinant-fixed weighted quotient and \textbf{twisted instanton moduli space} are
\begin{align*}
 \B_\beta^p(E,L;a_{P,\mathrm{par}})
 &:=\A_\beta^{1,p}/\G_\beta^{2,p},\\
 \M^{\KM}_\beta(E,L;a_{P,\mathrm{par}})
 &:=\{A\in\A_\beta^{1,p}:(F_A^\perp)^+=0\}/\G_\beta^{2,p}\subset \B_\beta^p(E,L;a_{P,\mathrm{par}}).
\end{align*}
\end{definition}

\begin{proposition}[Projectivization of the determinant-fixed model]
\label{prop:projectivizationweighted}
For every $A\in\A_\beta^{1,p}(E,L;a_{P,\mathrm{par}})$, the projective connection
$A^{\rm proj}$ belongs to
$\widehat\A_\beta^{1,p}(\PU(E),R_L)$ and $\det A=a_{P,\mathrm{par}}$.  Every
$u\in\G_\beta^{2,p}(E,L)$ induces
$[u]\in\widehat\G_\beta^{2,p}(\PU(E),R_L)$ and
$$
 (u\cdot A)^{\rm proj}=[u]\cdot A^{\rm proj}.
$$
Conversely, a projective comparison gauge lifts to a determinant-one comparison on a
given domain exactly when its central double-cover obstruction vanishes.
\end{proposition}
\begin{proof}
In the local trace-free representative of \cref{prop:detproj}, the difference
$A-A_0$ is already $\mathfrak{su}(E)$-valued, and projectivization simply retains this
adjoint-valued remainder.  Thus the weighted norm is exactly the projective norm of
\cref{sec:weighted}.  A determinant-one gauge has a projective image, and functoriality
of induced connections gives the displayed gauge identity.  The final assertion is
the lifting criterion in \cref{prop:lift} applied to the gauge transformation.
\end{proof}

\begin{corollary}[Trace-free action identity]\label{cor:tracefreeaction}
For every $A\in\M^{\KM}_\beta(E,L;a_{P,\mathrm{par}})$,
\begin{equation}\label{eq:tracefreeaction}
 \frac1{8\pi^2}\int_X|F_A^\perp|^2dV_\omega
 =\kappa_\beta(E,L):=k(E)+2\beta\;l(E,L)-\beta^2D^2.
\end{equation}
\end{corollary}

\begin{proof}
Combine \cref{prop:CW,lem:labels} with the curvature identity \Cref{prop:detproj} \eqref{eq:curvdetproj} and the normalization
$\rho^*\tr_3=\tr_{\CC^2}$ fixed in \Cref{sec:weighted}.
\end{proof}

\begin{lemma}[Weighted determinant-one gauges preserve the divisor line]
\label{lem:gaugeblock}
Every $u\in\G_\beta^{2,p}(E,L)$ extends continuously to $D$.  In the ordered
splitting $E|_N=L\oplus Q$, the off-diagonal entries of $u$ vanish on $D$.  Consequently
$u|_D$ is diagonal and preserves both $L$ and $Q$.
\end{lemma}

\begin{proof}
With respect to the splitting $E=L\oplus Q$, locally the singular part of the trace-free model connection is
\begin{align*}
B_\beta^{\mathrm{sing}}
=
\sqrt{-1}
\begin{pmatrix}
\beta&0\\
0&-\beta
\end{pmatrix}d\theta .
\end{align*}
The scalar determinant part of $A_0$ acts trivially on $\End(E)$ by
commutator.  Hence the connection induced by $A_0$ on $\End(E)$ has the same
singular part as the connection induced by $B_\beta$.

Write
\begin{align*}
u=
\begin{pmatrix}
u_{11}&0\\
0&u_{22}
\end{pmatrix}
+\begin{pmatrix}
0&u_{12}\\
u_{21}&0
\end{pmatrix}=:u_{\diag}+u_{\mathrm{off}},
\end{align*}
then the induced covariant derivative has the form
\begin{align*}
\nabla^{\End}_{A_0}u
=
du+
\left[
\sqrt{-1}
\begin{pmatrix}
\beta&0\\
0&-\beta
\end{pmatrix}d\theta,u
\right]
+\left[\begin{pmatrix}
    \gamma_1&\\&\gamma_2
\end{pmatrix},u\right] ,
\end{align*}
where $\gamma_i$ are smooth one-forms near
$D$.  Thus
\begin{align*}
    (\nabla^{\End}_{A_0}u)_{\diag}=du_{\diag},\qquad
(\nabla^{\End}_{A_0}u)_{\mathrm{off}}=du_{\mathrm{off}}\pm2\sqrt{-1}\beta\,u_{\mathrm{off}}\,d\theta\pm(\gamma_1-\gamma_2)u_{\mathrm{off}}.
\end{align*}

Since $0<2\beta<1$, the two off-diagonal scalar connections have nonintegral monodromy parameter.  By \cite[Proposition~3.7]{KM93} or Proof of \Cref{prop:localweighted}, we have
\begin{align*}
u_{\rm diag}\in W^{2,p},
\qquad
u_{\rm off}\in W_D^{2,p}.
\end{align*}
The arguments in \cref{prop:localweighted} gives a continuous extension of $u$ across $D$ and $u_{\rm off}|_D=0$. Hence $u|_D$ is diagonal.
\end{proof}

\begin{proposition}\label{prop:banach}
The space $\A_\beta^{1,p}(E,L;a_{P,\mathrm{par}})$ is a separable affine Banach space.  The group $\G_\beta^{2,p}(E,L)$ is a Banach Lie group, modeled on $W_{B_\beta}^{2,p}(\mathfrak{su}(E))$, and acts smoothly on $\A_\beta^{1,p}(E,L;a_{P,\mathrm{par}})$.  The curvature map
\begin{align*}
    \A_\beta^{1,p}(E,L;a_{P,\mathrm{par}})&\longrightarrow L^p\bigl(\Lambda^2T^*X\otimes\mathfrak u(E)\bigr)\\
    A&\longmapsto F_A
\end{align*}
is continuous.

If $B_\beta$ and $B'_\beta$ are two model connections associated with the same
$(PU(E),R_L,\beta)$ as in Definition~3.9, then the corresponding weighted
Sobolev norms are equivalent.  In particular, they define the same spaces
$\A_\beta^{1,p}(E,L;a_{P,\mathrm{par}})$ and
$\G_\beta^{2,p}(E,L)$ with equivalent Banach topologies.
\end{proposition}

\begin{proof}
Locally we write $A_0=\frac12a_{P,\mathrm{par}}\,\Id+A_0^\perp$.
The trace-free part $A_0^\perp$ is an $SU(2)$-connection of the form
\begin{align*}
A_0^\perp=\sqrt{-1}\beta
\begin{pmatrix}
1&0\\
0&-1
\end{pmatrix}\eta+\Gamma_0,
\end{align*}
where $\Gamma_0$ is smooth and preserves the splitting $L\oplus Q$ near $D$. This is precisely the local $SU(2)$ model $A_\alpha$ of \cite{KM93} with $\alpha=\beta$.

Moreover, the scalar term $\frac12a_{P,\mathrm{par}}\Id$ acts trivially by
commutator. Hence the covariant derivatives induced by $A_0$ on $\mathfrak{su}(E)$ and on $\End(E)$ coincide with those induced by $A_0^\perp$. Therefore the spaces in \Cref{def:gaugemoduli} are locally precisely the corresponding covariant Sobolev spaces appearing in
\cite{KM93}. Consequently, \cite[Proposition~3.7]{KM93} gives the required local weighted
Sobolev description, and \cite[Lemma~3.8]{KM93} gives the multiplication maps
needed for the gauge action and curvature map.  The proof of
\cite[Proposition~2.4(ii)]{KM93} therefore applies without change to the
trace-free part, giving the Banach Lie group structure and the smooth gauge
action. The same multiplication estimate show that $A\mapsto F_A$ is continuous.

Finally, if $B_\beta$ and $B'_\beta$ are two model connections associated
with the same $(PU(E),R_L,\beta)$, then \Cref{prop:weighted-norm-indep-model} shows their weighted Sobolev norms are equivalent.  The same
holds for the induced connections on $\End(E)$, which proves the last
assertion.
\end{proof}

\subsection{Topology of bundle pairs and weighted regularity}\label{subsec:pairtopology}
\begin{definition}\label{def:smoothpair}
A smooth \textbf{determinant-fixed bundle pair} is a triple $(F,M,\phi)$, where $F$ is a smooth Hermitian rank-two bundle on $\X$, $M\subset F|_D$ is a smooth complex line bundle, and $\phi:\det F\to P$ is a smooth unitary isomorphism. Two pairs are isomorphic if there exists an isomorphism $f:F\to F'$ such that $f(M)=M'$ and $\phi=\phi'\circ\det f$.
\end{definition}

\begin{lemma}[Lifting obstruction]\label{lem:comparisonlift}
Let $(F,M,\phi_F)$ and $(E,L,\phi_E)$ be determinant-fixed bundle pairs, and let $V\subset X$ be an open set.  Suppose
\begin{align*}
\widehat u:PU(F)|_V\longrightarrow PU(E)|_V
\end{align*}
is a continuous projective bundle isomorphism such that $\widehat u(M_x)=L_x$ for any $x\in V\cap D$.

Define a $\ZZ_2$-bundle $\mathcal{L}_{\widehat u}$ on $V$ such that for $x\in V$,
\begin{align*}
\mathcal L_{\widehat u,x}
:=
\left\{
u_x:F_x\to E_x:
u_x\text{ is unitary},\; PU(u_x)=\widehat u_x,\; \phi_E\circ\det u_x=\phi_F
\right\}.
\end{align*}
Let $o(\widehat u):=[\mathcal L_{\widehat u}]\in H^1(V;\ZZ_2)$ be its isomorphism class.  

Then there exists a pair isomorphism
\begin{align*}
u:(F,M,\phi_F)|_V\longrightarrow(E,L,\phi_E)|_V
\end{align*}
satisfying $PU(u)=\widehat u$ if and only if $o(\widehat u)=0$.

Let $\xi\in H^1(V;\ZZ_2)$ and let $L_\xi\to V$ be the Hermitian line bundle induced by $\xi$.  Thenb the obstruction to lifting $\widehat u$ to $F|_V\otimes L_\xi\to E|_V$ is $o(\widehat u)+\xi$ , under the identification $PU(F|_V\otimes L_\xi)\simeq PU(F)|_V$.

Finally, let $(A_0^F,B_\beta^F)$ and $(A_0^E,B_\beta^E)$ be the singular model connections for the two bundle pairs.  If 
\begin{align*}
\widehat u\in
W^{2,p}_{B_{\beta}^F,B_\beta^E}\!\big(V;\Hom(\mathfrak{pu}_{F},\mathfrak{pu}_{E})\big),
\end{align*}
then every continuous lift $u$ satisfies
\begin{align*}
u\in W^{2,p}_{A_0^F,A_0^E}\!\big(
V;\Hom(F,E)
\big).
\end{align*}
\end{lemma}

\begin{proof}
The homomorphism
\begin{align*}
U(2)\longrightarrow PU(2)\times U(1),
\qquad
g\longmapsto([g],\det g),
\end{align*}
is a two-fold covering with kernel $\{\pm I\}$.  It follows that fiberwise
$\mathcal L_{\widehat u,x}$ consists of two elements, interchanged by
multiplication by $-I$, and that local sections of
$\mathcal L_{\widehat u}\to V$ exist.  Hence
$\mathcal L_{\widehat u}$ is a principal $\ZZ_2$-bundle.

Choose a good cover $\{V_i\}$ of $V$ and local sections $u_i\in\Gamma(V_i,\mathcal L_{\widehat u})$. On $V_i\cap V_j$ there are locally constant functions $\epsilon_{ij}:V_i\cap V_j\longrightarrow\{\pm1\}$ such that $u_j=\epsilon_{ij}u_i$. The relations $u_k=\epsilon_{jk}\epsilon_{ij}u_i$ imply $\epsilon_{ij}\epsilon_{jk}\epsilon_{ki}=1$, so $\{\epsilon_{ij}\}$ represents $o(\widehat u)$.
The local lifts can be modified by signs so as to agree on all overlaps if and
only if this cocycle is a coboundary.  Thus a global determinant-preserving
unitary lift exists exactly when $o(\widehat u)=0$.

If such a lift $u$ exists, then for $x\in V\cap D$ its projectivization sends $M_x$ to $L_x$, hence $u(M_x)=L_x$. By the definition of $\mathcal{L}_{\widehat u}$, we have $\phi_E\circ\det u=\phi_F$, so $u$ is a pair isomorphism.

Now let $\{\delta_{ij}\}$ be a $\{\pm1\}$-valued cocycle representing $\xi$. After replacing $F|_V$ by $F|_V\otimes L_\xi$, the local lifting cocycle becomes $\{\epsilon_{ij}\delta_{ij}\}$. Its cohomology class is $[\epsilon_{ij}\delta_{ij}]=o(\widehat u)+\xi$.

For the last assertion, the map from determinant-preserving unitary isomorphisms to projective isomorphisms is locally the covering $SU(2)\longrightarrow PU(2)$. Locally we can take a smooth lift $s:PU(2)\to SU(2)$ such that $u=s\circ \widehat u$.  Since $A_0^F$ and $A_0^E$ project to $B_\beta^F$ and $B_\beta^E$, the corresponding covariant derivatives satisfy
the usual chain rules. In particular, locally,
\begin{align*}
|\nabla u|=|Ds(\widehat u)\cdot\nabla\widehat u|
\le C|\nabla\widehat u|,
\end{align*}
and
\begin{align*}
|\nabla^2u|=|Ds(\widehat u)\cdot \nabla^2u+D^2s(\widehat u)(\nabla\widehat u,\nabla\widehat u)|\le C\left(
|\nabla^2\widehat u|
+
|\nabla\widehat u|^2
\right).
\end{align*}
The quadratic term belongs to $L^p$ by
\eqref{eq:weightedproducts}.  Hence $u\in W^{2,p}$.
\end{proof}

\begin{lemma}[Classification of determinant-fixed bundle pairs]\label{lem:pairclassification}
The isomorphism class of smooth determinant-fixed bundle pairs $(F,M,\phi)$ is determined by
\begin{align*}
c_2(F)\in H^4(\X;\ZZ)
\qquad\text{and}\qquad
\deg M\in\ZZ.
\end{align*}
\end{lemma}

\begin{proof}
Let $(F_0,M_0,\phi_0)$ and $(F_1,M_1,\phi_1)$ be smooth determinant-fixed bundle pairs satisfying $c_2(F_0)=c_2(F_1)$ and $\deg M_0=\deg M_1$. We need to construct a unitary bundle isomorphism $u:F_0\xrightarrow{\sim} F_1$ such that $\phi_1\circ\det u=\phi_0$ and $u(M_0)=M_1$.

We first construct a smooth fiber bundle $\mathcal I\longrightarrow \X$ such that for each $x\in \X$,
\begin{align*}
\mathcal I_x:=\left\{
u_x:F_{0,x}\longrightarrow F_{1,x}:
u_x\text{ is unitary and }\phi_{1,x}\circ\det u_x=\phi_{0,x}
\right\}.
\end{align*}
Each fibre of $\mathcal I$ is $SU(2)\simeq S^3$. A section of $\mathcal I$ is precisely a unitary
isomorphism $u:F_0\xrightarrow{\sim} F_1$ satisfying $\phi_1\circ \det u=\phi_0$.

 Since $\pi_0(S^3)=\pi_1(S^3)=\pi_2(S^3)=0$, a section of $\mathcal I$ can be constructed over the three-skeleton $\X^{(3)}$ of $\X$.  The only possible obstruction $o(\mathcal{I})$ to extending it over $\X$ lies in $H^4(\X;\pi_3(S^3))=H^4(\X;\ZZ)$.
This is the standard obstruction to extending a section of a fibre bundle. The obstruction theory gives
\begin{align*}
o(\mathcal I)=\pm\bigl(c_2(F_1)-c_2(F_0)\bigr)\in H^4(\X;\ZZ).
\end{align*}
By hypothesis this class vanishes, therefore we have a global continuous section of $\mathcal I$, and hence a continuous unitary bundle isomorphism
\begin{align*}
f:F_0\xrightarrow{\sim} F_1,\qquad \phi_1\circ\det f=\phi_0.
\end{align*}
Moreover, we may take $f$ to be smooth by approximation.

It remains to identify the line subbundles. Put $M'_0:=f(M_0)\subset F_1|_D$. Since $f$ is a bundle isomorphism, $\deg M'_0=\deg M_0=\deg M_1$.
Thus it is enough to assume $F_0=F_1=:F$ and construct a determinant-one unitary automorphism carrying $M_0$ to
$M_1$.

Complex line bundles over the closed oriented surface $D$ are classified by
their first Chern class. Hence there is a unitary isomorphism $f_M:M_0\xrightarrow{\sim}M_1$. Take 
the orthogonal decompositions $F|_D=M_i\oplus Q_i$, then $\deg Q_0=\deg Q_1$. So we also have a unitary isomorphism $f_Q:Q_0\xrightarrow{\sim}Q_1$.
Thus we get a unitary isomorphism over $D$
\begin{align*}
f_D:=f_M\oplus f_Q:F|_D\longrightarrow F|_D,\qquad f_D(M_0)=M_1.
\end{align*}
Moreover, we can assume $\det f_D\equiv 1$ after multiplying a smooth function.

Since $D$ has real dimension two and $\pi_1(SU(2))=\pi_2(SU(2))=0$, $f_D\in SU(F|_D)$ is homotopic to the identity through $SU(F|_D)$.  Let $U_t$ be such a homotopy, with $U_0=\Id$ and $U_1=f_D$. Choose a tubular neighborhood $N$ of $D$, identify $F|_N$ with the pullback
of $F|_D$. Let $\chi$ be a smooth radial cutoff equal to $1$ near $D$ and $0$ near $\partial N$.  Then
\begin{align*}
x\longmapsto U_{\chi(r(x))}(\pi(x))
\end{align*}
defines an isomorphism $SU(F|_N)$, equal to $f_D$ on $D$ and to the identity near $\partial N$. Extending it by the identity outside $N$ gives an isomorphism $SU(F)$ carrying $M_0$ to $M_1$.
\end{proof}
\begin{lemma}\label{lem:pairsmoothing}
Let $C:(F,M,\phi_F)\to(E,L,\phi_E)$ be a continuous determinant-fixed bundle pair isomorphism. Then $C$ is homotopic through pair isomorphisms to a smooth one.
\end{lemma}

\begin{proof}
By the Whitney approximation theorem, there exists a smooth section $C':F\to E$ of $\Hom(F,E)$ which is uniformly arbitrarily close to $C$. 

We first modify $C'$ so that it preserves the line subbundle. Let $p_M,p_L$ be the orthogonal projections of $F|_D,E|_D$ to the line subbundle, respectively. On $D$ put
\begin{align*}
    R:=(\Id-P_L)\,C'\,P_M\in\Hom(F|_D,E|_D).
\end{align*}
Since $C(M)=L$, the section $R$ is arbitrarily small. Extend $R$ smoothly by a smooth cutoff function, and denote the resulting section of $\Hom(F,E)$ by $\widetilde R$. Replacing $C'$ by $C_1':=C'-\widetilde R$. we obtain a smooth bundle isomorphism, still arbitrarily close to $C$, such that $C_1'(M)=L$. Moreover, the maps
\begin{align*}
    C_s':=(1-s)C+sC_1'
\end{align*}
are invertible for all $s\in[0,1]$ and satisfies $C_s'(M)=L$, after choosing the approximation sufficiently close.

It remains to impose the determinant condition.  Define the smooth functions $d_s:\X\longrightarrow\CC^*$ by
\begin{align*}
\phi_E\circ\det C_s'=d_s\,\phi_F.
\end{align*}
Since $d_0\equiv1$ and $C_1'$ may be chosen uniformly arbitrarily close to $C$, all $d_s$ take values in a small neighborhood of $1\in\CC^*$. Hence we can define a family of smooth isomorphism
\begin{align*}
C_s:=\exp\!\left(-\frac12\log d_s\right)C_s'.
\end{align*}
such that
\begin{align*}
\phi_E\circ\det C_s=\phi_F,\quad C_s(M)=L,\quad C_0=C.
\end{align*}
\end{proof}

\begin{proposition}[Regularity and extension]\label{prop:gauge}
Let $A_i=A_{0,i}+a_i\in\A_\beta^{1,p}(E_i,L_i;a_{P,\mathrm{par}})$ $(i=1,2)$
be weighted unitary connections on bundle pairs $(E_i,L_i,\phi_i)$.  Let $S\subset \X$ be finite.

If a determinant-preserving unitary bundle isomorphism $u\in W^{1,1}_{\mathrm{loc}}\big(X\setminus S;\Hom(E_2,E_1)\big)$ satisfies
\begin{align*}
\nabla_{A_1}\circ u=u\circ\nabla_{A_2}.
\end{align*}
Then $u$ extends uniquely across $S$ and $D$ to a continuous bundle pair isomorphism with weighted regularity
\begin{align*}
u(L_2)=L_1,\qquad u\in W^{2,p}_{A_0}(\X;\Hom(E_2,E_1)).
\end{align*}

Moreover, for any $U_0\subset\subset U\subset \X$, there exists a constant $C$ such that
\begin{align*}
\|u\|_{W^{2,p}_{A_0}(U_0)}
\le
C\Bigl(1+\|a_1\|_{W^{1,p}_{A_{0,1}}(U)}+\|a_2\|_{W^{1,p}_{A_{0,2}}(U)}
\Bigr)^2.
\end{align*}
\end{proposition}

\begin{proof}
Let $\nabla_0$ be the connection on $\Hom(E_2,E_1)$ induced by $(E_i,A_{0,i})$, i.e. 
\begin{align*}
\nabla_0u:=\nabla_{A_{0,1}}\circ u-u\circ\nabla_{A_{0,2}}.
\end{align*}
Then the equation is equivalent to
\begin{align}\label{eq:intertwining}
    \nabla_0u=ua_2-a_1u\qquad \text{on }X\setminus S.
\end{align}

By \cref{prop:localweighted} and the Sobolev embedding $W^{1,p}\hookrightarrow L^{2p}$ in dimension four, we have $a_i\in L^{2p}(X)$. Since $u$ is unitary,
\begin{equation}\label{eq:gaugefirst}
\|\nabla_0u\|_{L^{2p}(X\setminus S)}
\le C\bigl(\|a_1\|_{L^{2p}}+\|a_2\|_{L^{2p}}\bigr).
\end{equation}

We show that \eqref{eq:intertwining} holds distributionally across $D\cup S$. Let $\Psi\in C^\infty(T^*\X\otimes\Hom(E_2,E_1))$ be a test section. Take a cutoff function $\chi_\varepsilon\in C^\infty_c(X\setminus S)$ such that $\chi_\varepsilon=0$ on the $\varepsilon$-neighborhood of $D\cup S$, $\chi_\varepsilon=1$ outside the $2\varepsilon$-neighborhood, and $|d\chi_\varepsilon|\le C\varepsilon^{-1}$. We therefore have $\|d\chi_\varepsilon\|_{L^1}\le C\varepsilon$.
Testing $\chi_\varepsilon\Psi$ to \eqref{eq:intertwining} gives
\begin{equation}\label{eq:cutoffdistribution}
    \begin{aligned}
        \int_{\X}\bigl\langle ua_2-a_1u,\chi_\varepsilon\Psi\bigr\rangle=\int_{\X}
\bigl\langle u,\nabla_0^*(\chi_\varepsilon\Psi)\bigr\rangle
=\int_{\X}\bigl\langle u,\chi_\varepsilon\nabla_0^*\Psi\bigr\rangle-\int_{\X}\bigl\langle u,i_{d\chi_\varepsilon}\Psi\bigr\rangle.
    \end{aligned}
\end{equation}
$\|d\chi_\varepsilon\|_{L^1}\to 0$ implies $\int_{\X}\bigl\langle u,i_{d\chi_\varepsilon}\Psi\bigr\rangle\to 0$ as $\varepsilon\to 0$. The singular terms of the singular connections $A_{0,i}$ are $O(r^{-1})$, where $r$ is the radial function to the divisor $D$. Thus 
\begin{align*}
    \left|\int_N\bigl\langle u,\chi_\varepsilon\nabla_0^*\Psi\bigr\rangle\right|\le C\int_0^1r^{-1}\;rdr<\infty.
\end{align*}
Take $\varepsilon\to 0$ to \eqref{eq:cutoffdistribution}, then \eqref{eq:intertwining} holds on all of $\X$ in the distributional sense by dominated convergence and $ua_2-a_1u\in L^{2p}(\X)$. Since $\nabla_0$ is smooth away from $D$,
\begin{align*}
    u\in W^{1,2p}_{\mathrm{loc}}(X),\qquad
\nabla_0u=ua_2-a_1u \quad\text{on }X.
\end{align*}

We differentiate \eqref{eq:intertwining} weakly to obtain
\begin{align*}
    \nabla_0^2u=(\nabla_0u)\otimes a_2+u\nabla_{A_{0,2}}a_2-(\nabla_{A_{0,1}}a_1)u
-a_1\otimes\nabla_0u.
\end{align*}
Here $a_i,\nabla_0u\in L^{2p}$,
$\nabla_{A_{0,i}}a_i\in L^p$, and $u\in L^\infty$.
Hölder's inequality therefore gives $\nabla_0^2u\in L^p(X)$. Together with \eqref{eq:gaugefirst} and the boundedness of $u$, this proves
\begin{align*}
    u\in W^{2,p}_{A_0}(X;\Hom(E_2,E_1)).
\end{align*}

We now work near $D$ in the unitary frames of $E_i$ adapted to the splitting $E_i|_{N}=L_i\oplus Q_i$. We say the $\Hom(L_2,L_1)\oplus\Hom(Q_2,Q_1)$ component to be the diagonal part $u_{\diag}$, and the $\Hom(L_2,Q_1)\oplus\Hom(Q_2,L_1)$ component to be the off-diagonal part $u_{\mathrm{off}}$. As the computation in Proof of \Cref{lem:gaugeblock}, we have
\begin{align*}
(\nabla_0u)_{\diag}&=\nabla_{\mathrm{sm}}u_{\mathrm{diag}},\\    (\nabla_0u)_{\mathrm{off}}&=\left(\nabla_{\mathrm{sm}}\pm2\ii\beta d\theta\right) u_{\mathrm{off}},
\end{align*}
where $\nabla_{\mathrm{sm}}$ is a smooth connection. Then $u\in W^{2,p}_{A_0}$ and \Cref{prop:localweighted} gives $u_{\mathrm{diag}}\in W^{2,p}$, $u_{\mathrm{off}}\in W_D^{2,p}$ and give a continuous extension
of $u$ across $D$, with $u_{\mathrm{off}}|_D=0$.
Hence $u(L_2)=L_1$.
Unitarity and the determinant condition extend
by continuity, so the extension is a bundle pair
isomorphism. Uniqueness follows from density.

Finally, let $U_0\subset\subset U\subset\X$ and set
$$
M_U:=\|a_1\|_{W^{1,p}_{A_{0,1}}(U)}
+\|a_2\|_{W^{1,p}_{A_{0,2}}(U)}.
$$
The local Sobolev embedding and
\eqref{eq:gaugefirst} give
$$
\|\nabla_0u\|_{L^{2p}(U_0)}\le CM_U.
$$
The differentiated equation then gives
$$
\|\nabla_0^2u\|_{L^p(U_0)}
\le C(M_U+M_U^2).
$$
Since $u$ is unitary, its $L^p$ norm is bounded.
Combining these estimates yields
$$
\|u\|_{W^{2,p}_{A_0}(U_0)}
\le C(1+M_U)^2.
$$
\end{proof}

\subsection{Lifting the projective limit}

The projective compactness theorem has been established in \Cref{sec:weighted}. We now lift its limit to the prescribed-determinant rank-two problem. The difference between $PU(2)$ case and $U(2)$ case is to keep track of the $H^1(X;\ZZ_2)$ obstruction.

\begin{theorem}\label{thm:detcompact}
Let $[A_n]\in\mathcal M_\beta^{\mathrm{KM}}(E,L;a_{P,\mathrm{par}})$. Then after passing to a subsequence, there exist
\begin{enumerate}
\item finite set $S:=S_X\sqcup S_D$, where $S_X\subset X$ and $S_D\subset D$;
\item a smooth determinant-fixed bundle pair $(E_\infty,L_\infty,\phi_\infty)$;
\item a connection $A_\infty\in\mathcal A_\beta^{1,p}
(E_\infty,L_\infty;a_{P,\mathrm{par}})$ satisfying $(F_{A_\infty}^{\perp})^+=0$;
\item determinant-preserving bundle-pair isomorphisms
\begin{align*}
    \Psi_n:(E_\infty,L_\infty,\phi_\infty)|_{\X\setminus S}\to(E,L,\phi_E)|_{\X\setminus S};
\end{align*}
\item integers $m_x\in\ZZ_{>0}$ for $x\in S_X$, and $(k_y,l_y)\in\ZZ^2\setminus\{(0,0)\}$ for $y\in S_D$, satisfying $k_y\geq0$, $k_y+l_y\geq0$.
\end{enumerate}
such that $\Psi_n^*A_n\to A_\infty$ in weighted \(W^{1,p}\) on every relatively compact subsets of \(\X\setminus S\).

Moreover,
\begin{align*}
c_2(E_\infty)&=c_2(E)-\sum_{x\in S_X}m_x-\sum_{y\in S_D}k_y,\\
\deg L_\infty&=\deg L+\sum_{y\in S_D}l_y,\\
|F_{A_n}^{\perp}|^2\,dV_\omega&\rightharpoonup|F_{A_\infty}^{\perp}|^2\,dV_\omega+8\pi^2\sum_{x\in S_X}m_x\delta_x+8\pi^2\sum_{y\in S_D}(k_y+2\beta l_y)\delta_y.
\end{align*}
\end{theorem}

\begin{proof}
Choose representatives $A_n\in \mathcal A_\beta^{1,p}
(E,L;a_{P,\mathrm{par}})$ of the given classes and set $B_n:=A_n^{\mathrm{proj}}$. By \Cref{prop:projectivizationweighted}, $B_n\in\widehat{\mathcal A}_\beta^{1,p}(\PU(E),R_L)$, $F_{B_n}=F_{A_n}^{\perp}$ and $F_{B_n}^+=0$. Then we can apply \Cref{thm:puKM}. After passing to a subsequence, there exist
a finite set \(S\subset\X\), a projective bundle pair
\((\widehat P_\infty,R_\infty)\), an anti-self-dual connection $B_\infty\in\widehat{\mathcal A}_\beta^{1,p}(\widehat P_\infty,R_\infty)$, and projective pair isomorphisms
\begin{align*}
\widehat\Psi_n:
(\widehat P_\infty,R_\infty)|_{\X\setminus S}
\longrightarrow
(\PU(E),R_L)|_{\X\setminus S}
\end{align*}
with the convergence asserted in \Cref{thm:puKM}. Let
\((k_z,l_z)\) denote its relative defect labels.

The projective conservation law gives
\begin{align*}
w_2(\widehat P_\infty)
=
w_2(\PU(E))
=
c_1(P)\pmod 2.
\end{align*}
Hence \Cref{prop:lift} gives a Hermitian rank-two bundle
\(F\to\X\), together with unitary identifications
\begin{align*}
\phi_F:\det F\xrightarrow{\sim}P,
\qquad
\Phi_F:\PU(F)\xrightarrow{\sim}\widehat P_\infty.
\end{align*}
By \Cref{prop:line-reduction}, the oriented reduction
\(R_\infty\) determines a line subbundle $M\subset F|_D$.

It remains to arrange determinant-preserving lifts of the comparison maps. For each \(n\), let $o_n:=o(\widehat\Psi_n)\in H^1(\X\setminus S;\ZZ_2)$
be the obstruction of \Cref{lem:comparisonlift}. Since \(S\) is finite and \(\X\) has real dimension four, restriction induces an isomorphism $H^1(\X;\ZZ_2)\simeq H^1(\X\setminus S;\ZZ_2)$. Since the group \(H^1(\X;\ZZ_2)\) is finite, after passing to a further subsequence, there exists $\xi\in H^1(\X;\ZZ_2)$ such that \(o_n=\xi\) for every \(n\).

Let \(L_\xi\to\X\) be the flat line bundle associated with \(\xi\), with a fixed trivialization $L_\xi^{\otimes2}\simeq\underline{\CC}$. Replacing $(F,M,\phi_F)$ by $(F\otimes L_\xi, M\otimes L_\xi|_D, \phi_F\otimes\Id_{L_\xi^{\otimes2}})$ doesn't change the projectivization and oriented reduction \((\widehat P_\infty,R_\infty)\). By
\Cref{lem:comparisonlift}, the obstruction of each
\(\widehat\Psi_n\) is now $o_n+\xi=0$.
Thus \(\widehat\Psi_n\) admits a determinant-preserving lift
\begin{align*}
\Psi_n:
(F,M,\phi_F)|_{\X\setminus S}
\longrightarrow
(E,L,\phi_E)|_{\X\setminus S}.
\end{align*}
Set $A_n'=\Psi_n^*A_n$ on $\X\setminus(S\cup D)$.
Since $\det A_n'=a_{P,\mathrm{par}}$ and
$(A_n')^{\mathrm{proj}}=\widehat\Psi_n^*B_n$,
therefore $A_n'$ has locally weighted $W^{1,p}$ regularity. Applying \Cref{prop:gauge} to $\Psi_n$ gives its weighted $W^{2,p}$ regularity.

We next construct the limiting rank-two connection. Transport
\(B_\infty\) to \(\PU(F)\) by \(\Phi_F\). By
\Cref{prop:detproj}, there exists a unique unitary connection
\(A_\infty\) on \(F|_X\) satisfying
\begin{align*}
\det A_\infty=a_{P,\mathrm{par}},
\qquad
A_\infty^{\mathrm{proj}}=B_\infty.
\end{align*}
Let \(A_{0,\infty}\) be the determinant-fixed model connection associated with \((F,M)\), and its projective part is the model \(B_{\beta,\infty}\) associated with \((\widehat P_\infty,R_\infty)\). Since $B_\infty-B_{\beta,\infty}\in W^{1,p}_{B_{\beta,\infty}}$, and \(A_\infty-A_{0,\infty}\) is its trace-free lift under
\(\mathfrak{pu}(2)\simeq\mathfrak{su}(2)\), one has $A_\infty-A_{0,\infty}
\in
W^{1,p}_{B_{\beta,\infty}}
(T^*X\otimes\mathfrak{su}(F))$. Therefore $A_\infty\in\mathcal A_\beta^{1,p}
(F,M;a_{P,\mathrm{par}})$.
Moreover, $F_{A_\infty}^{\perp}=F_{B_\infty}$, so $(F_{A_\infty}^{\perp})^+=0$.
Since every \(\Psi_n\) is determinant-preserving, $\Psi_n^*A_n-A_\infty=\widehat\Psi_n^*B_n-B_\infty$.
The convergence in \Cref{thm:puKM} therefore gives the asserted strong weighted \(W^{1,p}\) convergence and the smooth convergence away from \(D\cup S\).

Set $S_X:=S\setminus D$, $S_D:=S\cap D$. For \(x\in S_X\), write $m_x:=k_x>0$.
For \(y\in S_D\), retain the notation \((k_y,l_y)\). \Cref{prop:puadmissible} gives $k_y\geq0$, $k_y+l_y\geq0$.

Since \(\det E=\det F=P\), \Cref{lem:labels} gives
\begin{align*}
k(E)-k(F)=c_2(E)-c_2(F),\qquad l(E,L)-l(F,M)=\deg M-\deg L.
\end{align*}
\Cref{thm:puKM} therefore yields
\begin{align*}
    c_2(F)&=c_2(E)-\sum_{x\in S_X}m_x-\sum_{y\in S_D}k_y,\qquad\deg M=\deg L+\sum_{y\in S_D}l_y.
\end{align*}
Set $(E_\infty,L_\infty,\phi_\infty):=(F,M,\phi_F)$.

Finally, \(F_{B_n}=F_{A_n}^{\perp}\) and
\(F_{B_\infty}=F_{A_\infty}^{\perp}\), so the measure convergence in \Cref{thm:puKM} gives the stated convergence.
\end{proof}

\section{From parabolic Hermitian--Einstein to twisted ASD connection}\label{sec:HE-to-KM}
In this section, we aim to prove one direction: given a parabolic Hermitian--Einstein connection, it can be naturally mapped into the space of determinant-fixed weighted unitary connections $\A_\beta^{1,p}(E,L;a_{P,\mathrm{par}})$.
For this goal, we show that the curvature of any compatible parabolic Hermitian--Einstein metric belongs to $L^p$ for some $p>2$. 

\subsection{Finite energy from compatibility}\label{sec:finite}
In this subsection, we first show the curvature belongs to $L^2$.
\begin{theorem}[c.f. \Cref{thm:energy}]
For any smooth Hermitian--Einstein metric $H$ on $E|_X$ compatible with $\mathcal{E}_*$ (in the sense of \Cref{def:compatibility}), one has $F_H\in L^2(X,H,\omega)$.
\end{theorem}

It suffices to consider locally in a divisor chart. On a holomorphic bidisc $U=\Delta_z\times\Delta_w$ with $D\cap U=\{z=0\}$, write $r=|z|$ and use a holomorphic frame $(e_1,e_2)$ such that $e_1|_D$ spans $L\subset E|_D$. In this frame, the smooth metric $H$ is a matrix-valued function on $U\setminus D$.  Put $K=\diag(r^{2\alpha_2},r^{2\alpha_1})$ to be the local model metric.

\begin{lemma}\label{lem:model}
The compatibility of $H$ implies
$$
 C_0^{-1}K\leq H\leq C_0K,\qquad \nabla_K(K^{-1}H)\in L^2(K).
$$
Equivalently, in the multivalued holomorphic $K$-unitary frame $(z^{-\alpha_2}e_1,z^{-\alpha_1}e_2)$, the matrix representative $u$ of the metric $H$ satisfies
\begin{align}\label{eq:u-du-bound}
     C_0^{-1}I\leq u\leq C_0I,\qquad du\in L^2.
\end{align}
\end{lemma}

\begin{proof}
    By \cref{lem:parmodelindependence}, compatibility is unchanged after replacing the chosen model metric by another model metric $K$.  Thus $b:=K^{-1}H$ and $b^{-1}$ are bounded and $\bar\partial_{\End\E}b\in L^2(K)$. Since $b$ is $K$-self-adjoint, the Chern connection satisfies $\nabla^{1,0}_Kb=(\nabla^{0,1}_K b)^{*K}=(\bar\partial_{\End\E}b)^{*K}$, so $\nabla_Kb\in L^2(K)$.
\end{proof}
\begin{remark}
    The multivalued frame has constant unitary monodromy, so all norms and trace identities are single-valued.
\end{remark}

Now we have a smooth function $u\in C^\infty(U\setminus D;\mathrm{Herm}_+)$ satisfying \eqref{eq:u-du-bound} and the Hermitian--Einstein equation $\ii\Lambda_\omega\db(u^{-1}\dd u)=\lambda \Id$. It remains to prove $F_u:=\db(u^{-1}\dd u)\in L^p(U\setminus D)$.

Put $B_u:=u^{-1}\dd u\in L^2$, then $F_u=\db B_u$. We first show $F_u\in L^2$ by differentiating the Hermitian--Einstein equation and integration by part. However, the integrability of our setting is not good enough, so we need a difference quotient trick to get better regularity.

\begin{lemma}\label{lem:logdifference}
Suppose $u$ be as above. Let $Y\in\mathrm{span}_\RR\{\dd_x,\dd_y\}$ be a real tangent vector field, where $w=x+\ii y$. 
For sufficiently small real number $h\ne0$ and $\tau\in[0,1]$, define
$$u_h(x):=u(x+hY),\quad s_h:=\frac1h\log(u^{-1}u_h),\quad u_{\tau,h}:=u\exp(\tau h s_h).$$
On smaller bidiscs, the metrics $u_{\tau,h}$ and their inverses are uniformly bounded, independently of $h$ and $\tau\in[0,1]$. For fixed $h$, $s_h\in W^{1,2}\cap L^\infty$ and is self-adjoint relative to every $u_{\tau,h}$. Moreover, there exists a constant $C$ only depends on $C_0$ such that
\begin{align*}
    \|s_h\|_{L^2}&\leq C\|Yu\|_{L^2},\qquad s_h\longrightarrow u^{-1}Yu=:s_Y\quad\text{in }L^2.
\end{align*}
\end{lemma}
\begin{proof}
Consider Hermitian matrix $A_h:=u^{-1/2}u_hu^{-1/2}$, then eigenvalues of $A_h$ lie in $[C_0^{-2},C_0^2]$. Since
\begin{align*}
u_{\tau,h}=u\exp(\tau\log(u^{-1}u_h))=u\exp(\tau\log(u^{-\frac12}A_h u^{\frac12}))=u^{\frac12}A_h^\tau u^{\frac12},
\end{align*}
we have $C_0^{-3}I\le u_{\tau,h}\le C_0^3 I$.

Let $T_h:=u^{-1}u_h$. The fact $u^*=u$ and $u_h^*=u_h$ gives $T_h^*u=uT_h$. Since $s_h=\frac1h\log T_h$, we have $s_h^*u=us_h$ and $s_h^*u_{\tau,h}=u_{\tau,h}s_h$. This proves the self-adjointness. 

By assumption, we have $u,u_h\in W^{1,2}\cap L^\infty$. Since $C_0^{-1}I\le u\le C_0I$, we also have $u^{\pm\frac12}\in W^{1,2}\cap L^\infty$. Note that $W^{1,2}\cap L^\infty$ is an algebra due to $\|d(fg)\|_{L^2}\le \|f\|_{L^{\infty}}\|dg\|_{L^2}+\|df\|_{L^2}\|g\|_{L^\infty}$. Thus we have $A_h\in W^{1,2}\cap L^\infty$. Moreover, we have
\begin{align*}
    s_h=\frac1h u^{-\frac12}\log A_h u^{\frac12}\in W^{1,2}\cap L^\infty.
\end{align*}
In fact, the $L^\infty$-norm estimate is not uniform on $h$.

Define the standard differential quotient $\delta_hu:=h^{-1}(u_h-u)$. Since $C_0^{-2}I\le A_h\le C_0^2I$,
\begin{align*}
    \left|\log A_h\right|\le C |A_h-I|=C\left|u^{-\frac12}(u_h-u)u^{-\frac12}\right|\le C|u_h-u|.
\end{align*}
Thus one has 
$$|s_h|=\frac1h\left|u^{-\frac12}\log A_h u^{\frac12}\right|\le \frac{C}{h}|u_h-u|=C|\delta_hu|.$$
The standard Sobolev theory gives $\|\delta_hu\|_{L^2}\le C\|Yu\|_{L^2}$, then we have $\|s_h\|_{L^2}\le C\|Yu\|_{L^2}$.

Finally, write $s_h-u^{-1}Yu=(s_h-u^{-1}\delta_hu)+(u^{-1}\delta_hu-u^{-1}Yu)$. The standard Sobolev theory gives $\|\delta_hu-Yu\|_{L^2}\to 0$. Since $u^{-1}\in L^\infty$, the second term goes to $0$ in $L^2$. For the first term, take $X_h:=u^{-\frac12}(u_h-u)^{-\frac12}$, one has
\begin{align*}
    |s_h-u^{-1}\delta_hu|&=\frac1h\left|u^{-\frac12}\log A_h\;u^{\frac12}-u^{-1}(u_h-u)\right|\\
    &=\frac1h\left|u^{-\frac12}\big(\log(I+X_h)-X_h\big)u^{\frac12}\right|\\
    &\le \frac{C}{h}|X_h|^2\\
    &\le C|u_h-u|\cdot|\delta_hu|.
\end{align*}
Since $\delta_hu\to Yu$ in $L^2$, there exists $M$ such that $\|\delta_hu\|_{L^2}\le M$. Then
\begin{align*}
    \int_U|s_h-u^{-1}\delta_hu|^2&\le C_0\left(\int_{|u_h-u|\le \varepsilon}|u_h-u|^2\cdot|\delta_hu|^2+\int_{|u_h-u|> \varepsilon}|u_h-u|^2\cdot|\delta_hu|^2\right)\\
    &\le C_0(\varepsilon^2M+10C^2\|Yu\|_{L^2(|u_h-u|>\varepsilon)}).
\end{align*}
Since $Yu\in L^2$ and $u_h-u\to 0$ in measure, we have $s_h-u^{-1}\delta_hu\to 0$ in $L^2$.
\end{proof}
\begin{remark}
    Note that the assumption only implies $s_Y\in L^2$, while we have $s_h\in W^{1,2}\cap L^\infty$ and $s_h\to s_Y$ in $L^2$.
\end{remark}

\begin{proposition}
\label{prop:weakvariation}
In the setting of \cref{lem:logdifference}, suppose $\ii\Lambda_gF_u=\lambda\Id$ in $U\setminus D$. Let $U_0\subset\subset U_1\subset\subset U$, and $\chi$ be a smooth cutoff function satisfies $\chi|_{U_0}\equiv 1$, $\supp\chi\subset U_1$.

Then there exists a constant $C=C(C_0,\|\chi\|_{C^1},\|g\|_{C^2})$ such that
\begin{equation}\label{eq:caccioppoli}
\|\chi\bar\partial s_h\|_{L^2(U_1)}\leq
 C\bigl(\|B_{u_h}\|_{L^2(U_1)}+\|s_h\|_{L^2(U_1)}\bigr).
\end{equation}
\end{proposition}

\begin{proof}
Let $g_h(x)=g(x+hY)$, and $\nabla^{1,0}_{u_{\tau,h}}=\dd+[B_{u_{\tau,h}},\cdot]$ be $(1,0)$-part of the Chern connection of Hermitian metric $u_{\tau,h}$ induced on $\End E$. Since $\dd_\tau u_{\tau,h}=\dd_\tau\left(u\exp(\tau hs_h)\right)=hu_{\tau,h}s_h$,
\begin{align*}
    \dd_\tau B_{u_{\tau,h}}=\dd_\tau\left(u_{\tau,h}^{-1}\dd u_{\tau,h}\right)=h\left(\dd s_{h}+B_{u_{\tau,h}}s_h-s_h B_{u_{\tau,h}}\right)=h\nabla^{1,0}_{u_{\tau,h}}s_h.
\end{align*}
Integrating along $u_{1,h}=u_h$ and $u_{0,h}=u$ gives
\begin{equation}\label{eq:connectionvariation}
\frac{B_{u_h}-B_u}{h}=\int_0^1\nabla^{1,0}_{u_{\tau,h}}s_h\,d\tau.
\end{equation}
Actually this is the difference quotient version of the identity $YB_u=\nabla^{1,0}_u s_Y$.

The Hermitian--Einstein equation $\ii\Lambda_gF_u=\ii\Lambda_{g_h}F_{u_h}=\lambda\Id$ gives,
\begin{align*}
    \ii\Lambda_g\left(F_{u_h}-F_u\right)=\ii\left(\Lambda_g-\Lambda_{g_h}\right)F_{u_h}=hC_hF_{u_h},
\end{align*}
where $C_h=h^{-1}\left(\ii\Lambda_g-\ii\Lambda_{g_h}\right)$. Since the metric $g$ is smooth, the coefficients of $C_h$ have uniform $C^k$ bounds for small $h$. Then we get the weak variation identity by differentiating \eqref{eq:connectionvariation}
\begin{align}\label{eq:weak-variation-identity}
    \int_0^1\ii\Lambda_g\db\nabla^{1,0}_{u_{\tau,h}}s_h d\tau=\ii\Lambda_g\db\left(\frac{B_{u_h}-B_u}{h}\right)=C_hF_{u_h}=C_h\db B_{u_h}.
\end{align}
Take the cutoff function $\eta_\varepsilon$ constructed in \eqref{eq:logcutoff}. Recall that one has $\eta_\varepsilon=0$ for $r\leq\varepsilon^2$, $\eta_\varepsilon=1$ for $r\geq\varepsilon$, and has
$$|d\eta_\varepsilon|\leq \frac{C}{r|\log\varepsilon|},\qquad
\|d\eta_\varepsilon\|_{L^2}^2\leq\frac{C}{|\log\varepsilon|}.$$
We use the function $v:=\eta_\varepsilon^2\chi^2s_h$ to be a test function of \eqref{eq:weak-variation-identity}, then we have
\begin{align}\label{eq:testv}
    \int_0^1\int_U\re\tr\left(\ii\Lambda_g\db\nabla^{1,0}_{u_{\tau,h}}s_h \cdot v\right) dV d\tau=\int_U\re\tr\left(C_h\db B_{u_h}\cdot  v\right) dV.
\end{align}
The left hand side of \eqref{eq:testv} has
\begin{align*}
    L:=\int_0^1\int_U\re\tr\left(\ii\Lambda_g\db\nabla^{1,0}_{u_{\tau,h}}s_h \cdot v\right) dV d\tau&=-\int_0^1\int_U g^{i\bar{j}}\re\tr\left(\dd_{\bar j}\left(\nabla^{1,0}_{u_{\tau,h}}s_h\right)_i \cdot v\right) dV d\tau\\
    &=\int_0^1\int_U g^{i\bar{j}}\re\tr\left(\left(\nabla^{1,0}_{u_{\tau,h}}s_h\right)_i \cdot \dd_{\bar j}v\right)dV d\tau\\
    &\quad +\int_0^1\int_U \dd_{\bar j}g^{i\bar{j}}\re\tr\left(\left(\nabla^{1,0}_{u_{\tau,h}}s_h\right)_i \cdot v\right)dV d\tau.
\end{align*}
Let $v=\eta_\varepsilon^2\chi^2s_h$, we divide $L$ into four parts: $L=L_1+L_2+L_3+L_4$, where
\begin{align*}
    L_1&=\int_0^1\int_U \chi^2\eta^2g^{i\bar{j}}\re\tr\left(\left(\nabla^{1,0}_{u_{\tau,h}}s_h\right)_i \cdot \dd_{\bar j}s_h\right)dV d\tau,\\
    L_2&=2\int_0^1\int_U \chi^2\eta\dd_{\bar j}\eta \;g^{i\bar{j}}\re\tr\left(\left(\nabla^{1,0}_{u_{\tau,h}}s_h\right)_i \cdot s_h\right)dV d\tau,\\
    L_3&=2\int_0^1\int_U \eta^2\chi\dd_{\bar j}\chi\;  g^{i\bar{j}}\re\tr\left(\left(\nabla^{1,0}_{u_{\tau,h}}s_h\right)_i \cdot s_h\right)dV d\tau,\\
    L_4&=\int_0^1\int_U \chi^2\eta^2\dd_{\bar j}g^{i\bar{j}}\re\tr\left(\left(\nabla^{1,0}_{u_{\tau,h}}s_h\right)_i \cdot s_h\right)dV d\tau.
\end{align*}
\Cref{lem:logdifference} shows that $s_h$ is self-adjoint to $u_{\tau,h}$, i.e. $s_h^*u_{\tau,h}=u_{\tau,h}s_h$. Note that this gives $\left(\nabla^{1,0}_{u_{\tau,h}}s_h\right)_i=u_{\tau,h}^{-1}(\dd_{\bar i} s_h)^*u_{\tau,h}$. Hence $\left|\nabla^{1,0}_{u_{\tau,h}}s_h\right|\le C_0|\db s_h|$.

Since $g$ is smooth and $C^{-1}I\le u_{\tau,h}\le CI$,
\begin{align*}
    L_1&=\int_0^1\int_U \chi^2\eta^2g^{i\bar{j}}\tr\left(u_{\tau,h}^{-1}(\dd_{\bar i} s_h)^*u_{\tau,h} \dd_{\bar j}s_h\right)dV d\tau\ge C^{-1}\int_U \chi^2\eta^2\left|\db s_h\right|^2 dV=C^{-1}\|\chi\eta\db s_h\|^2_{L^2}.
\end{align*}
Since $\left|\nabla^{1,0}_{u_{\tau,h}}s_h\right|\le C|\db s_h|$ and $s_h\in W^{1,2}\cap L^\infty$ for any fixed $h$,
\begin{align*}
    |L_2|\le C\int_U |\db \eta|\cdot |\db s_h|\cdot |s_h|dV \le C \|s_h\|_{L^\infty}\|\db s_h\|_{L^2}\|\db\eta_{\varepsilon}\|_{L^2}.
\end{align*}
Recall $\|d\eta_\varepsilon\|_{L^2}^2\le C|\log \varepsilon|^{-1}$, then we have $L_2=\Psi(\varepsilon|h)$, which means for any fixed $h$, this quantity tends to $0$ as $\varepsilon\to0$.

Since $|d\chi|\le C$ and $|dg|\le C$,
\begin{align*}
    |L_3|+|L_4|\le C_0\int_U \eta\chi |\db s_h| \cdot |s_h|dV\le \delta_1\|\chi\eta\db s_h\|_{L^2}^2+C_\delta\|s_h\|_{L^2}^2.
\end{align*}

On the other hand, the right hand side of \eqref{eq:testv} gives
\begin{align*}
    R:=\int_U\re\tr\left(C_h\db B_{u_h}\cdot  v\right) dV&=-\int_U C_h^{i\bar{j}}\re\tr\left(\dd_{\bar j} (B_{u_h})_i\cdot  v\right) dV\\
    &=\int_U C_h^{i\bar{j}}\re\tr\left( (B_{u_h})_i\cdot \dd_{\bar j}v\right) dV+\int_U \dd_{\bar j} C_h^{i\bar{j}}\re\tr\left((B_{u_h})_i\cdot  v\right) dV.
\end{align*}
Similarly, we can divide $R$ into four parts, where
\begin{align*}
    R_1&=\int_U \chi^2\eta^2 C_h^{i\bar{j}}\re\tr\left( (B_{u_h})_i\cdot \dd_{\bar j}s_h\right) dV,\qquad R_2=2\int_U \chi^2\eta\dd_{\bar j}\eta \;C_h^{i\bar{j}}\re\tr\left( (B_{u_h})_i\cdot s_h\right) dV,\\
    R_3&=2\int_U \eta^2\chi\dd_{\bar j}\chi\;C_h^{i\bar{j}}\re\tr\left( (B_{u_h})_i\cdot s_h\right) dV,\quad\; R_4=\int_U \chi^2\eta^2 \dd_{\bar j} C_h^{i\bar{j}}\re\tr\left((B_{u_h})_i\cdot  s_h\right) dV.
\end{align*}
By similar arguments, we have
\begin{align*}
    |R_1|&\le \delta_1\|\chi\eta\db s_h\|_{L^2}^2+C_\delta\|B_{u_h}\|_{L^2}^2,\qquad |R_2|\le C\|s_h\|_{L^\infty}\|B_{u_h}\|_{L^2}\|d\eta_\varepsilon\|_{L^2}=\Psi(\varepsilon|h),\\
    |R_3|&\le C(\|B_{u_h}\|_{L^2}^2+\|s_h\|_{L^2}^2),\qquad\qquad |R_4|\le C(\|B_{u_h}\|_{L^2}^2+\|s_h\|_{L^2}^2).
\end{align*}

Combining the estimation for $L$ and $R$ together and take $0<\delta_1<10^{-2}C^{-1}$ gives
\begin{align*}
    \|\chi\eta_\varepsilon\db s_h\|_{L^2}^2\le C\left(\|B_{u_h}\|^2_{L^2}+\|s_h\|_{L^2}^2+\Psi(\varepsilon|h)\right).
\end{align*}
Finally taking $\varepsilon\to 0$ gives the result.
\end{proof}
\begin{corollary}\label{cor:tangent-curvature}
    In the setting of \cref{lem:logdifference}, take $s_Y:=u^{-1}Yu\in L^2(U)$. Then for any $U_0\subset\subset U$, there exists a constant $C=C(C_0,U_0,U,g)$, 
    \begin{align*}
        \|\bar\partial s_Y\|_{L^2(U_0)}\le C\,\|du\|_{L^2(U)}.
    \end{align*}
\end{corollary}
\begin{proof}
    \Cref{lem:logdifference} gives $\|s_h\|_{L^2(U_1)}\le C\|Yu\|_{L^2(U_1)}\le C\|du\|_{L^2(U)}$. Since $B_{u_h}=u_{h}^{-1}\dd u_{h}$ and $C_0^{-1}I\le u\le C_0I$, one has $\|B_{u_h}\|_{L^2(U_1)}\le C\|du_h\|_{L^2(U_1)}\le C\|du\|_{L^2(U)}$. Thus \Cref{prop:weakvariation} gives
    \begin{align*}
        \|\db s_h\|_{L^2(U_0)}\le \|\chi\db s_h\|_{L^2(U_1)}\le C\left(\|B_{u_h}\|_{L^2(U_1)}+\|s_h\|_{L^2(U_1)}\right)\le C\|du\|_{L^2(U)}.
    \end{align*}
    So $\db s_h$ are uniformly bounded in $L^2$. Then there exists weak convergent subsequence $\db s_{h_i}\rightharpoonup q\in L^2(U_0)$. Since $s_h\to s_Y$ strongly in $L^2$, one has $\db s_Y=q\in L^2(U_0)$. By the weak lower semi-continuity in the Hilbert space $L^2(U_0)$, one has
    \begin{align*}
        \|\db s_Y\|_{L^2(U_0)}\le \liminf_{i\to\infty} \|\db s_{h_i}\|_{L^2(U_0)}\le C\|du\|_{L^2(U)}.
    \end{align*}
\end{proof}

\begin{theorem}[Finite energy]\label{thm:energy}
For any smooth Hermitian--Einstein metric $H$ on $E|_X$ compatible with $\mathcal{E}_*$, one has $F_H\in L^2(X,H,\omega)$.
\end{theorem}
\begin{proof}
It suffices to show \eqref{eq:u-du-bound} and Hermitian--Einstein equation implies $F_u\in L^2$. \Cref{cor:tangent-curvature} gives $\db(u^{-1}\dd_xu),\db(u^{-1}\dd_yu)\in L^2$. Since $\dd_w=\frac12(\dd_x-\ii\dd_y)$, we also have $\db(u^{-1}\dd_wu)\in L^2$, i.e. $F_{w\bar{z}}, F_{w\bar{w}}\in L^2$. Moreover, the Hermitian symmetry of the curvature gives $F_{z\bar{w}}=F_{w\bar{z}}^{*u}\in L^2$ as well, because $u,u^{-1}$ are bounded.

The Hermitian--Einstein equation says
\begin{align*}
    g^{z\bar z}F_{z\bar z}+g^{z\bar w}F_{z\bar w}
+g^{w\bar z}F_{w\bar z}+g^{w\bar w}F_{w\bar w}=\lambda\Id.
\end{align*}
Since $g$ is smooth and all terms except $F_{z\bar z}$ lie in $L^2$; hence so does this last component. Hence we have $F_u\in L^2$.
\end{proof}

\subsection{Normal form of singular connection with finite energy}
We introduce a normal form of singular connection with finite energy due to Sibner--Sibner \cite{Sibner92}, see also \cite{Rade94I,Rade94II,Rade95}, and then apply to our setting.
\begin{lemma}[Finite-action normal form]
\label{lem:critical-gauge}
Let \(A\) be a smooth unitary connection on \(E|_{U\setminus D}\), where \(U=\Delta_z\times\Delta_w\) is a divisor chart. In a unitary trivialization of $E$ on \(U\), write $A=\frac12 \tr A\cdot\Id+A^\perp$.
Suppose $F_A^\perp\in L^2(U)$.

Then there exists a holonomy parameter $b_0\in [0,\frac12]$ such that for almost every
\(x\in D\cap U\), the limit of holonomy conjugacy class exists and
\begin{align*}
\lim_{r\to0}\left[\operatorname{Hol}_{\gamma_{r,x}}(A^\perp)\right]
=\left[
\begin{pmatrix}
e^{-2\pi\sqrt{-1}b_0}&0\\
0&e^{2\pi\sqrt{-1}b_0}
\end{pmatrix}
\right],
\end{align*}
where $\gamma_{r,x}(\theta)=(re^{\ii\theta},w_0)$ for $x=(0,w_0)\in U\cap D$.

Moreover, after shrinking \(U\), there exist a gauge transformation $g\in \Gamma(U\setminus D,SU(E))$ such that
\begin{align*}
g\cdot A^\perp=\widetilde B_{b_0}+a,\qquad a\in W^{1,2}_{\widetilde B_{b_0}},
\end{align*}
where $\widetilde B_{b_0}:=d+\sqrt{-1}\diag(b_0,-b_0)d\theta$. Moreover,
\begin{align*}
    \left\|a\right\|_{W^{1,2}_{\widetilde B_{b_0}}}\le C\|F_A^\perp\|_{L^2(U)}.
\end{align*}
\end{lemma}

\begin{proof}
Locally $A^\perp$ is an \(SU(2)\)-connection and $F_{A^\perp}=F_A^\perp$. Hence \(F_{A^\perp}\in L^2\). The connection \(A^\perp\) is smooth on \(U\setminus D\), hence belongs to \(W^{1,2}_{\mathrm{loc}}\) there.  After shrinking \(U\), absolute continuity of the \(L^2\)-integral makes its curvature norm smaller than the constant in Sibner--Sibner~\cite[Theorem~5.1]{Sibner92}.
Then the results follow from the limiting-holonomy theorem \cite[Theorem~4.1]{Sibner92} and classification theorem ~\cite[Theorem~5.1]{Sibner92}. Choosing the standard representative of an \(SU(2)\) conjugacy class gives \(b_0\in[0,\frac12]\).
\end{proof}

Since the compatiblity condition gives the leading term of the singular connections, we can get explicit value of the holonomy parameter after controlling the effect of error terms.

\begin{lemma}\label{lem:limiting-residue}
Let \(H\) be a Hermitian--Einstein metric on \(E|_X\) compatible with $\mathcal{E}_*$, and \(A\) be any \(h\)-unitary representative of the Chern connection of \(H\). On a divisor chart, let $b_0\in[0,\tfrac12]$ be the holonomy parameter in \cref{lem:critical-gauge}. Then $b_0=\beta$.
\end{lemma}
\begin{proof}
We take a explicit model connection to compute the holonomy and estimate the error term. Choose a holomorphic frame \((e_1,e_2)\) on \(U\) such that \(e_1|_D\) spans \(L\). Take the model metric $K=\|\sigma\|^{2\alpha_2}h_2\oplus \|\sigma\|^{2\alpha_1}h_1$ with respect to the holomorphic splitting $E|_U=\mathbb C e_1\oplus\mathbb C e_2$, where \(h_2\) and \(h_1\) are smooth Hermitian metrics on the two line summands, chosen such that $\det K=h_{P,\mathrm{par}}$. Put $h_0:=h_2\oplus h_1$, then $\det h_0=h_P$.

Let $b:=K^{-1}H$, then the compatibility and self-adjointness give $b,b^{-1}\in L^\infty$ and $\nabla_Kb\in L^2$. By \Cref{prop:metricdet}, one has $\det H=h_{P,\mathrm{par}}=\det K$, thus one also has $\det b=1$. Moreover, one has $\nabla_H=\nabla_K+b^{-1}\nabla^{1,0}_Kb$.

Set $q:=b^{1/2}:(E,H)\to (E,K)$ to be the canonical isometry. The properties of $b$ gives
\begin{align*}
q,q^{-1}\in L^\infty,\qquad \nabla_Kq\in L^2,\qquad \det q=1.
\end{align*}

Put $\widetilde \nabla:=q \nabla_H q^{-1}$ to be pullback of $\nabla_H$ by $q^{-1}$, then
\begin{align*}
    \widetilde \nabla-\nabla_K=q(\nabla_H-\nabla_K) q^{-1}-(\nabla_Kq) q^{-1}
 =
q b^{-1}(\nabla^{1,0}_Kb) q^{-1}-(\nabla_Kq) q^{-1}\in L^2.
\end{align*}

Next we define the isometry
\begin{align*}
    G:=\|\sigma\|^{\alpha_2}\Id_{\mathbb C e_1}
   \oplus \|\sigma\|^{\alpha_1}\Id_{\mathbb C e_2}:(E,K)\longrightarrow(E,h_0).
\end{align*}
 Put $A_*:=G\widetilde \nabla G^{-1}$ and $A_K:=G\nabla_KG^{-1}$. Since \(G\) is an isometry, one has $A_*-A_K\in L^2$. Direct computation gives
\begin{align*}
A_K=A_{h_0}+
\begin{pmatrix}
\alpha_2&0\\
0&\alpha_1
\end{pmatrix}
(\partial-\bar\partial)\log \|\sigma\|,
\end{align*}
where \(A_{h_0}\) is a smooth diagonal unitary connection. Since $\|\sigma\|=r\rho$ for a smooth function $\rho>0$ and $(\partial-\bar\partial)\log r=\sqrt{-1}\,d\theta$, one has
\begin{align*}
A_K^\perp=\sqrt{-1}\,\beta
\begin{pmatrix}
1&0\\
0&-1
\end{pmatrix}d\theta+\text{smooth term}=\widetilde B_{\beta}+\text{smooth term}.
\end{align*}
$A_*-A_K\in L^2$ gives $A_*^\perp-\widetilde B_\beta\in L^2$.

Since $\det h_0=\det h=h_P$, the canonical positive isometry $T:(E,h_0)\longrightarrow(E,h)$ has determinant one. Hence
\begin{align*}
T\circ G\circ q:(E,H)\longrightarrow(E,h)
\end{align*}
is an isometry whose determinant is the canonical
isometry from \(h_{P,\mathrm{par}}\) to \(h_P\). It therefore differs from the canonical \(h\)-unitary representative by an element of
\(SU(E,h)\), and hence is an \(h\)-unitary representative. Any two \(h\)-unitary representatives differ by an \(SU(E,h)\)-gauge transformation, which have the same conjugacy classes. Thus the parameter \(b_0\) can be computed using \(A_*^\perp\).

Actually it suffices to show an $L^2$ term doesn't change the holonomy class. Assume $A_*^\perp-\widetilde B_\beta=c\in L^2$. Since \(c\in L^2(U)\), Fubini's theorem shows that for almost every \(w\),
\begin{align*}
\int_0^\epsilon E(r)\,dr<\infty,
\qquad\text{ where }
E(r):=\int_0^{2\pi}|c(r,\theta,w)|^2\,r\,d\theta.
\end{align*}
Choose \(w\) for which the above inequality holds and for which the holonomy limit exists. Since \(E\in L^1(0,\epsilon)\),
there is a sequence \(r_j\to0\) such that $r_jE(r_j)\to 0$.

Since $|\partial_\theta|\leq Cr$,
\begin{align}\label{eq:gronwall-estimate}
    \int_0^{2\pi}
\bigl|c(\partial_\theta)(r_j,\theta,w)\bigr|\,d\theta\leq
Cr_j\int_0^{2\pi}|c(r_j,\theta,w)|\,d\theta\leq C\bigl(r_jE(r_j)\bigr)^{1/2}
\to0.
\end{align}

Let \(U_j(\theta)\) denote parallel transport for \(A_*^\perp\) along \(\gamma_{r_j,x}\), normalized by \(U_j(0)=\Id\). Let $X_\beta:=\ii \diag(\beta,-\beta)$,
\begin{align*}
U_j'=-\bigl(X_\beta+c(\partial_\theta)\bigr)U_j.
\end{align*}
The parallel transport for the model connection \(X_\beta d\theta\) is $U_0(\theta)=e^{-\theta X_\beta}$. Set $V_j:=U_0^{-1}U_j$, then
\begin{align*}
V_j'=-U_0^{-1}c(\partial_\theta)U_0\,V_j,
\qquad
V_j(0)=\Id.
\end{align*}
Since \(U_0\) is unitary, Gronwall's inequality and \eqref{eq:gronwall-estimate} give
\begin{align*}
V_j(2\pi)\to\Id.
\end{align*}
Therefore
\begin{align*}
\operatorname{Hol}_{\gamma_{r_j,x}}(A_*^\perp)
=
U_j(2\pi)
\to
e^{-2\pi X_\beta}
=
\begin{pmatrix}
e^{-2\pi\sqrt{-1}\beta}&0\\
0&e^{2\pi\sqrt{-1}\beta}
\end{pmatrix}.
\end{align*}

On the other hand, by the defining property of \(b_0\),
\begin{align*}
\lim_{r\to0}
\left[
\operatorname{Hol}_{\gamma_{r,x}}(A^\perp)
\right]
=
\left[
\begin{pmatrix}
e^{-2\pi\sqrt{-1}b_0}&0\\
0&e^{2\pi\sqrt{-1}b_0}
\end{pmatrix}
\right].
\end{align*}
Since \(A_*^\perp\) and \(A^\perp\) differ by an \(SU(2)\)-gauge
transformation, they have the same holonomy conjugacy classes. Hence $b_0=\beta$.
\end{proof}

\begin{corollary}\label{cor:critical-weighted-realization}
Let \(H\) be a Hermitian--Einstein metric on \(E|_X\) compatible
with \(\mathcal E_*\), and let \(A\) be any \(h\)-unitary
representative of its Chern connection. Then, on every sufficiently
small divisor chart \(U\), there exists an \(SU(E)\)-gauge
transformation \(u\) on \(U\setminus D\) such that \(u\cdot A\)
satisfies the defining condition of \cref{def:gaugemoduli} with
\(p=2\); namely,
\begin{align*}
\det(u\cdot A)=a_{P,\mathrm{par}},
\qquad
u\cdot A-A_0
\in
W^{1,2}_{B_\beta}
\bigl(U;T^*U\otimes\mathfrak{su}(E)\bigr).
\end{align*}
\end{corollary}

\begin{proof}
By \cref{lem:critical-gauge}, after an \(SU(E)\)-gauge transformation, the trace-free part of \(A\) has the form $u\cdot A^\perp=\widetilde B_{b_0}+a$ with $a\in W^{1,2}_{\widetilde B_{b_0}}$, where $\widetilde B_{b_0}:=d+\sqrt{-1}\diag(b_0,-b_0)d\theta$.

After possibly composing the gauge $u$ with a
constant \(SU(2)\) transformation, we may assume the unitary trivialization satisfying \(L|_{D\cap U}=\mathbb C e_1\). By \Cref{prop:banach}, the trace-free part of the reference connection in \cref{def:gaugemoduli} is locally
\begin{align*}
A_0^\perp
=
\sqrt{-1}\,\beta
\begin{pmatrix}
1&0\\
0&-1
\end{pmatrix}\eta+\Gamma_0,
\end{align*}
where \(\Gamma_0\) is smooth and preserves the splitting
\(L\oplus Q\). Since \(\eta-d\theta\) is smooth, the two model
connections differ by a smooth splitting-preserving one-form.
Hence the critical covariant \(W^{1,2}\) condition above is precisely
\begin{align*}
u\cdot A^\perp-A_0^\perp
\in W^{1,2}_{B_\beta}
\bigl(U;T^*U\otimes\mathfrak{su}(E)\bigr).
\end{align*}
Finally, the gauge transformation is \(SU(E)\)-valued, so it preserves the determinant connection. Thus $\det A=a_{P,\mathrm{par}}$, and consequently $A-A_0\in W^{1,2}_{B_\beta}\bigl(U;T^*U\otimes\mathfrak{su}(E)\bigr)$. This is exactly the local defining condition of \cref{def:gaugemoduli} at \(p=2\).
\end{proof}

\subsection{Improved regularity}
The purpose of this subsection is to pass from the weighted regularity $W^{1,2}_{B_\beta}$ obtained in
\cref{cor:critical-weighted-realization} to the exponent $p>2$ used in \cref{def:gaugemoduli}. The techniques we use here are due to Kronheimer--Mrowka \cite[\S4(iii), \S7(i)]{KM93}. 

We first recall some notations before and in \cite{KM93}. Fix the exponent $p>2$ chosen in \eqref{eq:p-assumption}, then the local elliptic theory of \cite[Proposition~4.11]{KM93} is available. Since the regularity problem is local, so we consider a divisor chart $U=\Delta_z\times\Delta_w$, $D\cap U=\{z=0\}$. We identify it with the standard pair \((B^4,B^2)\), and use the local mixed gauge on $S^4$ introduced by Kronheimer--Mrowka in \cite[\S4(iii)]{KM93}. Their construction is carried out after conformally identifying $(B^4,B^2)$ with a hemisphere of \((S^4,S^2)\). In the adapted unitary trivialization on the divisor chart, $A_0^\perp=\widetilde B_\beta+\Gamma$, where $\widetilde B_\beta=d+\ii\diag(\beta,-\beta)d\theta$ and $\Gamma$ is smooth and preserves the splitting.

Choose a rational approximation $2\beta_0=m/\nu$ of $2\beta$ satisfying the conditions of
\cite[Proposition~4.8 and Lemma~4.9]{KM93}, and let $\check g$ be the auxiliary
orbifold metric of cone angle $2\pi/\nu$ occurring there.  Define the mixed operator \cite[(4.13)]{KM93} to $\mathfrak{su}(E)$-valued tensor field
\begin{align}
    \label{eq:mixedgauge}
 d_0a:=d^*_{g_{\mathrm{rd}}}a_{\mathrm{diag}}+(d_{2\beta_0})^{*}_{\check g}\;a_{\mathrm{off}},
 \qquad
 a=a_{\mathrm{diag}}+a_{\mathrm{off}},
\end{align}
where $g_{\mathrm{rd}}$ is the round metric on $S^4$, and $d_{2\beta_0}\xi:=(d+2\sqrt{-1}\beta_0\,d\theta)\xi$ on the off-diagonal summand. 

For a Riemannian metric \(g\), denote by $\pi_g^+:\Lambda^2T^*U\longrightarrow\Lambda_g^+$
the orthogonal projection to $g$-self-dual two-forms.  If
\(B\) is a unitary connection, set
\begin{align*}
 d^+_{B,g}
 :=
 \pi_g^+\circ d_B:
 \Omega^1(\mathfrak{su}(E))
 \longrightarrow
 \Omega_g^+(\mathfrak{su}(E)).
\end{align*}
In fact \(d^+_{B,g}\) is the linearization of the self-dual curvature map at \(B\), since we have
\begin{align*}
 (F_{B+a})^{+}_g=(F_B)^{+}_{g}+d^+_{B,g}a+(a\wedge a)^{+}_{g}.
\end{align*}

\begin{proposition}[Local mixed left inverse]
\label{prop:mixedlinear}
Let
\begin{equation}
\label{eq:mixed-first-order-operator}
 \mathcal D_\omega:=
 \bigl(d_0,d^+_{A_0^\perp,\omega}\bigr):
 W^{1,q}_{B_\beta}
 \bigl(U;T^*U\otimes\mathfrak{su}(E)\bigr)
 \to
 L^q\bigl(U;\mathfrak{su}(E)\bigr)
 \oplus
 L^q\bigl(U;\Lambda_\omega^+\otimes\mathfrak{su}(E)\bigr).
\end{equation}

Then after possibly shrinking the divisor chart, there is a bounded linear operator \(P_\omega\) for every \(q\in[2,p]\)
\begin{align*}
 P_\omega:
 L^q\bigl(
 U;\mathfrak{su}(E)\oplus
 \Lambda_\omega^+\otimes\mathfrak{su}(E)
 \bigr)
 \longrightarrow
 W^{1,q}_{B_\beta}
 \bigl(
 U;T^*U\otimes\mathfrak{su}(E)
 \bigr),
\end{align*}
such that for $a\in W^{1,q}_{B_\beta}$ with $\supp a\subset U$, we have
\begin{equation}
\label{eq:mixed-left-inverse}
 P_\omega\mathcal D_\omega a=a.
\end{equation}
Moreover, $P_\omega$ is compatible with $q$: if
\(q_1,q_2\in[2,p]\) and \(f\in L^{q_1}\cap L^{q_2}\), then $P_\omega^{(q_1)}f=P_\omega^{(q_2)}f$.
\end{proposition}

\begin{proof}
We first consider the local model $\mathcal D_0:=\left(d_0,d^+_{\widetilde B_\beta,g_{\rm rd}}\right)$ on $(S^4,S^2)$.
Set
\begin{align*}
 T_q:=\ker d_0\cap
 W^{1,q}_{\widetilde B_\beta}(S^4;T^*S^4\otimes\mathfrak{su}(2)).
\end{align*}
\cite[\S4(i)--(iii), Proposition~4.11]{KM93}
shows that
\begin{align*}
 d^+_{\widetilde B_\beta,g_{\rm rd}}
 :
 T_q
 \longrightarrow
 L^q(S^4;\Lambda^+_{g_{\rm rd}}\otimes\mathfrak{su}(2))
\end{align*}
is an isomorphism. Let $Q$ denote its inverse. By \cite[Proposition~4.11]{KM93}, \(Q\) and \(T_q\) are independent of \(q\).

By standard elliptic theory applied to the two summands of \(d_0\), together with the norm comparison \cite[Lemma~3.9]{KM93}, there is a bounded operator
\begin{align*}
 R:L^q(S^4;\mathfrak{su}(2))
 \to
 W^{1,q}_{\widetilde B_\beta}
 (S^4;T^*S^4\otimes\mathfrak{su}(2)),
\end{align*}
such that
\begin{align*}
 d_0Rf=f-\bar f_{\diag},
\end{align*}
where $\bar f_{\diag}$ is the integral average of $f_{\diag}$. Here the Laplacian on diagonal part is invertible modulo constants, whereas the Laplacian on the off-diagonal part is invertible. The construction is independent of \(q\). 

 Define
\begin{align*}
 \mathcal P_{0}(f,\eta):=Rf+Q\bigl(\eta-d^+_{\widetilde B_\beta,g_{\rm rd}}Rf\bigr).
\end{align*}

For $a\in W^{1,q}_{\widetilde B_\beta}(S^4;T^*S^4\otimes\mathfrak{su}(2))$, we have $\overline{(d_0a)}_{\diag}=0$, since $\int (d_0a)_{\diag}=\langle d^*a_{\rm diag},1\rangle=\langle a_{\rm diag},d1\rangle=0$. Hence
\begin{align*}
 d_0\bigl(a-Rd_0a\bigr)=d_0a-d_0R(d_0a)=\overline{(d_0a)}_{\diag}=0.
\end{align*}
Therefore \(a-Rd_0a\in T_q\). Thus
\begin{align*}
    \mathcal{P}_0\mathcal{D}_0a=Rd_0a+Qd^+_{\widetilde B_\beta,g_{\rm rd}}\bigl(a-Rd_0a\bigr)=Rd_0a+a-Rd_0a=a. 
\end{align*}
Thus $\mathcal{P}_0$ is a bounded left inverse of
\(\mathcal D_{0}\) on \(S^4\). Since both \(R\) and \(Q\) are independent of \(q\), $\mathcal{P}_0$ is also compatible with $q$.

We now pass from the model $(S^4,\mathcal{D}_0)$ to the divisor chart $(U,\mathcal{D}_\omega)$.
After shrinking \(U\), identify \(\Lambda_\omega^+\) with \(\Lambda_{g_{\rm rd}}^+\) by $\iota_\omega:=\pi^+_{g_{\rm rd}}\big|_{\Lambda_\omega^+}$. Write
\begin{align*}
 \mathcal D'_\omega:=\bigl(d_0,\iota_\omega d^+_{A_0^\perp,\omega}\bigr)=:\mathcal{D}_0+K_\omega
\end{align*}
on \(U\), compared with $\mathcal{D}_\omega=\bigl(d_0,d^+_{A_0^\perp,\omega}\bigr)$. Extend \(K_\omega\) to \(S^4\) using a cutoff function, and denote the resulting operator by \(\widehat K_\omega\).
The local comparison of the two metrics, together with the smoothness of \(A_0^\perp-\widetilde B_\beta\) and the embedding $W^{1,q}_{\widetilde B_\beta}\hookrightarrow L^{4q/(4-q)}$ gives
\begin{align*}
 \|\widehat K_\omega a\|_{L^q}
 \le
 \varepsilon(U)
 \|a\|_{W^{1,q}_{\widetilde B_\beta}},
\end{align*}
where $\varepsilon(U)\to0$ as \(U\) shrinks, uniformly for \(q\in[2,p]\).
Hence, after shrinking \(U\) once more, we can assume $\|\mathcal P_{0}^{(q)}\widehat K_\omega\|<1$.
Then $\widehat P_\omega^{(q)}:=\bigl(\Id+\mathcal P_0^{(q)}\widehat K_\omega\bigr)^{-1} \mathcal P_{0}^{(q)}$ satisfies
\begin{align*}
 \widehat P_\omega^{(q)}
 \bigl(\mathcal{D}_0+\widehat K_\omega\bigr)
 =\Id .
\end{align*}

Let \(E_U\) denote extension by zero from \(U\) to \(S^4\), and define
\begin{align*}
 P_\omega^{(q)}(f,\eta)
 :=
 \left.
 \widehat P_\omega^{(q)}
 E_U\bigl(f,\iota_\omega\eta\bigr)
 \right|_U .
\end{align*}
If \(a\in W^{1,q}_{B_\beta}(U)\) and
\(\supp a\Subset U\), then \(E_Ua\in
W^{1,q}_{\widetilde B_\beta}(S^4)\). Therefore we have $P_\omega^{(q)}\mathcal D_\omega a=a$. 
Finally, \(\mathcal P_{0}\) and \(\widehat K_\omega\) are independent of \(q\), hence \(P_\omega^{(q)}\) is independent of $q$.
\end{proof}

The following proposition is a slight modification of \cite[Proposition~7.2(i)--(ii)]{KM93}.
\begin{proposition}\label{prop:criticalmixed}
Let $A=A_0+a$ with $a\in W^{1,2}_{B_\beta}\bigl(U;T^*U\otimes\mathfrak{su}(E)\bigr)$. If $\|F_A^\perp\|_{L^2(U)}$ is sufficiently small, then after shrinking $U'\subset\subset U$, there is $u\in W^{2,2}_{B_\beta,\mathrm{loc}}(U';SU(E))$ such that
\begin{align*}
    u\cdot A=A_0+\widetilde a,\qquad d_0\widetilde a=0,\qquad \|\widetilde a\|_{W^{1,2}_{B_\beta}(U')}
 \le C\|F_A^\perp\|_{L^2(U)}.
\end{align*}
\end{proposition}

\begin{proof}
Choose smooth $a_j$ converging to $a$ in $W^{1,2}_{B_\beta}$, and let $A_j:=A_0+a_j$.  Since
$W^{1,2}_{B_\beta}\hookrightarrow L^4$,
\begin{align*}
 F_{A_j}^\perp-F_A^\perp=d_{A_0^\perp}(a_j-a)
 +(a_j-a)\wedge a_j+a\wedge(a_j-a)
 \to 0 \quad\text{ in }L^2.
\end{align*}
For each $j$, \cite[Proposition~7.2]{KM93} gives gauges $u_j\in W^{2,p}_{B_\beta}(U';SU(E))$ such that
\begin{align*}
    u_j\cdot A_j=A_0+\widetilde a_j,\quad d_0\widetilde a_j=0,\quad \|\widetilde a_j\|_{W^{1,2}_{B_\beta}(U')}\le C\|F_{A_j}^\perp\|_{L^2(U)}.
\end{align*}
Let \(\nabla_0\) denote the covariant derivative induced by \(A_0\)
on \(\End(E)\) and the associated tensor bundles.  The gauge identity $u_j\cdot(A_0+a_j)=A_0+\widetilde a_j$ is equivalent to
\begin{align*}
    \nabla_0u_j=u_ja_j-\widetilde a_j u_j .
\end{align*}
Since \(u_j\) is unitary and $a_j$, $\widetilde a_j$ are uniformly bounded in $W^{1,2}_{B_\beta}(U')$, the embedding \(W^{1,2}_{B_\beta}\hookrightarrow L^4\) gives
\begin{align*}
 \|\nabla_0u_j\|_{L^4(U')}
 \leq
 C\bigl(
 \|a_j\|_{L^4(U')}
 +
 \|\widetilde a_j\|_{L^4(U')}
 \bigr)
 \leq C .
\end{align*}
Taking the covariant derivative yields,
\begin{align*}
    \nabla_0^2u_j=(\nabla_0u_j)\otimes a_j+u_j\nabla_{A_0}a_j-(\nabla_{A_0}\widetilde a_j)u_j-\widetilde a_j\otimes\nabla_0u_j.
\end{align*}
Hence, by Hölder's inequality,
\begin{align*}
 \|\nabla_0^2u_j\|_{L^2(U')}
 \leq
 C\left(
 \|\nabla_0u_j\|_{L^4}
 \bigl(\|a_j\|_{L^4}
      +\|\widetilde a_j\|_{L^4}\bigr)
 +
 \|\nabla_{A_0}a_j\|_{L^2}
 +
 \|\nabla_{A_0}\widetilde a_j\|_{L^2}
 \right)
 \leq C .
\end{align*}
Since \(|u_j|\) is uniformly bounded and the \(A_0\)- and \(B_\beta\)-Sobolev norms on \(\End(E)\) are equivalent on \(U'\), we have $\|u_j\|_{W^{2,2}_{B_\beta}(U')}\leq C$.
After passing to a subsequence, we have the weak convergence
\begin{align*}
    u_j\rightharpoonup u
 \quad\text{in }W^{2,2}_{B_\beta}(U'),\qquad \widetilde a_j\rightharpoonup\widetilde a
 \quad\text{in }W^{1,2}_{B_\beta}(U').
\end{align*}

Since \(a_j\to a\) strongly in \(W^{1,2}_{B_\beta}(U')\), in particular $a_j\to a$ in $L^4$. Moreover, $\widetilde a_j\rightharpoonup\widetilde a$ in $L^4$.
Passing to the limit in $\nabla_0u_j=u_ja_j-\widetilde a_ju_j$ gives $\nabla_0u=ua-\widetilde au$ in $\mathcal{D}'(U')$. Hence
\begin{align*}
 u\cdot A=A_0+\widetilde a .
\end{align*}

Since $d_0:W^{1,2}_{B_\beta}(U')\to L^2(U')$ is bounded and \(d_0\widetilde a_j=0\), weak convergence gives $d_0\widetilde a=0$.
Finally, weak lower semicontinuity and the strong \(L^2\)-convergence of the curvatures give
\begin{align*}
    \|\widetilde a\|_{W^{1,2}_{B_\beta}(U')}\le \liminf_{j\to\infty}\|\widetilde a_j\|_{W^{1,2}_{B_\beta}(U')}\le C\lim_{j\to\infty}\|F_{A_0+a_j}^{\perp}\|_{L^2(U)}=C\|F_A^\perp\|_{L^2(U)}.
\end{align*}
This proves the proposition.
\end{proof}

\begin{proposition}[Improved regularity]\label{prop:upgrade}
In the setting of \cref{prop:criticalmixed}. If
$(F_A^\perp)^+=0$, then after shrinking the divisor chart, $\widetilde a\in W^{1,p}_{B_\beta}\bigl(T^*U\otimes\mathfrak{su}(E)\bigr)$.
\end{proposition}
\begin{proof}
By \cref{prop:criticalmixed}, we have $u\cdot A=A_0+\widetilde a$, $d_0\widetilde a=0$ and $\widetilde a\in W^{1,2}_{B_\beta}(U)$.
Since the ASD equation is gauge invariant,
\begin{align*}
 0=(F_{A_0+\widetilde a}^{\perp})^+
   =(F_{A_0^\perp})^+
     +d^+_{A_0^\perp,\omega}\widetilde a
     +(\widetilde a\wedge\widetilde a)^+ .
\end{align*}
Hence we have
\begin{align}
    \label{eq:critical-mixed-equation}
 \mathcal D_\omega\widetilde a=\bigl(0,-(F_{A_0^\perp})^+
          -(\widetilde a\wedge\widetilde a)^+\bigr),
\end{align}
where $\mathcal D_\omega=(d_0,d^+_{A_0^\perp,\omega})$ defined in \Cref{prop:mixedlinear}.

Choose \(V\subset\subset U\) and a cutoff function
\(\chi\in C_c^\infty(U)\) such that \(\chi=1\) on \(V\), and set $b:=\chi\widetilde a$.
Then \(\supp b\subset\subset U\), and \eqref{eq:critical-mixed-equation}
gives
\begin{equation}
\label{eq:localized-mixed-equation}
 \mathcal D_\omega b=-N_{\widetilde a}b+f,
\end{equation}
where $N_{\widetilde a}v:=\bigl(0,(\widetilde a\wedge v)^+\bigr)$ and $f:=[\mathcal D_\omega,\chi]\widetilde a-\bigl(0,\chi(F_{A_0^\perp})^+\bigr)$.
Here $[\mathcal D_\omega,\chi]\widetilde a:=\mathcal D_\omega(\chi\widetilde a)-\chi\mathcal D_\omega\widetilde a$ is a zeroth-order operator with $\left|[\mathcal D_\omega,\chi]\widetilde a\right|\le C|d\chi|\cdot |\widetilde a|$.
Since \(\widetilde a\in W^{1,2}_{B_\beta}\subset L^4\) and \(F_{A_0^\perp}\) is smooth, we have $f\in L^4(U)$.

Put $q^*=\frac{4q}{4-q}$ for $q=2,p$.
The weighted Sobolev embedding and Hölder's inequality give
\begin{align*}
 \|N_{\widetilde a}v\|_{L^q}
 \le
 C\|\widetilde a\|_{L^4}
   \|v\|_{L^{q^*}}
 \le
 C\|\widetilde a\|_{L^4}
   \|v\|_{W^{1,q}_{B_\beta}}.
\end{align*}
After shrinking the chart, we may assume $\|P_\omega N_{\widetilde a}\|_{W^{1,q}_{B_\beta}\to W^{1,q}_{B_\beta}}<\frac12$. Thus $\Id+P_\omega N_{\widetilde a}$ is invertible.

Since $\supp b\subset U$, the left-inverse identity
for \(P_\omega\) and \eqref{eq:localized-mixed-equation} imply
\begin{align}
    \label{eq:critical-fixed-point}
 (\Id+P_\omega N_{\widetilde a})b=P_\omega f\in W^{1,q}_{B_\beta}(U).
\end{align}
Since \(\Id+P_\omega N_{\widetilde a}\) is invertible on
\(W^{1,p}_{B_\beta}(U)\) and \(f\in L^p\),
there exists
\begin{align*}
 c:=(\Id+P_\omega N_{\widetilde a})^{-1}P_\omega f
 \in W^{1,p}_{B_\beta}(U).
\end{align*}
The compatibility of \(P_\omega\) for the exponents \(2\) and \(p\) shows that \(c\) satisfies
\eqref{eq:critical-fixed-point} also in
\(W^{1,2}_{B_\beta}(U)\).  Since
\(\Id+PN_{\widetilde a}\) is invertible on this space, \(c=b\).  Therefore $b\in W^{1,p}_{B_\beta}(U)$.

Since \(b=\widetilde a\) on \(V\), we have $\widetilde a\in
 W^{1,p}_{B_\beta}
 \bigl(V;T^*V\otimes\mathfrak{su}(E)\bigr)$.
\end{proof}

\begin{remark}\label{rem:example-not-W2p}
The weighted $W^{1,p}$ regularity in the above
proposition cannot, in general, be improved to
$W^{2,p}$. Consider the Euclidean divisor model with
coordinates $(r,\theta,u,v)$, and choose
$T_1,T_2,T_3\in\mathfrak{su}(2)$ satisfying
$[T_i,T_{i+1}]=2T_{i+2}$ for $i=1,2,3$.
Let $A_\beta=-\beta T_1\,d\theta$.
For $c>0$, set
\begin{align*}
\begin{gathered}
A:=fT_1\,d\theta+G(T_2\,du+T_3\,dv),\\
f:=-\frac{\beta+(1+\beta)\tau}{1-\tau},
\qquad
G:=\frac{c(1+2\beta)r^{2\beta}}{1-\tau},
\qquad
\tau:=c^2r^{2+4\beta}.
\end{gathered}
\end{align*}
Then $A$ is ASD, and $a:=A-A_\beta$ satisfies
\begin{align*}
d_{A_\beta}^*a=0,\qquad
a\in W^{1,p}_{A_\beta}
\quad\text{for }2<p<\frac{2}{1-2\beta}.
\end{align*}
However, $G''(r)\sim c(2\beta+1)(2\beta)(2\beta-1)r^{2\beta-2}\notin L^p$ for any $p>2$.
Hence $a\notin W^{2,p}_{A_\beta}$.
\end{remark}

\begin{corollary}\label{cor:curvatureLp}
There exists a constant $p_*(\beta)>2$ such that for any smooth Hermitian--Einstein metric $H$ on $E|_X$ compatible with $\mathcal{E}_*$, one has $F_H\in L^p(X,H,\omega)$ for any $2\le p<p_*(\beta)$.
\end{corollary}

\begin{proof}
Let $A$ be any $h$-unitary representative of the Chern connection of $H$. 
For \(x\in D\), \Cref{prop:upgrade} shows that on a sufficiently small divisor chart \(U\ni x\), there is \(u\in\Gamma(U\setminus D,SU(E))\) such that $u\cdot A=A_0+a$ and $a\in W^{1,p}_{B_\beta}\bigl(U;T^*U\otimes\mathfrak{su}(E)\bigr)$.
Therefore
\begin{align*}
F_{u\cdot A}
=
F_{A_0}+d_{A_0}a+a\wedge a.
\end{align*}
The model curvature \(F_{A_0}\) extends smoothly across \(D\), and \(d_{A_0}a\in L^p(U)\).  By the weighted multiplication estimate of \cref{prop:banach}, $a\wedge a\in L^p(U)$.
Thus \(F_{u\cdot A}\in L^p(U)\).

Since \(u\) is unitary, one has $F_{u\cdot A}=uF_Au^{-1}$ and $|F_{u\cdot A}|=|F_A|$ pointwise on \(U\setminus D\). Hence \(F_A\in L^p(U\setminus D)\).
A finite collection of such charts covers \(D\), while \(A\) is
smooth on compact subsets of \(X\). Therefore \(F_A\in L^p(X)\).
\end{proof}

Now we can patch the local gauges to obtain a global gauge without changing the topology of the bundle $(E,L)$.
\begin{proposition}
\label{prop:globalrealization}
Let \(H\) be a Hermitian--Einstein metric on \(E|_X\) compatible
with \(\mathcal E_*\), and let \(A\) be any \(h\)-unitary representative of its Chern connection.  Then there exists $u\in\Gamma(X,SU(E))$ such that $u\cdot A\in\mathcal A_\beta^{1,p}(E,L;a_{P,\mathrm{par}})$.
\end{proposition}

\begin{proof}
By \cref{prop:metricdet}, we can assume $\det H=h_{P,\mathrm{par}}$, thus $\det A=a_{P,\mathrm{par}}$. By \cref{cor:curvatureLp}, $F_A\in L^p(X)$ for some $p>2$.

Choose a finite collection of divisor charts
\(\{U_i\}_{i=1}^N\) covering \(D\).  By
\cref{cor:critical-weighted-realization,prop:upgrade}, after shrinking the \(U_i\) if necessary, there exist $u_i\in\Gamma(U_i\setminus D,SU(E))$ such that
\begin{align*}
A_i:=u_i\cdot A=A_{0}+a_i,
\qquad
a_i\in W^{1,p}_{B_\beta}\bigl(U_i;T^*U_i\otimes\mathfrak{su}(E)\bigr).
\end{align*}
On \(U_i\cap U_j\setminus D\), define $g_{ij}:=u_i u_j^{-1}$. Then
\begin{align*}
g_{ij}\cdot A_j=A_i,
\qquad
g_{ij}g_{jk}=g_{ik}.
\end{align*}
By \cref{prop:gauge}, each \(g_{ij}\) extends uniquely to a determinant-preserving unitary isomorphism on \(U_i\cap U_j\), satisfies $g_{ij}\in W^{2,p}_{B_\beta}$ and preserves \(L\) over \(D\cap U_i\cap U_j\).

By \cite[Proposition~3.5]{Biquard96}, the above transition functions $\{g_{ij}\}$ determine an extension \(F\to\X\) of \(E|_X\to X\), and the restrictions \(L|_{D\cap U_i}\) determine a line subbundle $M\subset F|_D$.
More precisely, the local identifications
\begin{align*}
\Phi_i:=u_i^{-1}:(E|_{U_i\setminus D},A_i)
\longrightarrow(E|_{U_i\setminus D},A)
\end{align*}
satisfy the compatibility relation $\Phi_i\circ g_{ij}=\Phi_j$, therefore gives an isomorphism $\Phi:F|_X\to E|_X$. Thus the local connections \(A_i\) define a connection \(A_F:=\Phi^*A\) on \(F|_X\).  In the local trivialization over \(U_i\), we also have $A_F=A_{0,i}+a_i$ with $a_i\in W^{1,p}_{B_\beta}$.
Since every \(g_{ij}\) has determinant one, the fixed identifications $\det(E|_{U_i})\simeq P|_{U_i}$ descend to an identification $\phi_F:\det F\xrightarrow{\sim}P$ and $\det A_F=a_{P,\mathrm{par}}$.

We next identify the topological type of \((F,M)\).
Let $\Psi:(\mathcal E|_X,H)\to(E|_X,h)$ be the unitary isomorphism defining the representative \(A\).
Then \(\Psi\) is holomorphic from
\((\mathcal E|_X,\bar\partial_E)\) to
\((E|_X,\bar\partial_A)\).
Moreover, by compatibility, $C^{-1}h_{\mathrm{par}}\le H\le C h_{\mathrm{par}}$.
Hence, for every
local holomorphic section \(s\) and every \(b\in\mathbb R\),
\begin{align*}
|s|_H=O(r^b)
\quad\Longleftrightarrow\quad
|s|_{h_{\mathrm{par}}}=O(r^b).
\end{align*}
Therefore the parabolic filtration associated by
\cite[Theorem~2.1]{Biquard96} to the Hermitian holomorphic bundle \((\mathcal E|_X,H)\) is precisely the given filtration of \(\mathcal E_*\).

The curvature hypothesis required in
\cite[Theorem~2.1]{Biquard96} holds by
\cref{cor:curvatureLp}.  Hence
\cite[Proposition~3.6]{Biquard96} applies and shows that the bundle pair \((F,M)\) is topologically isomorphic to the underlying bundle pair \((E,L)\). Consequently $c_2(F)=c_2(E)$ and $\deg M=\deg L$.

Since \(p>2\), the \(W^{2,p}_{B_\beta}\) transition functions $g_{ij}$ are continuous.  By the standard smoothing argument, there exist a smooth determinant-fixed bundle pair $(F',M',\phi_{F'})$ and a determinant-preserving unitary isomorphism $V:(F,M,\phi_F)\longrightarrow(F',M',\phi_{F'})$ which is \(W^{2,p}_{B_\beta}\) up to \(D\) and smooth on \(X\).
Replacing \(F,M,A_F\) by $(F',M',V_*A_F)$, we may therefore assume that \((F,M,\phi_F)\) is smooth.
The weighted gauge-action estimates of \cref{prop:banach} give $A_F\in\mathcal A_\beta^{1,p}(F,M;a_{P,\mathrm{par}})$.

Both \((F,M,\phi_F)\) and \((E,L,\phi_E)\) are now smooth determinant-fixed bundle pairs with same topological class.  Hence
\cref{lem:pairclassification} gives a smooth
determinant-preserving unitary pair isomorphism
\begin{align*}
C:(F,M,\phi_F)\xrightarrow{\sim}(E,L,\phi_E).
\end{align*}
Therefore $C_*A_F\in\mathcal A_\beta^{1,p}(E,L;a_{P,\mathrm{par}})$.
Let $u:=C\circ \Phi^{-1}\in \Gamma(X,SU(E))$, then
\begin{align*}
    u\cdot A=C_*\Phi^*A=C_*A_F\in \mathcal A_\beta^{1,p}(E,L;a_{P,\mathrm{par}}).
\end{align*}

\end{proof}

\begin{remark}\label{rem:example-not-canonical}
Here we need an extra gauge $u$ to map into the weighted space, since the canonical $h$-unitary representative may not belong to $\A_\beta^{1,p}$, even for a flat compatible Chern connection.
Indeed, on a divisor chart, let
$\gamma=\alpha_2-\alpha_1\in(0,1)$ and set
\begin{align*}
G=\diag(r^{\alpha_2},r^{\alpha_1}),\qquad
K=G^2,\qquad
H=C^*KC,\qquad
C=\begin{pmatrix}1&c\\0&1\end{pmatrix},
\quad c\neq0\in\mathbb R.
\end{align*}
It's easy to check the metric $H$ is compatible with the model metric $K$ and flat.
The isometry $R=GC$ transforms its Chern connection
into the model connection $A_0$.
However, writing $R=US$ with $S=H^{1/2}$,
the unitary factor $U$ is a rotation with angle
\begin{align*}
\vartheta(r)=
\arctan\frac{cr^\gamma}{1+r^\gamma},
\qquad
\vartheta''(r)
=c\gamma(\gamma-1)r^{\gamma-2}
+o(r^{\gamma-2}).
\end{align*}
The connection defined by $S$ is $U^{-1}A_0U$.
Its radial remainder is $\vartheta'(r)J\,dr$,
where $J$ is a fixed nonzero skew-Hermitian matrix.
Since $\vartheta''\notin L^p$ for every $p>2$,
this canonical representative does not belong to $\A_\beta^{1,p}$.
\end{remark}

\begin{lemma}
\label{prop:realizationunique}
Let \(A\) be an \(h\)-unitary representative of the Chern connection of a compatible Hermitian--Einstein metric of $\mathcal{E}_*$. Suppose $u_1,u_2\in\Gamma(X,SU(E))$ satisfy $A_i:=u_i\cdot A\in\mathcal A_\beta^{1,p}(E,L;a_{P,\mathrm{par}})$ for $i=1,2$.

Then
\begin{align*}
    g:=u_1u_2^{-1}\in\mathcal G_\beta^{2,p}(E,L),\qquad g\cdot A_2=A_1. 
\end{align*}
In particular,
\begin{align*}
[A_1]=[A_2]
\quad\text{in}\quad
\mathcal B_\beta^p(E,L;a_{P,\mathrm{par}}).
\end{align*}
\end{lemma}

\begin{proof}
Set $g:=u_1u_2^{-1}$, then \(g\in\Gamma(X,SU(E))\) and $g\cdot A_2=(u_1u_2^{-1})\cdot(u_2\cdot A)=u_1\cdot A=A_1$.
Both \(A_1\) and \(A_2\) belong to \(\mathcal A_\beta^{1,p}(E,L;a_{P,\mathrm{par}})\). Hence \cref{prop:gauge} applied to the determinant-preserving unitary isomorphism \(g\) gives $g\in\mathcal G_\beta^{2,p}(E,L)$.
\end{proof}

\subsection{Parabolic Hermitian--Einstein moduli space}
\begin{definition}
\label{def:parHEmoduli}
Define the space of \textbf{parabolic Hermitian--Einstein connections} $\mathcal A^{\mathrm{par,HE}}_h(E,L;\mathcal{P})$ to be the set of smooth \(h\)-unitary connections \(A\) on
\(E|_X\) for which there exist $(\db_E,H)$ such that \(A\) is an \(h\)-unitary representative (\cref{def:hunitaryrep}) of the Chern connection of
\((E|_X,\bar\partial_E,H)\), where
\begin{enumerate}
\item $\db_E:\Gamma(E)\to \Omega^{0,1}_{\X}(E)$ is an integrable holomorphic structure satisfying $\det\bar\partial_E=\bar\partial_P$ and $L\subset (E|_D,\bar\partial_E|_D)$ is a holomorphic line subbundle;

\item $H$ is a Hermitian--Einstein metric on
\((E|_X,\bar\partial_E)\), compatible with the parabolic bundle $\mathcal E_*=(E,\bar\partial_E,L,\alpha_1,\alpha_2)$.
\end{enumerate}

Take the gauge group $\mathcal G_h(X):=\Gamma(X,SU(E,h))$ acting on \(\mathcal A^{\mathrm{par,HE}}_h(E,L;\mathcal P)\). Define the \textbf{parabolic Hermitian--Einstein connection moduli space}
\begin{align*}
\mathcal M^{\mathrm{par,HE}}_h(E,L;\mathcal{P})
:=
\mathcal A^{\mathrm{par,HE}}_h(E,L;\mathcal{P})
/\mathcal G_h(X).
\end{align*}
\end{definition}

\begin{remark}
\label{rem:parHEobjectmoduli}
The moduli set of \cref{def:parHEmoduli} is naturally identified
with the set of isomorphism classes of compatible parabolic Hermitian--Einstein pairs $(\mathcal E_*,H)$ normalized by $\det H=h_{P,\mathrm{par}}$ (see \Cref{prop:metricdet}) satisfying $\det\db_{\mathcal E}=\db_{\mathcal P}$. More precisely, let $(\mathcal E_{i,*},H_i)=(E,\bar\partial_i,L,\alpha_1,\alpha_2,H_i)$ be two such pairs.  We call them isomorphic if there exists a holomorphic bundle isomorphism $f:(E,\db_1)\to (E,\db_2)$ such that
\begin{align*}
    f(L)=L,\qquad \det f=\Id_P,\qquad f^*H_2=H_1.
\end{align*}

Let \(A_i^h\) denote the canonical \(h\)-unitary representative of the Chern connection of \((E|_X,\bar\partial_i,H_i)\).
Then
\begin{align*}
[(\mathcal E_*,H)]\longmapsto[A^h]
\end{align*}
defines a bijection from these isomorphism classes onto \(\mathcal M^{\mathrm{par,HE}}_h(E,L;P)\).

Indeed, suppose first that
\(f:(\mathcal E_{1,*},H_1)\to(\mathcal E_{2,*},H_2)\)
is an isomorphism as above.  Let $S_i:(E|_X,H_i)\longrightarrow(E|_X,h)$ be the positive \(h\)-self-adjoint isometry used in \cref{def:hunitaryrep}. Set $g:=S_2 f S_1^{-1}$.
Since \(f^*H_2=H_1\), the map \(g\) is \(h\)-unitary.  Since \(\det f=\Id_P\) and \(S_1,S_2\) induce the same fixed determinant isometry $(P,h_{P,\mathrm{par}})\longrightarrow(P,h_P)$, one has $\det g=1$. Thus $g\in\mathcal G_h(X)$. Because \(f\) intertwines the two Chern connections, $A_2^h=g\cdot A_1^h$. Hence the map is well-defined.

Conversely, suppose $A_2^h=g\cdot A_1^h$ for some \(g\in\mathcal G_h(X)\). Define
\begin{align*}
f:=S_2^{-1}gS_1:
(E|_X,H_1)\longrightarrow(E|_X,H_2).
\end{align*}
Then $f^*H_2=H_1$, $\det f=\Id_P$ and $\bar\partial_2\circ f=(f\otimes\Id)\circ\bar\partial_1$ on $X$.
 Since \(f\) and \(f^{-1}\) are isometries, their
operator norms with respect to \(H_1,H_2\) are bounded.  
Moreover, by \cref{cor:curvatureLp}, $F_{H_i}\in L^p(X)$ for some $p>2$. Compatibility gives $C^{-1}h_{\mathrm{par},i}\le H_i\le C h_{\mathrm{par},i}$ and hence, for every local holomorphic section \(s\) and every
\(b\in\mathbb R\), $|s|_{H_i}=O(r^b)\Longleftrightarrow
|s|_{h_{\mathrm{par}}}=O(r^b)$. Thus the given parabolic structure \(\mathcal E_{i,*}\) is the
parabolic extension characterized by the growth condition in
\cite[Theorem~2.1]{Biquard96}.  The uniqueness statement of
\cite[Theorem~2.1]{Biquard96} therefore implies that $f:E|_X\longrightarrow E|_X$ extends uniquely to a parabolic holomorphic isomorphism $f:\mathcal E_{1,*}\longrightarrow\mathcal E_{2,*}$. The identity \(\det f=\Id_P\) extends across \(D\) as well. Thus the two parabolic Hermitian--Einstein pairs are isomorphic.

Finally, every \(A\in\mathcal A^{\mathrm{par,HE}}_h(E,L;P)\) differs from the canonical \(h\)-unitary representative of the same Chern connection by an element of \(\mathcal G_h(X)\). Hence the above map is also surjective.
\end{remark}

\begin{theorem}\label{thm:parHEtoKM}
    There is a well-defined map
\begin{align*}
    \mathfrak R:\mathcal M_h^{\mathrm{par,HE}}(E,L;\mathcal P)&\longrightarrow\mathcal M_\beta^{\mathrm{KM}}(E,L;a_{P,\mathrm{par}})\\
    [A]&\longmapsto[u\cdot A],
\end{align*}
where $u\in\Gamma(X,SU(E))$ is given in \cref{prop:globalrealization}.
\end{theorem}

\begin{proof}
Let $[A]\in\mathcal M_h^{\mathrm{par,HE}}(E,L;\mathcal P)$. By definition, \(A\) is an \(h\)-unitary representative of the Chern connection of a compatible Hermitian--Einstein metric. By \cref{prop:globalrealization}, there exists $u\in\Gamma(X,SU(E))$, such that $B:=u\cdot A\in\mathcal A_\beta^{1,p}(E,L;a_{P,\mathrm{par}})$.

By \cref{prop:hrepcurvature}, the trace-free Hermitian--Einstein equation is preserved under passing to an \(h\)-unitary representative. It is also invariant under unitary gauge transformations.  Hence $F_B^{0,2}=0$ and $\sqrt{-1}\Lambda_\omega F_B^\perp=0$. Since $\det B=a_{P,\mathrm{par}}$, \cref{lem:scalarHE} gives $(F_B^\perp)^+=0$. Therefore $[B]\in\mathcal M_\beta^{\mathrm{KM}}(E,L;a_{P,\mathrm{par}})$. \Cref{prop:realizationunique} shows that the class $[B]$ is independent of the choice of $u$.

It remains to prove independence of the representative of the
class \([A]\).  Let
\begin{align*}
A'=q\cdot A,
\qquad
q\in\mathcal G_h(X),
\end{align*}
and choose $u,u'\in\Gamma(X,SU(E))$ such that $B:=u\cdot A$, $B':=u'\cdot A'$ belong to
\(\mathcal A_\beta^{1,p}(E,L;a_{P,\mathrm{par}})\).
Set $g:=u'qu^{-1}\in\Gamma(X,SU(E))$. Then
\begin{align*}
B'=u'\cdot(q\cdot A)=g\cdot(u\cdot A)=g\cdot B.
\end{align*}
Since both \(B\) and \(B'\) belong to
\(\mathcal A_\beta^{1,p}(E,L;a_{P,\mathrm{par}})\),
\cref{prop:gauge} implies $g\in\mathcal G_\beta^{2,p}(E,L)$. Consequently $[B']=[B]\in\mathcal M_\beta^{\mathrm{KM}}(E,L;a_{P,\mathrm{par}})$.

Thus the map \(\mathfrak R\) is independent of all choices and is well-defined.
\end{proof}

\section{From twisted ASD connections to parabolic Hermitian--Einstein}\label{sec:KM-to-HE}

We now construct the inverse map from the twisted
instanton moduli space to the parabolic Hermitian--Einstein moduli space.

\begin{lemma}\label{lem:interior-smoothing}
Let $A\in\mathcal A_\beta^{1,p}(E,L;a_{P,\mathrm{par}})$ with $p>2$, and suppose $(F_A^\perp)^+=0$.
Then there exists a gauge transformation $u\in W_{\mathrm{loc}}^{2,p}\bigl(X,SU(E,h)\bigr)$
such that $\widetilde{A}:=u\cdot A$ is a smooth unitary connection on the smooth Hermitian bundle $(E|_X,h)$.
\end{lemma}

\begin{proof}
Choose a locally finite cover \(\{U_i\}_{i\in I}\) of \(X\) by relatively compact coordinate balls over which \(E\) is smoothly
trivial.  After shrinking the balls if necessary, the standard interior regularity theorem for weak anti-self-dual connections gives gauges $u_i\in W^{2,p}(U_i,SU(E,h))$ such that $A_i:=u_i\cdot A$ is smooth on \(U_i\).

Let $c_{ij}:U_i\cap U_j\to U(2)$ be the smooth transition functions of the original smooth bundle
\(E|_X\) in the chosen trivializations. On \(U_i\cap U_j\) define $c_{ij}' :=u_i c_{ij}u_j^{-1}$. Then $c_{ij}'\in W_{\mathrm{loc}}^{2,p}(U_i\cap U_j)$ and $A_i=c_{ij}'\cdot A_j$.
Equivalently, in the chosen trivializations, $dc_{ij}'=c_{ij}'A_j-A_i c_{ij}'$. Since \(A_i\) and \(A_j\) are smooth, the equation gives $c_{ij}'\in C^\infty(U_i\cap U_j)$.
The cocycle identities $c_{ij}'c_{jk}'=c_{ik}'$ hold by construction, thus \(\{c_{ij}'\}\) defines a smooth Hermitian bundle \(E'\to X\).  The smooth local connections \(A_i\) satisfy the
corresponding transformation law and therefore glue to a smooth
unitary connection \(A'\) on \(E'\).

The local maps \(u_i\) satisfy $c_{ij}'u_j=u_i c_{ij}$, and hence define a global determinant-preserving unitary bundle isomorphism $\Psi:E|_X\to E'$ of class \(W_{\mathrm{loc}}^{2,p}\), satisfying $\Psi_*A=A'$. Since \(p>2\), \(\Psi\) is in particular a continuous bundle isomorphism. By the standard smoothing argument, there exists a smooth
determinant-preserving unitary bundle isomorphism $C:E'\to E|_X$.

Set $u:=C\circ\Psi$, then $u\in W_{\mathrm{loc}}^{2,p}(X,SU(E,h))$ and $\widetilde{A}:=u\cdot A=C_*A'$ is a smooth unitary connection on the smooth bundle \(E|_X\).
\end{proof}

\begin{proposition}
\label{prop:biquardextension}
In the setting of \Cref{lem:interior-smoothing}.
Then the smooth Hermitian holomorphic bundle $(E|_X,\bar\partial_{\widetilde{A}},h)$ admits a parabolic holomorphic extension $\mathcal F_*=(\mathcal F,M,\alpha_1,\alpha_2)$ over \(\X\).

Let 
\begin{align*}
\psi:\mathcal F|_X\longrightarrow
(E|_X,\bar\partial_{\widetilde{A}})
\end{align*}
be the holomorphic identification induced by the extension, and set $H_F:=\psi^*h$. 

Then \(H_F\) is a Hermitian--Einstein metric compatible with \(\mathcal F_*\) in the sense of \Cref{def:compatibility}. Moreover, the underlying bundle pair $(F,M)$ is topologically isomorphic to $(E,L)$.
\end{proposition}

\begin{proof}
Let $A\in\mathcal A_\beta^{1,p}(E,L;a_{P,\mathrm{par}})$ satisfy $(F_A^\perp)^+=0$, and let $u\in W_{\mathrm{loc}}^{2,p}(X,SU(E,h))$, $\widetilde A:=u\cdot A$ be as in \Cref{lem:interior-smoothing}. Thus \(\widetilde A\) is a
smooth \(h\)-unitary connection on \(E|_X\). By \Cref{prop:banach}, $F_A\in
L^p$. Since \(u\) is unitary, $F_{\widetilde A}=uF_A u^{-1}\in L^p(X)$. By \Cref{lem:scalarHE}, $(F_A^\perp)^+=0$ implies $F_{\widetilde A}^{0,2}=0$ and $\ii\Lambda_\omega F_{\widetilde A}=\lambda\Id$.
Since \(\widetilde A\) is smooth, $(E|_X,\bar\partial_{\widetilde A})$ is a smooth holomorphic bundle. Thus \(\widetilde A\) is the Chern connection of the Hermitian holomorphic bundle $(E|_X,\bar\partial_{\widetilde A},h)$, whose curvature belongs to \(L^p\), with \(p>2\). Therefore
\cite[Theorem~2.1]{Biquard96} applies and gives a parabolic holomorphic extension $\mathcal F_*=(\mathcal F,M,\gamma_1,\gamma_2)$ of \((E|_X,\bar\partial_{\widetilde A})\) over \(\X\), where $0\leq\gamma_1\leq\gamma_2<1$. We next identify the weights and the underlying bundle pair.

By the definition of \(\mathcal A_\beta^{1,p}(E,L;a_{P,\mathrm{par}})\), in every
smooth \(h\)-unitary frame adapted to the smooth splitting $E|_U=Q\oplus L$, the connection \(A\) has the form $A=d_{\alpha}+c$, where $d_\alpha=d+\ii\diag(\alpha_1,\alpha_2)d\theta$, $c\in L^p$ and $r^{-1}c_{\mathrm{off}}\in L^p$. Thus these adapted frames are admissible local unitary gauges in the construction preceding \cite[Proposition~3.5]{Biquard96}. On the overlap of two such frames, the transition function is exactly the transition function of the original bundle \(E\).
Since both frames respect the splitting \(Q\oplus L\), this transition function preserves \(Q\) and \(L\).  
Consequently the topological extension induced by \(A\) is  isomorphic to \(E\), and its direct sum decomposition over \(D\) is $E|_D=Q\oplus L$, corresponding respectively to the residues \(\alpha_1\) and \(\alpha_2\). 

Since $\widetilde A=u\cdot A$ with $u\in W_{\mathrm{loc}}^{2,p}(X,SU(E,h))$, the map \(u\) is a unitary isomorphism intertwining the two
connections. By the argument preceding \cite[Proposition~3.5]{Biquard96}, \(A\) and \(\widetilde A\) have the same limiting weights, and \(u\) extends to an isomorphism between the topological extensions induced by the two connections, preserving their direct-sum decompositions over \(D\). Thus the topological extension induced by $\widetilde A$ is also $(E,Q\oplus L)$.

By \cite[Proposition~3.6]{Biquard96}, the parabolic holomorphic extension $\mathcal{F}_*$ is topologically isomorphic, together with its filtration, to the topological extension induced by its Chern connection. Applying this to
\(\widetilde A\), we obtain $(F,M)$ is topologically isomorphic to $(E,L)$ and $\gamma_i=\alpha_i$.

It remains to prove compatibility of the metric $H_F:=\psi^*h$ with $\mathcal{F}_*$. In fact, it follows from \cite[Theorem~2.1, Proposition~2.2]{Biquard96}. He starts with a Hermitian holomorphic bundle $(V,\bar\partial_V,h_V)$ over $U\setminus D$
whose Chern curvature belongs to \(L^p\), \(p>2\).  In a unitary trivialization, its Chern connection $D$ has the form $d_{\alpha}+a$. After shrinking the divisor chart, there exists a complex gauge transformation $g\in W^{1,p}_{d_\alpha}\cap C^0$ such that $g\cdot \db_D=\db_\alpha$. The proof of \cite[Theorem~2.1]{Biquard96} then constructs the
holomorphic extension by declaring $s_i:=r^{\alpha_i}g^{-1}e_i$ to be an adapted holomorphic frame, where \((e_i)\) is the above
\(h_V\)-unitary frame.

In our setting, we take $(V,\db_V,h_V):=(E|_X,\db_{\widetilde A},h)$. Since $\psi:\mathcal F|_X\to(E|_X,\db_{\widetilde A})$ is the identification induced by the extension, Biquard's result gives $\psi(s_i)=r^{\alpha_i}g^{-1}e_i$, where
$s_i$ is the holomorphic frame of $\mathcal{F}|_U$ extending over $D$ and $e_i$ is $h$-unitary frame of $E|_{U\setminus D}$. Define the local parabolic model metric \(K\) by $K(s_i,s_j)=\delta_{ij}r^{2\alpha_i}$, and take the local $K$-unitary frame $f_i:=r^{-\alpha_i}s_i$. Then with respect to the frame $(f_i)$, the relative
endomorphism $b:=K^{-1}H_F$ has the matrix $[b]_f=(g^{-1})^*g^{-1}$.

Now $g,g^{-1}\in C^0\subset L^\infty$ implies $b,b^{-1}\in L^\infty(K)$. $g\in W^{1,p}_{d_\alpha}$ implies $\nabla_{d_\alpha}g^{\pm1}\in L^p$. Since the Chern connection of \(K\) in the frame \((f_i)\) is exactly $d_\alpha$, then one has $D_Kb\in L^p$. Therefore $\db_{\End\mathcal{F}}b\in L^p\subset L^2$.
The Hermitian--Einstein condition implies $|\Lambda_\omega F_{H_F}|_{H_F}\in L^1(X)$. Thus we show $H_F$ is compatible with the local parabolic model metric \(K\). By \Cref{lem:parmodelindependence}, \(H_F\) is compatible with \(\mathcal F_*\) in the sense of
\Cref{def:compatibility}. This completes the proof.
\end{proof}

\begin{lemma}\label{lem:biquarddeterminant}
    In the setting of \Cref{prop:biquardextension}, there exists a holomorphic isomorphism
\begin{align*}
\phi_F:\det\mathcal F\longrightarrow\mathcal P.
\end{align*}
Consequently, we have $\det H_F=\phi_F^*h_{P,\mathrm{par}}$.
\end{lemma}
\begin{proof}
Let $S_P:(P|_X,h_{P,\mathrm{par}})\to(P|_X,h_P)$ be the canonical self-adjoint isometry.  On \(X\), define
\begin{align*}
\phi_F^X:=S_P^{-1}\circ\phi_E\circ\det\psi:
\det\mathcal F|_X\longrightarrow\mathcal P|_X,
\end{align*}
where $\psi:(\mathcal{F}|_X,H_F)\to (E|_X,\db_{\widetilde A},h)$ is the holomorphic identification defined in \Cref{prop:biquardextension}. 

We first prove that \(\phi_F^X\) is holomorphic. Let $A_F$ be the Chern connection of $(\mathcal{F},H_F)$. Then $H_F=\psi^*h$ gives $\psi_*A_F=\widetilde A$, and thus $(\det\psi)_*\det A_F=\det \widetilde A=\det A$. Under the fixed unitary identification \(\phi_E:\det E\to P\), the latter is the connection \(a_{P,\mathrm{par}}\) on \((P,h_P)\). By \eqref{eq:detconnection}, $a_{P,\mathrm{par}}=(S_P)_*\nabla^{h_{P,\mathrm{par}}}$, where \(\nabla^{h_{P,\mathrm{par}}}\) is the Chern connection of the fixed holomorphic line bundle $(\mathcal{P},h_{P,\mathrm{par}})$. Hence $(\phi^X_F)_*\det A_F=\nabla^{h_{P,\mathrm{par}}}$.
Taking \((0,1)\)-parts shows that \(\phi_F^X\) is holomorphic.

Moreover, \(\det\psi\), \(\phi_E\), and \(S_P^{-1}\) are isometries for the corresponding determinant metrics. Thus
\begin{align*}
\det H_F=(\phi_F^X)^*h_{P,\mathrm{par}}.
\end{align*}

It remains to extend \(\phi_F^X\) across \(D\). Fix a divisor chart \(U\), let \((s_1,s_2)\) be a holomorphic frame of \(\mathcal F|_U\) adapted to the parabolic filtration. Let \(e_P\) be a holomorphic frame of \(\mathcal P|_U\). Since \(\phi_F^X\) is a holomorphic line-bundle isomorphism on \(U\setminus D\), there exists a nowhere-vanishing holomorphic function \(f\) on \(U\setminus D\) such that $\phi_F^X(s_1\wedge s_2)=f\,e_P$.

By \Cref{prop:biquardextension}, \(H_F\) is compatible with $\mathcal{F}_*$. Hence $|s_1\wedge s_2|_{\det H_F}\asymp r^{\alpha_1+\alpha_2}$. On the other hand, by the definition \(h_{P,\mathrm{par}}=\|\sigma\|^{2(\alpha_1+\alpha_2)}h_P\), one also has $|e_P|_{h_{P,\mathrm{par}}}\asymp r^{\alpha_1+\alpha_2}$. Thus both \(f\) and \(f^{-1}\) are bounded holomorphic functions on \(U\setminus D\), hence they extend holomorphically across \(D\). Thus one has a global holomorphic isomorphism $\phi_F:\det\mathcal F\longrightarrow\mathcal P$.
\end{proof}

\begin{theorem}\label{thm:inverse-map}
    There exists a well-defined map
    \begin{align*}
        \mathfrak S:\mathcal M_\beta^{\mathrm{KM}}(E,L;a_{P,\mathrm{par}})\longrightarrow\mathcal M_h^{\mathrm{par,HE}}(E,L;\mathcal P)
    \end{align*}
    satisfying $\mathfrak S=\mathfrak R^{-1}$, which is defined in \Cref{thm:parHEtoKM}.

    Consequently, there exists a bijection between two moduli spaces
    \begin{align*}
        \mathfrak R:\mathcal M_h^{\mathrm{par,HE}}(E,L;\mathcal P)\xrightarrow{\;\sim\;}\mathcal M_\beta^{\mathrm{KM}}(E,L;a_{P,\mathrm{par}}).
    \end{align*}
\end{theorem}
\begin{proof}
Let $[A]\in\mathcal M_\beta^{\mathrm{KM}}(E,L;a_{P,\mathrm{par}})$ and choose a representative
\(A\in\mathcal A_\beta^{1,p}(E,L;a_{P,\mathrm{par}})\)
satisfying \((F_A^\perp)^+=0\).  By
\Cref{lem:interior-smoothing}, there exists
\(u\in W_{\mathrm{loc}}^{2,p}(X,SU(E,h))\) such that
\(\widetilde A:=u\cdot A\) is smooth on \(X\).
By \Cref{prop:biquardextension}, the Hermitian holomorphic bundle
\((E|_X,\bar\partial_{\widetilde A},h)\) admits a parabolic
holomorphic extension
\(\mathcal F_*=(\mathcal F,M,\alpha_1,\alpha_2)\) with a compatible
Hermitian--Einstein metric \(H_F\), and the underlying bundle pair
\((F,M)\) is topologically isomorphic to \((E,L)\).  By
\Cref{lem:biquarddeterminant}, there is a holomorphic isomorphism
\(\phi_F:\det\mathcal F\to\mathcal P\) satisfying
\(\det H_F=\phi_F^*h_{P,\mathrm{par}}\).

Since \((F,M)\) and \((E,L)\) are topologically isomorphic,
\(c_2(F)=c_2(E)\) and \(\deg M=\deg L\).  Hence
\Cref{lem:pairclassification} gives a smooth determinant-preserving bundle-pair isomorphism $C:(F,M,\phi_F)\longrightarrow(E,L,\phi_E)$.
Set
\begin{align*}
    \db_{E,A}:=C_*\db_F,\qquad H_A:=C_*H_F,
\end{align*}
then the parabolic bundle $\mathcal E_{A,*}:=(E,\bar\partial_{E,A},L,\alpha_1,\alpha_2)$
has determinant \(\mathcal P\), and \(H_A\) is a
compatible Hermitian--Einstein metric.  We define
\begin{align*}
\mathfrak S([A]):=[(\mathcal E_{A,*},H_A)].
\end{align*}

We first prove that this definition is independent of all choices. Suppose \([A_0]=[A_1]\in \M_\beta^{\mathrm{KM}}\).  Then
\(A_1=g\cdot A_0\) for some
\(g\in\mathcal G_\beta^{2,p}(E,L)\).  Choose smoothing gauges
\(u_i\in W_{\mathrm{loc}}^{2,p}(X,SU(E,h))\) and put
\(\widetilde A_i=u_i\cdot A_i\).  The map
\(q:=u_1gu_0^{-1}\) is a determinant-one unitary isomorphism satisfying $\nabla_{\widetilde A_1}\circ q=(q\otimes\Id)\circ\nabla_{\widetilde A_0}$.
Since \(\widetilde A_0\) and \(\widetilde A_1\) are smooth, the identity
\(dq=q\widetilde A_0-\widetilde A_1q\) and Sobolev bootstrapping imply
\(q\in C^\infty(X,SU(E,h))\).  In particular, \(q\) is a holomorphic
isometry
\begin{align*}
q:(E|_X,\bar\partial_{\widetilde A_0},h)
\longrightarrow
(E|_X,\bar\partial_{\widetilde A_1},h).
\end{align*}
Both \(q\) and \(q^{-1}\) are bounded.  The uniqueness assertion in
\cite[Theorem~2.1]{Biquard96} therefore extends \(q\) uniquely to a parabolic holomorphic isomorphism between their parabolic extensions.  Since \(q\) is an isometry on \(X\), this extension
identifies the metrics \(H_{F,0}\) and \(H_{F,1}\).  Since
\(\det q=1\), it also identifies the determinant isomorphisms supplied
by \Cref{lem:biquarddeterminant}.  After transport to the fixed pair
\((E,L,\phi_E)\), the two constructions therefore determine the same
element of \(\mathcal M_h^{\mathrm{par,HE}}(E,L;\mathcal P)\).
Taking \(g=\Id\) proves independence of the smoothing gauge $u_i$.  A
different choice of the determinant-preserving pair isomorphism
\(C\) changes the transported object by a smooth determinant-preserving
bundle-pair automorphism, and hence does not change its moduli class.
Thus \(\mathfrak S\) is well-defined.

We next prove \(\mathfrak R\circ\mathfrak S=\Id\). Let \([A]\in\mathcal M_\beta^{\mathrm{KM}}\). Define
\begin{align*}
T:=\psi\circ C^{-1}:(E|_X,H_A)\longrightarrow(E|_X,h).
\end{align*}
Since \(H_A=C_*H_F\) and \(H_F=\psi^*h\), \(T\) is an isometry and $T_*\nabla^{H_A}=\widetilde A$.

We claim that \(\widetilde A\) is an \(h\)-unitary representative of \(\nabla^{H_A}\). By the definition of \(\phi_F\) in \Cref{lem:biquarddeterminant},
\begin{align*}
\phi_E\circ\det\psi=S_P\circ\phi_F
.
\end{align*}
Since \(C\) is determinant-preserving, \(\phi_E\circ\det C=\phi_F\), and hence
\begin{align*}
\phi_E\circ\det T
=\phi_E\circ\det\psi\circ\det C^{-1}
=S_P\circ\phi_E.
\end{align*}
Moreover, \(\det H_A=\phi_E^*h_{P,\mathrm{par}}\) and \(\det h=\phi_E^*h_P\). Therefore the canonical positive isometry
\(\widehat S_{H_A}:(E|_X,H_A)\to(E|_X,h)\) satisfies
\begin{align*}
\phi_E\circ\det\widehat S_{H_A}
=S_P\circ\phi_E.
\end{align*}
Thus $\det T=\det\widehat S_{H_A}$. Therefore \(\widetilde A=T_*\nabla^{H_A}\) is an \(h\)-unitary representative of \(\nabla^{H_A}\).

By \Cref{thm:parHEtoKM}, there exists
\(w\in\Gamma(X,SU(E,h))\) such that $B:=w\cdot\widetilde A
\in\mathcal A_\beta^{1,p}$, and $\mathfrak R\circ\mathfrak S([A])=[B]$.
Since \(\widetilde A=u\cdot A\), one has
\begin{align*}
B=(wu)\cdot A.
\end{align*}
The connections \(A\) and \(B\) both belong to
\(\mathcal A_\beta^{1,p}(E,L;a_{P,\mathrm{par}})\), and
\(wu\in W_{\mathrm{loc}}^{1,1}(X,SU(E,h))\). Hence \Cref{prop:gauge} implies
\begin{align*}
wu\in\mathcal G_\beta^{2,p}(E,L).
\end{align*}
Thus \([B]=[A]\), and therefore
\begin{align*}
\mathfrak R(\mathfrak S([A]))=[A].
\end{align*}

Finally, we prove \(\mathfrak S\circ\mathfrak R=\Id\). Let $[(\mathcal E_*,H)]\in
\mathcal M_h^{\mathrm{par,HE}}$. Choose an \(h\)-unitary representative $A_0=T_*\nabla^H$, where $T:(E|_X,H)\to (E|_X,h)$.
By \Cref{thm:parHEtoKM}, there exists
\(v\in\Gamma(X,SU(E,h))\) such that
\begin{align*}
A:=v\cdot A_0\in\mathcal A_\beta^{1,p}(E,L;a_{P,\mathrm{par}}),\qquad \mathfrak R([(\mathcal{E}_*,H)])=[A].
\end{align*}
Apply \Cref{lem:interior-smoothing} to \(A\). There exists
\(u\in W_{\mathrm{loc}}^{2,p}(X,SU(E,h))\) such that
\(\widetilde A:=u\cdot A\) is smooth on \(X\). Set
\begin{align*}
Q:=uvT:(E|_X,H)\longrightarrow(E|_X,h).
\end{align*}
Then
\begin{align*}
Q_*\nabla^H=u_*(v_*A_0)=u\cdot A=\widetilde A.
\end{align*}
Since both \(\nabla^H\) and \(\widetilde A\) are smooth, elliptic bootstrapping gives
\(Q\in C^\infty(X)\). Taking \((0,1)\)-parts yields
 \(Q\) is a holomorphic isometry $Q:(\mathcal E|_X,H)\to (E|_X,\bar\partial_{\widetilde A},h)$.
Let $\psi:(\mathcal F|_X,H_F)\to (E|_X,\bar\partial_{\widetilde A},h)$ be the holomorphic isometry in
\Cref{prop:biquardextension}, and define
\begin{align*}
\Phi:=\psi^{-1}\circ Q:
(\mathcal E|_X,H)\longrightarrow(\mathcal F|_X,H_F).
\end{align*}
Then \(\Phi\) is a holomorphic isometry.

Since \(u,v\in SU(E,h)\), one has \(\det Q=\det T\). Since \(A_0=T_*\nabla^H\) is an \(h\)-unitary representative, $\phi_E\circ\det T=S_P\circ\phi_E$.
Hence $\phi_E\circ\det Q=S_P\circ\phi_E$. By the construction of \(\phi_F\) in \Cref{lem:biquarddeterminant}, $\phi_E\circ\det\psi=S_P\circ\phi_F$.
Since \(\Phi=\psi^{-1}Q\), one has $\phi_F\circ\det\Phi=\phi_E$.

Since $\Phi$ is an isometry and both metrics are compatible with corresponding parabolic structures,
\cite[Theorem 2.1]{Biquard96} extends \(\Phi\) uniquely to a parabolic holomorphic isomorphism
\begin{align*}
\overline\Phi:\mathcal E_*\xrightarrow{\;\sim\;}\mathcal F_*.
\end{align*}
Moreover, $\phi_F\circ\det\overline\Phi=\phi_E$ by extending across $D$.

Let $C:(F,M,\phi_F)\to (E,L,\phi_E)$ be the smooth bundle pair isomorphism used in the
definition of \(\mathfrak S([A])\), and let
\(\mathcal E_{A,*}:=C_*\mathcal F_*\) and \(H_A:=C_*H_F\).
Then
\begin{align*}
C\circ\overline\Phi:
(\mathcal E_*,H)\xrightarrow{\;\sim\;}(\mathcal E_{A,*},H_A)
\end{align*}
is a determinant-preserving holomorphic parabolic isometry. Hence $\mathfrak S\circ\mathfrak R([(\mathcal{E}_*,H)])=[(\mathcal{E}_*,H)]$.

Thus $\mathfrak R\circ\mathfrak S=\Id$ and $\mathfrak S\circ\mathfrak R=\Id$.
Therefore \(\mathfrak S=\mathfrak R^{-1}\), and
\begin{align*}
\mathfrak R:
\mathcal M_h^{\mathrm{par,HE}}(E,L;\mathcal P)
\xrightarrow{\sim}
\mathcal M_\beta^{\mathrm{KM}}(E,L;a_{P,\mathrm{par}})
\end{align*}
is a bijection.
\end{proof}

\section{Parabolic polystability and Hermitian--Einstein metrics}
\label{sec:parabolic-polystability}

We have already established a correspondence between compatible parabolic Hermitian--Einstein connections and determinant-fixed twisted ASD connections. In this section, by parabolic Kobayashi--Hitchin correspondence, we identify the moduli space of polystable parabolic bundles with these two moduli spaces we considered.

\begin{definition}\label{def:polystable-moduli}
Let $\varphi_E:\det E\to P$ be the fixed smooth determinant
identification. Define
$\M_{\mathrm{pst}}^{\mathrm{par}}(E,L;\mathcal P)$
to be the set of isomorphism classes of pairs
$(\mathcal F_*,\varphi)$ satisfying the following conditions:

\begin{enumerate}
\item
$\mathcal F_*=(\mathcal F,\mathcal M,\alpha_1,\alpha_2)$ is a rank-two parabolic polystable holomorphic bundle on $(\X,D)$.

\item
$\varphi:\det\mathcal F\xrightarrow{\sim}\mathcal P$
is a holomorphic line bundle isomorphism.

\item Let $(F,M)$ be the underlying smooth complex bundles of $(\mathcal F,\mathcal{M})$. There exists a smooth complex bundle isomorphism $\Psi:F\to E$ satisfying $\Psi(M)=L$ and $\varphi_E\circ\det\Psi=\varphi$.
\end{enumerate}

Two pairs $(\mathcal F_*,\varphi)$ and $(\mathcal F'_*,\varphi')$ are isomorphic if there exists a holomorphic bundle isomorphism $f:\mathcal F\xrightarrow{\sim}\mathcal F'$ such that $f(\mathcal M)=\mathcal M'$ and $\varphi'\circ\det f=\varphi$.
\end{definition}

\begin{proposition}\label{prop:HE-implies-polystable}
Let $\mathcal E_*=(\mathcal E,\mathcal L,\alpha_1,\alpha_2)$ be a rank-two parabolic holomorphic bundle on $(\X,D)$. If $\mathcal E|_X$ admits a compatible Hermitian--Einstein metric, then $\mathcal E_*$ is polystable.
\end{proposition}

\begin{proof}
Let $H$ be a compatible Hermitian--Einstein metric.
If $\mathcal E$ is indecomposable, then
$\mathcal E_*$ is stable by \cite[Theorem~2.14]{LiTian}.

Suppose that $\mathcal E$ is decomposable and not stable. Then the Chern--Weil argument in the proof of
\cite[Theorem~2.14]{LiTian}, together with \cite[Proposition~3.9]{LiTian}, shows that
$\mathcal E_*$ is semistable and that there exists
a saturated rank-one subsheaf $\mathcal S\subset\mathcal E$ satisfying $\pardeg_\omega(\mathcal S_*)=\parmu_\omega(\mathcal E_*)$. Let $S:=\mathrm{Sing}(\mathcal S)$ be the singular set of $\mathcal S$, which is a finite subset of $\X$.

Let $\pi$ be the $H$-orthogonal projection onto
$\mathcal S$ on $X\setminus S$. The equality case of the Chern--Weil formula gives $\bar\partial\pi=0$. Since $\pi$ is $H$-self-adjoint, it is parallel with respect to the Chern connection of $H$. The holomorphic projection $\pi$ extends across $S\cap X$ by Hartogs' theorem.
Moreover, $\|\pi\|_H=1$, and compatibility of $H$
implies that $\pi$ has bounded operator norm
with respect to the parabolic model metric.
Thus $\pi$ extends across $D$ to a holomorphic endomorphism $\widetilde\pi$ of $\mathcal E$ preserving $\mathcal L$.

The identities $\widetilde\pi^2=\widetilde\pi$ and
$\operatorname{tr}\widetilde\pi=1$ follow from continuity. Consequently, $\mathcal E
=\im\widetilde\pi\oplus\ker\widetilde\pi$ is a decomposition into holomorphic line bundles.
Since $\widetilde\pi$ preserves $\mathcal L$,
and $D$ is connected, this decomposition
respects the parabolic structure.

The image of $\widetilde\pi$ agrees with
$\mathcal S$ generically and hence everywhere
by saturation. Thus its induced parabolic degree
equals $\parmu_\omega(\mathcal E_*)$.
Additivity of parabolic degree shows that the
other summand has the same parabolic slope.
Therefore $\mathcal E_*$ is a direct sum of
two parabolic line bundles of equal slope,
and hence is polystable.
\end{proof}

\begin{theorem}[Parabolic Kobayashi--Hitchin correspondence]
\label{thm:parabolic-KH}
The forgetful map
\begin{align*}
    \mathfrak F:\M_h^{\mathrm{par,HE}}(E,L;\mathcal P)&\xrightarrow{\;\sim\;} \M_{\mathrm{pst}}^{\mathrm{par}}(E,L;\mathcal P),\\
[(\mathcal E_*,H)]&\longmapsto[(\mathcal E_*,\varphi_E)]
\end{align*}
is a bijection.
\end{theorem}

\begin{proof}
The map is well defined by \cref{prop:HE-implies-polystable}. We first prove surjectivity. Let $(\mathcal F_*,\varphi)$ represent an element of
$\M_{\mathrm{pst}}^{\mathrm{par}}(E,L;\mathcal P)$.
If $\mathcal F_*$ is stable, the existence theorem
\cite[Theorem~2.14]{LiTian} provides a compatible
Hermitian--Einstein metric.
Otherwise, polystability gives a decomposition
\begin{align*}
\mathcal F_*
\simeq \mathcal S_{1,*}\oplus\mathcal S_{2,*},
\qquad
\pardeg_\omega(\mathcal S_{1,*})
=\pardeg_\omega(\mathcal S_{2,*}).
\end{align*}
Applying the scalar construction of
\cref{prop:detnormalization} to each parabolic line
bundle gives compatible Hermitian--Einstein metrics
$H_1,H_2$. Their Einstein constants are
\begin{align*}
\lambda_i=\frac{2\pi}{\vol_\omega(\X)}
  \pardeg_\omega(\mathcal S_{i,*}),
\end{align*}
and hence coincide. Thus $H_1\oplus H_2$ is a
compatible Hermitian--Einstein metric on $\mathcal F_*$. The smooth bundle-pair identification $\Psi$ in \cref{def:polystable-moduli} then transports $(\mathcal F_*,H,\varphi)$ to an
element of $\M_h^{\mathrm{par,HE}}(E,L;\mathcal P)$.
This proves surjectivity.

We next prove injectivity. Suppose two normalized compatible Hermitian--Einstein metrics $H_0,H_1$ are defined on the same parabolic bundle $\mathcal E_*$,
with $\det H_0=\det H_1=\varphi_E^*h_{P,\mathrm{par}}$. In particular, their Einstein constants coincide.

Consider the parabolic bundle $\mathcal V_*=\mathcal E_*\oplus\mathcal E_*$ equipped with the compatible Hermitian--Einstein metric $K=H_0\oplus H_1$. Let $\Gamma\subset\mathcal V$ be the graph of the
identity map of $\mathcal E$. It is a holomorphic parabolic subbundle with $\Gamma_*\simeq\mathcal E_*$.

Let $\pi$ be the $K$-orthogonal projection onto
$\Gamma$. The degree identities \cite[Propositions~3.4 and~3.9]{LiTian}, together with the Chern--Weil projection formula, give
\begin{align*}
\pardeg_\omega(\Gamma_*)=\pardeg_\omega(\mathcal E_*)
-\frac{1}{2\pi}\int_X|\bar\partial\pi|_K^2\,dV_\omega.
\end{align*}
Hence $\bar\partial\pi=0$.
As $\pi$ is self-adjoint, it is parallel with
respect to the Chern connection of $K$.
Consequently, the graph of the identity is
parallel, which implies $\nabla^{H_0}=\nabla^{H_1}$ on $X$.

Therefore $s:=H_0^{-1}H_1$ is a positive parallel
endomorphism with $\det s=1$.
The automorphism $f:=s^{-1/2}$ is holomorphic on
$X$, satisfies $\det f=1$ and $f^*H_1=H_0$.
Compatibility implies that $f$ and $f^{-1}$ are
bounded with respect to the parabolic model metric.
By the uniqueness of compatible parabolic extension
\cite[Theorem~2.1]{Biquard96}, $f$ extends to a
holomorphic parabolic automorphism of $\mathcal E_*$.
Thus the two Hermitian--Einstein pairs are
isomorphic, proving injectivity.
\end{proof}

\section{Analytic compactification of parabolic Hermitian–Einstein moduli}\label{sec:ideal}

\subsection{Ideal connections and convergence}
\label{subsec:ideal}

The compactness theorem and the correspondence theorem shows that a sequence of parabolic Hermitian--Einstein connections on a fixed bundle pair $(E,L)$ may converge, away from finitely many points, to a parabolic Hermitian--Einstein connection on a possibly different bundle pair $(F,M)$.
At a concentration point $x\in X$, the local defect has the form $(m_x,0)$, where $m_x$ is a positive integer. At a concentration point $y\in D$, the defect is a pair $(k_y,l_y)\in\ZZ^2$ satisfying
\begin{align*}
    k_y\geq 0,\qquad k_y+l_y\geq 0.
\end{align*}
Consequently, both $k_y$ and $k_y+l_y$ are nonnegative integers. The corresponding curvature-energy loss is
\begin{align*}
    k_y+2\beta l_y
    =(1-2\beta)k_y+2\beta(k_y+l_y).
\end{align*}
This suggests recording the two nonnegative multiplicities separately. Away from $D$, they both equal the ordinary instanton multiplicity $m_x$. 
We use two effective zero-cycles to encode these
multiplicities.

\begin{definition}
\label{def:effectivecycles}
An \textbf{effective zero-cycle} on $\X$ is a finite formal sum $Z=\sum_{x\in\X}n_x x$ with $n_x\in\ZZ_{\ge 0}$,
where only finitely many coefficients $n_x$ are nonzero.

Its \textbf{support} and \textbf{degree} are defined by
\begin{align*}
    \operatorname{supp} Z
    :=\{x\in\X:n_x>0\},
    \qquad
    \deg Z:=\sum_{x\in\X}n_x.
\end{align*}
For $x\in\X$, the \textbf{multiplicity} is $\operatorname{mult}_x Z:=n_x$.
For a subset $U\subset\X$, define $Z|_U:=\sum_{x\in U}n_x x$.
The \textbf{counting measure associated with $Z$} is
\begin{align*}
    [Z]:=\sum_{x\in\X}n_x\delta_x.
\end{align*}

For $d\in\mathbb Z_{\geq0}$, the space of effective zero-cycles of degree $d$ is exactly $\operatorname{Sym}^d(\X)$.

\end{definition}

\begin{definition}\label{def:idealparHE}

An \textbf{ideal parabolic Hermitian--Einstein pair of initial type $(E,L;\mathcal P)$} is an equivalence class of data $(\mathcal F_*,H,\phi_F,Z_0,Z_1)$ satisfying the following conditions:
\begin{enumerate}
\item $\mathcal F_*=(\mathcal F,M,\alpha_1,\alpha_2)$
is a rank-two parabolic holomorphic bundle on $(\X,D)$. Let $F$ denote its underlying smooth complex vector bundle.
\item $\phi_F:\det\mathcal F\xrightarrow{\sim}\mathcal P$ is a holomorphic line-bundle isomorphism.
\item $H$ is a smooth Hermitian--Einstein metric on
$\mathcal F|_X$, compatible with $\mathcal F_*$,
and normalized by $\det H=\phi_F^*h_{P,\mathrm{par}}$.
\item $Z_0$ and $Z_1$ are effective zero-cycles on
$\X$ satisfying
\begin{align*}
    Z_0|_X=Z_1|_X,\quad c_2(E)=c_2(F)+\deg Z_0,\quad \deg L&=\deg M+\deg Z_0-\deg Z_1,
\end{align*}
where $\deg M:=\langle c_1(M),[D]\rangle$.
\end{enumerate}

Two such data
$(\mathcal F_*,H,\phi_F,Z_0,Z_1)$ and
$(\mathcal F'_*,H',\phi_{F'},Z'_0,Z'_1)$ are \textbf{equivalent} if $Z_0=Z'_0$, $Z_1=Z_1'$ 
and there exists a holomorphic bundle isomorphism $f:\mathcal F\xrightarrow{\sim}\mathcal F'$ such that
\begin{align*}
    f(M)=M',\quad \phi_{F'}\circ\det f&=\phi_F,\quad f^*H'=H.
\end{align*}

The set of equivalence classes is denoted by $\mathfrak I^{\mathrm{par,HE}}_\alpha(E,L;\mathcal P)$. We write its elements as $\mathfrak a=\bigl([(\mathcal F_*,H,\phi_F)],Z_0,Z_1\bigr)$.
\end{definition}

\begin{definition}
\label{def:totalmeasure}
Let $\mathfrak a=\bigl([(\mathcal F_*,H,\phi_F)],Z_0,Z_1\bigr)\in\mathfrak I^{\mathrm{par,HE}}_\alpha(E,L;\mathcal P)$.
Define the \textbf{total curvature measure} on $\X$ of $\mathfrak a$ by
\begin{align*}
    \nu_{\mathfrak a}:=\frac{1}{8\pi^2}|F_H^\perp|_H^2\,dV_\omega+(1-2\beta)[Z_0]+2\beta[Z_1].
\end{align*}
\end{definition}

Similar to the standard compactification theory of instanton moduli spaces, we can define the convergence among ideal twisted ASD connections. Then using the correspondence established in \Cref{sec:HE-to-KM} and \Cref{sec:KM-to-HE}, we can define the convergence among ideal parabolic Hermitian--Einstein pairs. Kronheimer--Mrowka proved the compactness theorem for twisted ASD $SU(2)$-connections, but they didn't give an explicit construction of the compactified space for the applications in their papers. We supply the details here. 

There is a subtlety between standard ASD connections and twisted ASD connection: the total curvature measure convergence is not enough to record the bubbling, since the energy defect $k+2\beta l$ cannot determine the pair $(k,l)$.

\begin{definition}\label{def:idealconv}

Let $\mathfrak a_n=\bigl([(\mathcal F_{n,*},H_n,\phi_{F_n})],Z_{0,n},Z_{1,n}\bigr)$ and $\mathfrak a=\bigl([(\mathcal F_*,H,\phi_F)],Z_0,Z_1\bigr)$ be elements of
$\mathfrak I^{\mathrm{par,HE}}_\alpha(E,L;\mathcal P)$.
Put $[A_n]:=\mathfrak R\bigl([(\mathcal F_{n,*},H_n,\phi_{F_n})]\bigr)$, $[A]:=\mathfrak R
\bigl([(\mathcal F_*,H,\phi_F)]\bigr)$ and $S:=\operatorname{supp}(Z_0+Z_1)$.

We say that $\mathfrak a_n$ \textbf{converges} to $\mathfrak a$ if the following conditions hold:
\begin{enumerate}
\item
For sufficiently large $n$, there exist representatives $A_n,A$ and bundle-pair isomorphisms
\begin{align*}
\Psi_n:
(F,M,\phi_F)|_{\X\setminus S}
\longrightarrow
(F_n,M_n,\phi_{F_n})|_{\X\setminus S},
\end{align*}
such that $\Psi_n^*A_n\to A$ in weighted $W^{1,p}$ on every relatively compact subsets of $\X\setminus S$.
\item
The total curvature measures satisfy $\nu_{\mathfrak a_n}\rightharpoonup\nu_{\mathfrak a}$ as finite Borel measures on $\X$.
\item
For every sufficiently small ball $B\subset\X$ whose boundary is disjoint from $\operatorname{supp}S$, one has, for sufficiently large $n$,
\begin{align*}
    \deg(Z_0|_B)&=\deg(Z_{0,n}|_B)+k_B(\Psi_n),\\
\deg(Z_1|_B)&=\deg(Z_{1,n}|_B)+k_B(\Psi_n)+l_B(\Psi_n),
\end{align*}
where $(k_B(\Psi_n),l_B(\Psi_n))$ are the relative
local labels defined in \Cref{def:relative}.
\end{enumerate}

\end{definition}

\begin{proposition}[Sequential compactness]
\label{prop:idealsequential}
Every sequence in $\mathfrak I^{\mathrm{par,HE}}_\alpha(E,L;\mathcal P)$ has a convergent subsequence whose limit also belongs to $\mathfrak I^{\mathrm{par,HE}}_\alpha(E,L;\mathcal P)$.
\end{proposition}

\begin{proof}
Let $\mathfrak a_n=\bigl([(\mathcal F_{n,*},H_n,\phi_{F_n})],Z_{0,n},Z_{1,n}\bigr)$
be a sequence in $\mathfrak I^{\mathrm{par,HE}}_\alpha(E,L;\mathcal P)$. Let $(F_n,M_n)$ be the underlying smooth bundle pairs.

By the singular Chern--Weil formula \Cref{cor:tracefreeaction} and the topological constraints in \Cref{def:idealparHE},
\begin{align*}
    \kappa_\beta(E,L)&=\frac{1}{8\pi^2}\int_X|F_{H_n}^\perp|_{H_n}^2\,dV_\omega+(1-2\beta)\deg Z_{0,n}+2\beta\deg Z_{1,n}.
\end{align*}
Since $0<\beta<1/2$, the degrees of both cycles are
uniformly bounded.
After passing to a subsequence, there exist
$d_0,d_1\in\ZZ_{\geq0}$ such that $\deg Z_{i,n}=d_i$ for all $n$ and $i=0,1$.
By compactness of $\operatorname{Sym}^{d_i}(\X)$,
we may further assume that $Z_{i,n}\to Y_i$ in $\operatorname{Sym}^{d_i}(\X)$. The condition $Z_{0,n}|_X=Z_{1,n}|_X$ passes to the limit, so we have $Y_0|_X=Y_1|_X$.

The topological constraints also give
\begin{align*}
c_2(F_n)=c_2(E)-d_0,\qquad \deg M_n=\deg L+d_1-d_0.
\end{align*}
By \Cref{lem:pairclassification}, after choosing
determinant-preserving smooth bundle-pair
isomorphisms, we may identify all
$(F_n,M_n,\phi_{F_n})$ with a fixed smooth
determinant-fixed pair $(F^0,M^0,\phi_0)$.

Put $[A_n]:=\mathfrak R\bigl([(\mathcal F_{n,*},H_n,\phi_{F_n})]\bigr)$, and transport these weighted instanton classes to $(F^0,M^0,\phi_0)$.
Applying the compactness theorem \Cref{thm:detcompact} and passing to a further
subsequence, we obtain a smooth determinant-fixed
pair $(F_\infty,M_\infty,\phi_\infty)$, a weighted
instanton class
\begin{align*}
[A_\infty]\in\mathcal M_\beta^{\mathrm{KM}}
(F_\infty,M_\infty;a_{P,\mathrm{par}}),
\end{align*}
a finite concentration set $T=T_X\sqcup T_D$, and comparison isomorphisms
\begin{align*}
\Psi_n:
(F_\infty,M_\infty,\phi_\infty)|_{\X\setminus T}
\longrightarrow
(F_n,M_n,\phi_{F_n})|_{\X\setminus T},
\end{align*}
satisfying the convergence asserted in that theorem.

Let $m_x$ be the ordinary defect at $x\in T_X$,
and let $(k_y,l_y)$ be the divisor defect at
$y\in T_D$. Define
\begin{align*}
    W_0:=\sum_{x\in T_X}m_xx+\sum_{y\in T_D}k_yy,\qquad W_1:=\sum_{x\in T_X}m_xx+\sum_{y\in T_D}(k_y+l_y)y.
\end{align*}
The admissibility inequalities imply that $W_0$ and
$W_1$ are effective zero-cycles, and
$W_0|_X=W_1|_X$.

Define $Z_i:=Y_i+W_i$, $i=0,1$. By the inverse correspondence, let
\begin{align*}
    [(\mathcal F_{\infty,*},H_\infty,\phi_\infty)]:=\mathfrak S([A_\infty]).
\end{align*}
We claim that $\mathfrak a:=\bigl(
[(\mathcal F_{\infty,*},H_\infty,\phi_\infty)],
Z_0,Z_1\bigr)$ belongs to $\mathfrak I^{\mathrm{par,HE}}_\alpha(E,L;\mathcal P)$
and that $\mathfrak a_n\to\mathfrak a$.

First, $Z_0|_X=Z_1|_X$. By \Cref{thm:detcompact},
\begin{align*}
    c_2(F_\infty)=c_2(F^0)-\deg W_0,\qquad \deg M_\infty=\deg M^0+\deg W_1-\deg W_0.
\end{align*}
Since
\begin{align*}
c_2(E)=c_2(F^0)+d_0,\qquad \deg L=\deg M^0-d_1+d_0,
\end{align*}
and $\deg Y_i=d_i$, it follows that
\begin{align*}
    c_2(E)=c_2(F_\infty)+\deg Z_0,\qquad \deg L=\deg M_\infty+\deg Z_0-\deg Z_1.
\end{align*}
Thus $\mathfrak a\in \mathfrak I^{\mathrm{par,HE}}_\alpha(E,L;\mathcal P)$.

Set $S:=\operatorname{supp}(Z_0+Z_1)$. Since $T\subset S$, the comparison isomorphisms $\Psi_n$ supplied by \Cref{thm:detcompact} restrict
to $\X\setminus S$. Their weighted $W^{1,p}$
convergence and smooth convergence away from $D$
give condition \textup{(1)} of \Cref{def:idealconv}.

For condition \textup{(2)}, the correspondence
$\mathfrak R$ preserves the pointwise norm of the
trace-free curvature. Hence the curvature measure
convergence in \Cref{thm:detcompact} gives
\begin{align*}
    \frac{1}{8\pi^2}|F_{H_n}^\perp|_{H_n}^2\,dV_\omega
\rightharpoonup\;\frac{1}{8\pi^2}|F_{H_\infty}^\perp|_{H_\infty}^2\,dV_\omega+(1-2\beta)[W_0]+2\beta[W_1].
\end{align*}
On the other hand, convergence of the cycles gives
$[Z_{i,n}]\rightharpoonup[Y_i]$.
Therefore
\begin{align*}
    \nu_{\mathfrak a_n}\rightharpoonup\;\frac{1}{8\pi^2}|F_{H_\infty}^\perp|_{H_\infty}^2+(1-2\beta)[Y_0+W_0]+2\beta[Y_1+W_1]=\nu_{\mathfrak a}.
\end{align*}

Finally, let $B$ be a sufficiently small adapted ball satisfying $\dd B\cap S=\varnothing$. Since $Z_{i,n}\to Y_i$ in $\operatorname{Sym}^{d_i}(\X)$ and $\partial B\cap\operatorname{supp}Y_i=\varnothing$,
we have, for sufficiently large $n$, $\deg(Z_{i,n}|_B)=\deg(Y_i|_B)$.

By the definition of the local defect labels and their additivity
under excision in \Cref{lem:localintegral},
the relative labels of $\Psi_n$ on $\partial B$
satisfy
\begin{align*}
    k_B(\Psi_n)=\deg(W_0|_B),\qquad k_B(\Psi_n)+l_B(\Psi_n)=\deg(W_1|_B)
\end{align*}
for sufficiently large $n$.

Consequently,
\begin{align*}
    &\deg(Z_{0,n}|_B)+k_B(\Psi_n)=\deg(Y_0|_B)+\deg(W_0|_B)=\deg(Z_0|_B),\\
    &\deg(Z_{1,n}|_B)+k_B(\Psi_n)+l_B(\Psi_n)=\deg(Y_1|_B)+\deg(W_1|_B)=\deg(Z_1|_B).
\end{align*}
This proves condition \textup{(3)}. Therefore $\mathfrak a_n\to\mathfrak a$ in the sense
of \Cref{def:idealconv}.
\end{proof}

\begin{proposition}[Uniqueness of ideal limits]
\label{prop:uniquelimit}
A sequence in $\mathfrak I^{\mathrm{par,HE}}_\alpha(E,L;\mathcal P)$ has at most one limit in the sense of
\Cref{def:idealconv}.
\end{proposition}

\begin{proof}
Suppose that a sequence $\mathfrak a_n$ converges to
both $\mathfrak a=\bigl([(\mathcal F_*,H,\phi_F)],Z_0,Z_1\bigr)$ and $\mathfrak b=\bigl([(\mathcal G_*,K,\phi_G)],Y_0,Y_1\bigr)$.
Let $(F,M,\phi_F)$ and $(G,N,\phi_G)$ be the
underlying determinant-fixed smooth bundle pairs.
By condition \textup{(2)} of \Cref{def:idealconv}, $\nu_{\mathfrak a}=\nu_{\mathfrak b}$. Consider the atomic part of these measures, then one has $\operatorname{supp}(Z_0+Z_1)
=\operatorname{supp}(Y_0+Y_1)=:S$.

Let $[A]=\mathfrak R([(\mathcal F_*,H,\phi_F)])$, $[B]=\mathfrak R([(\mathcal G_*,K,\phi_G)])$.
By condition \textup{(1)} of
\Cref{def:idealconv}, there are bundle-pair isomorphisms
\begin{align*}
    \Psi_n:(F,M,\phi_F)|_{\X\setminus S}\to(F_n,M_n,\phi_{F_n})|_{\X\setminus S},\qquad
    \Phi_n:(G,N,\phi_G)|_{\X\setminus S}\to(F_n,M_n,\phi_{F_n})|_{\X\setminus S}
\end{align*}
such that, after choosing a common representative
$A_n$ of the weighted instanton class associated
with $\mathfrak a_n$,
\begin{align*}
\Psi_n^*A_n\longrightarrow A,
\qquad
\Phi_n^*A_n\longrightarrow B
\end{align*}
in weighted $W^{1,p}$ on every relatively
compact subset of $\X\setminus S$.
Put $u_n:=\Phi_n^{-1}\circ\Psi_n:(F,M,\phi_F)|_{\X\setminus S}\to (G,N,\phi_G)|_{\X\setminus S}$, then $u_n^*(\Phi_n^*A_n)=\Psi_n^*A_n$. Applying the local estimate of \Cref{prop:gauge}, we obtain uniform weighted $W^{2,p}$ estimate on $u_n$ in every compact subset of $\X\setminus S$. By Sobolev embedding theorem and a diagonal subsequence argument, there exists $u:(F,M,\phi_F)|_{\X\setminus S}\to (G,N,\phi_G)|_{\X\setminus S}$ such that $u_n\to u$ in $C^0_{\mathrm{loc}}(\X\setminus S)$ and weakly in weighted $W^{2,p}_{\mathrm{loc}}$.
By passing to the limit, $u$ is a continuous unitary bundle-pair isomorphism on $\X\setminus S$. Moreover, $u^*B=A$ in the distribution sense.
Hence \Cref{prop:gauge}, applied with the
finite set $S$, extends $u$ uniquely to a weighted
$W^{2,p}$ bundle-pair isomorphism $u:(F,M,\phi_F)\xrightarrow{\sim}(G,N,\phi_G)$ on the whole $\X$ satisfying $u^*B=A$.
By the inverse correspondence $\mathfrak S=\mathfrak R^{-1}$, the two parabolic Hermitian--Einstein pairs are therefore equivalent $[(\mathcal F_*,H,\phi_F)]=[(\mathcal G_*,K,\phi_G)]$.

It remains to prove that $Z_i=Y_i$ for $i=0,1$.
Fix $z\in S$, and choose a sufficiently small
adapted ball $B_z$ such that $\overline{B_z}\cap S=\{z\}$. Let $(k_{n,z},l_{n,z})$ and $(k'_{n,z},l'_{n,z})$ be the relative topological labels associated with $\Psi_n$ and $\Phi_n$ on a collar of $\partial B_z$.
Since $u_n=\Phi_n^{-1}\circ\Psi_n\to u$ uniformly on this collar, the bundle automorphisms $u^{-1}u_n$ converge uniformly to the identity there. By taking the logarithm, $u_n$ and $u$ are homotopic
through bundle-pair isomorphisms on the collar. Since $u$ extends over the whole ball-pair
$(B_z,B_z\cap D)$, the homotopy invariance
of the relative characteristic classes in
\Cref{def:relative} gives $k_{n,z}=k'_{n,z}$ and $l_{n,z}=l'_{n,z}$ for sufficiently large $n$.
Applying condition \textup{(3)} of \Cref{def:idealconv} to the two limits, we obtain $\deg(Z_i|_{B_z})=\deg(Y_i|_{B_z})$ for $i=0,1$.

Since $z\in S$ is arbitrary, we have $Z_i=Y_i$. Hence we have $\mathfrak a=\mathfrak b$.
\end{proof}

\subsection{The topology of the moduli space}
\label{subsec:idealtopology}

In this subsection, we equip $\mathfrak I^{\mathrm{par,HE}}_\alpha(E,L;\mathcal P)$ with a topology whose convergent sequences are those of \Cref{def:idealconv}.

Let $d_{\mathrm w}$ be a metric inducing weak convergence
of finite positive Borel measures on $\X$.
For $\mathfrak a=
\bigl([(\mathcal F_*,H,\phi_F)],Z_0,Z_1\bigr)$,
put $S=\operatorname{supp}(Z_0+Z_1)$.
Choose disjoint adapted balls $B_j$ centered at the points of $S$ and balls $B_j^-\subset\subset B_j$ containing their centers. Set $K:=\X\setminus\bigcup_j B_j^-$, with $K=\X$ if $S=\varnothing$.

\begin{definition}[Uhlenbeck neighborhoods]
\label{def:idealbase}

Let $\mathfrak a=\bigl([(\mathcal F_*,H,\phi_F)],Z_0,Z_1\bigr)\in\mathfrak I^{\mathrm{par,HE}}_\alpha(E,L;\mathcal P)$, and choose $B_j^-,B_j,K$ as above. Fix a representative $A\in \mathfrak R([(\mathcal F_*,H,\phi_F)])$.

For $\varepsilon>0$, define $U(\mathfrak a;B^-,B,\varepsilon)$ to be the set of all $\mathfrak b=\bigl([(\mathcal G_*,H',\phi_G)],Y_0,Y_1\bigr)\in\mathfrak I^{\mathrm{par,HE}}_\alpha(E,L;\mathcal P)$ satisfying the following conditions:

\begin{enumerate}
\item $\operatorname{supp}(Y_0+Y_1)\subset\bigcup_{j} B_j^-$.

\item There exist representative $B\in
\mathfrak R([(\mathcal G_*,H',\phi_G)])$, and a unitary bundle-pair isomorphism $\Phi:
(F,M,\phi_F)|_{\X\setminus S}\to (G,N,\phi_G)|_{\X\setminus S}$ such that
\begin{align*}
    \|\Phi^*B-A\|_{W^{1,p}_\beta(K)}<\varepsilon.
\end{align*}

\item $d_{\mathrm w}(\nu_{\mathfrak a},\nu_{\mathfrak b})<\varepsilon$.

\item
For each $j=1,\ldots,r$, let
$(k_j,l_j)$ be the relative local
labels of $\Phi$ on a collar of
$\partial B_j$, as in \Cref{def:relative}.
Then
\begin{align*}
    \deg(Z_0|_{B_j})=\deg(Y_0|_{B_j})+k_j,\qquad \deg(Z_1|_{B_j})=\deg(Y_1|_{B_j})+k_j+l_j.
\end{align*}
\end{enumerate}
\end{definition}

\begin{lemma}\label{lem:idealbase}
The sets $U(\mathfrak a;B^-,B,\varepsilon)$ form
neighborhood bases for a first-countable topology
on $\mathfrak I^{\mathrm{par,HE}}_\alpha(E,L;\mathcal P)$.
\end{lemma}

\begin{proof}
Let $\mathfrak b\in U(\mathfrak a;B^-,B,\varepsilon)$, and let $\Phi$ be the bundle pair isomorphism on $\X\setminus S_\mathfrak{a}$ in \Cref{def:idealbase}. Write $S_{\mathfrak b}=\operatorname{supp}(Y_0+Y_1)$.
Choose disjoint adapted balls $C_k^-\subset\subset C_k$ centered at the points of $S_{\mathfrak b}$,
with each $\overline{C_k}$ contained in some $B_j^-$.
We need to show that, for sufficiently small $\delta>0$, $U(\mathfrak b;C^-,C,\delta)\subset U(\mathfrak a;B^-,B,\varepsilon)$. Let $\mathfrak c\in U(\mathfrak b;C^-,C,\delta)$, and let $\Theta$ be its comparison bundle pair isomorphism with $\mathfrak b$ on $\X\setminus S_{\mathfrak b}$.

We may let the composition $\Theta\circ\Phi$ to be the comparison isomorphism between $\mathfrak c$ and $\mathfrak a$. However, it is only defined on $\X\setminus(S_\mathfrak{a}\cup S_\mathfrak{b})$, while we need to define it on $\X\setminus S_\mathfrak{a}$. Note that the composition is defined on $K=\X\setminus\bigcup_j B_j^-$, and the weighted norm $W^{1,p}_\beta(K)$ and relative local labels on $\dd B_j$ only depends on the behaviour on $K$, we modify $\Theta\circ\Phi$ in $B_j^-$ and extend it to a bundle-pair isomorphism $\Psi$ on $\X\setminus S_{\mathfrak a}$. Since each punctured adapted ball $B_j^-\setminus\{z_j\}$ retracts, together with its divisor, onto its boundary pair.
Thus the comparison on the boundary collar
extends over the punctured ball.

Let $A,B,C$ denote the weighted connection
representatives of $\mathfrak a,\mathfrak b,
\mathfrak c$, respectively. By continuity of the weighted gauge action,
\begin{align*}
    \|\Psi^*C-A\|_{W^{1,p}_\beta(K)}=\|\Phi^*(\Theta^*C-B)+(\Phi^*B-A)\|_{W^{1,p}_{\beta}(K)}\leq C_\Phi\delta+\|\Phi^*B-A\|_{W^{1,p}_\beta(K)}.
\end{align*}
The second term is strictly less than
$\varepsilon$, so the required inequality
holds for sufficiently small $\delta$.

Let $r_0(-;B):=k_B(-)$ and $r_1(-;B):=k_B(-)+l_B(-)$. By composition and excision of relative characteristic classes,
\begin{align*}
r_i(\Psi;B_j)=r_i(\Phi;B_j)+\sum_{C_k\subset B_j^-}r_i(\Theta;C_k).
\end{align*}
Combining the identities for $\mathfrak a,\mathfrak b$ and $\mathfrak b,\mathfrak c$ gives
\begin{align*}
\deg(Z_i|_{B_j})=\deg(T_i|_{B_j})+r_i(\Psi;B_j),
\qquad i=0,1,
\end{align*}
where $T_i$ are the cycles of $\mathfrak c$.
Their supports are contained in $\bigcup_k C_k^-\subset\bigcup_j B_j^-$. Thus $\mathfrak c\in U(\mathfrak a;B^-,B,\varepsilon)$.

This proves the neighborhood-refinement
property. Two neighborhoods of the same center admit a common refinement by choosing smaller nested balls and a smaller $\varepsilon$. Hence the sets form neighborhood bases.

Since the support of each effective cycle is finite,
nested balls with radii $r_i=i^{-1}$, together with $\varepsilon_i=i^{-1}$, give a countable local base.
Thus the topology is first-countable.
\end{proof}

\begin{proposition}
\label{prop:idealconvergence}
A sequence in $\mathfrak I^{\mathrm{par,HE}}_\alpha(E,L;\mathcal P)$ converges in the topology of \Cref{lem:idealbase} if and only if it converges in the sense of \Cref{def:idealconv}.
\end{proposition}

\begin{proof}
Write $\mathfrak a_n=([A_n],Z_{0,n},Z_{1,n})$, $\mathfrak a=([A],Z_0,Z_1)$, and set $S=\operatorname{supp}(Z_0+Z_1)$.

Suppose first that $\mathfrak a_n\to\mathfrak a$ in the sense of \Cref{def:idealconv}. Let $U(\mathfrak a;B^-,B,\varepsilon)$
be a neighborhood, and put $K=\X\setminus\bigcup_j B_j^-$. The strong weighted $W^{1,p}$ convergence
in \Cref{def:idealconv} gives condition~(2),
while weak convergence of the total curvature
measures gives condition~(3). The local label identities give condition~(4).

To verify condition~(1), write $\nu_{\mathfrak a_n}=\mu_n+\sigma_n$, where $\mu_n$ is the curvature measure, and $\sigma_n=(1-2\beta)[Z_{0,n}]+2\beta[Z_{1,n}]$.
By strong local weighted $W^{1,p}$ convergence
and the weighted Sobolev multiplication estimates,
$\mu_n\to\mu$ locally in $L^1$ on $\X\setminus S$.
Since $\nu_{\mathfrak a_n}\rightharpoonup
\nu_{\mathfrak a}$ and the atomic part of $\nu_\mathfrak a$ is supported on $S$, it follows that $\int_\X\chi\,d\sigma_n\to 0$ for any $\chi\in C_c(\X\setminus S)$.
Choose a test function $0\leq\chi\leq1$ which equals to $1$ on $K$. Every atom of $\sigma_n$ has mass at least $\min\{1-2\beta,2\beta\}>0$.
Consequently, $\operatorname{supp}(Z_{0,n}+Z_{1,n})\cap K=\varnothing$ for all sufficiently large $n$.
Thus $\mathfrak a_n$ eventually belongs to
every neighborhood of $\mathfrak a$.

Conversely, suppose $\mathfrak a_n\to\mathfrak a$
in the topology of \Cref{lem:idealbase}.
Choose nested adapted balls $B_{j,m}^-\subset\subset B_{j,m}$ shrinking to each point $z_j$ of $S$, and put $K_m=\X\setminus\bigcup_j B_{j,m}^-$. By convergence, a diagonal choice gives integers $m(n)\to\infty$ such that $\mathfrak a_n\in
U(\mathfrak a;B_{m(n)}^-,B_{m(n)},1/m(n))$ for sufficiently large $n$.
Choose representatives $A_n$, $A$ and comparison
isomorphisms $\Psi_n$.
Then
\begin{align*}
\|\Psi_n^*A_n-A\|_{W^{1,p}_\beta(K_{m(n)})}<\frac1{m(n)},
\qquad
d_{\mathrm w}(\nu_{\mathfrak a_n},\nu_{\mathfrak a})<\frac1{m(n)}.
\end{align*}
Since the $K_m$ exhaust $\X\setminus S$, these inequalities give conditions~(1) and~(2)
of \Cref{def:idealconv}.

Finally, let $B$ be an adapted ball with
$\partial B\cap S=\varnothing$.
For sufficiently large $n$, each
$B_{j,m(n)}$ lies either inside $B$ or
outside $\overline B$, and all effective cycles
of $\mathfrak a_n$ are supported in the
smaller balls $B_{j,m(n)}^-$.
By excision of the relative local labels,
\begin{align*}
    k_B(\Psi_n)=\sum_{z_j\in B}k_{B_{j,m(n)}}(\Psi_n),\qquad l_B(\Psi_n)=\sum_{z_j\in B}l_{B_{j,m(n)}}(\Psi_n).
\end{align*}
Summing the identities in \Cref{def:idealbase} over $z_j\in B$ therefore gives
\begin{align*}
    \deg(Z_0|_B)&=\deg(Z_{0,n}|_B)+k_B(\Psi_n),\\
\deg(Z_1|_B)&=\deg(Z_{1,n}|_B)+k_B(\Psi_n)+l_B(\Psi_n).
\end{align*}
This is condition~(3) of \Cref{def:idealconv}, thus completing the proof.
\end{proof}

\begin{theorem}\label{thm:ideal-topology}
$\mathcal I_{\alpha}^{\mathrm{par,HE}}(E,L;\mathcal{P})$ equipped with the topology of \Cref{lem:idealbase} is a compact metrizable topological space.
\end{theorem}

\begin{proof}
By \Cref{lem:idealbase}, $\mathcal I$ is first-countable. \Cref{prop:uniquelimit,prop:idealconvergence} imply that convergent sequences have unique
limits. Thus $\mathcal I$ is Hausdorff.

We next construct a finite family of second-countable
spaces admitting continuous maps whose images cover
$\mathcal I$, which can induce a countable network. For $a=([(\mathcal F_*,H,\varphi_F)],Z_0,Z_1)\in\mathcal I$,
put $d_i=\deg Z_i$. The trace-free Chern--Weil identity give $(1-2\beta)d_0+2\beta d_1 \leq \kappa_\beta(E,L)$. Since $0<\beta<1/2$, only finitely many pairs $(d_0,d_1)\in\mathbb Z_{\geq0}^2$ can occur.
For each occurring pair $d=(d_0,d_1)$, choose a smooth determinant-fixed bundle pair $(F_d,M_d,\varphi_d)$
satisfying
\begin{align*}
c_2(F_d)=c_2(E)-d_0,\qquad \deg M_d=\deg L+d_1-d_0.
\end{align*}
Its isomorphism class is uniquely determined by $d$,
by \Cref{lem:pairclassification}.

Define
\begin{align*}
\mathcal C_d:=\left\{(Z_0,Z_1)\in\operatorname{Sym}^{d_0}(X)\times\operatorname{Sym}^{d_1}(X):Z_0|_{X\setminus D}=Z_1|_{X\setminus D}
\right\}
\end{align*}
with the subspace topology, and put
\begin{align*}
\mathcal P_d:=\mathcal M_\beta^{\mathrm{KM}}
(F_d,M_d;a_{P,\mathrm{par}})\times\mathcal C_d.
\end{align*}
The correspondence of \Cref{thm:inverse-map} defines a map
$f_d:\mathcal P_d\to\mathcal I$. The images of these
finitely many maps cover $\mathcal I$.

Each $\mathcal P_d$ is second-countable.
Indeed, the weighted instanton space is a subspace of
a separable Banach space. Its gauge quotient is
second-countable because the orbit projection is open.
The symmetric powers of the compact metrizable space
$X$ are second-countable, as are their subspaces and
finite products.

We verify that $f_d$ is continuous. Fix a point of
$\mathcal P_d$ and a neighborhood of its image
as in \Cref{def:idealbase}. Nearby cycles remain inside the prescribed smaller balls, with their multiplicities in each larger ball unchanged. For nearby connection classes, the quotient topology allows representatives converging in weighted $W^{1,p}$ on the fixed bundle pair. Using the identity comparison, all relative local labels vanish. Moreover, continuity of the curvature map into $L^p$, with $p>2$, implies $L^1$ convergence of the curvature-energy densities. Together with weak continuity of the counting measures on symmetric powers, this gives convergence of the total curvature measures. Hence all four conditions
of \Cref{def:idealbase} hold, proving continuity.

Let $\mathcal B_d$ be a countable base of $\mathcal P_d$. The family $\mathcal N
:=\bigcup_d\{f_d(V):V\in\mathcal B_d\}$ is a countable network for $\mathcal I$: if
$x\in U$ with $U$ open, choose $t\in\mathcal P_d$
such that $f_d(t)=x$, and then
$V\in\mathcal B_d$ with
$t\in V\subset f_d^{-1}(U)$.
Thus $x\in f_d(V)\subset U$.
In particular, $\mathcal I$ is Lindel\"{o}f.

By \Cref{prop:idealsequential,prop:idealconvergence},
$\mathcal I$ is sequentially compact. The sequentially-compactness and Lindel\"{o}f property implies $\mathcal I$ is compact. The compactness, Hausdorff property and the countable network implies
$\mathcal I$ is metrizable.
\end{proof}

\begin{proposition}\label{prop:open-topology}
    The moduli space $\mathcal M^{\mathrm{par,HE}}_h(E,L;\mathcal P)$ is naturally identified with the zero-defect locus $\mathcal I^0
:=\{([(\mathcal F_*,H,\varphi_F)],\varnothing,\varnothing)\in\mathcal I\}$ under the natural inclusion
\begin{align*}
    \mathcal M^{\mathrm{par,HE}}_h(E,L;\mathcal P)&\longrightarrow \mathfrak I^{\mathrm{par,HE}}_\alpha(E,L;\mathcal P)\\
    [(\mathcal E_*,H)]&\longmapsto ([(\mathcal E_*,H,\phi_E)],\varnothing,\varnothing).
\end{align*} 
Moreover, $\mathcal I^0$ is open in $\mathcal I$. Under the correspondence \Cref{thm:inverse-map}, it is also homeomorphic to
$\mathcal M_\beta^{\mathrm{KM}}
(E,L;a_{P,\mathrm{par}})$
equipped with its weighted gauge quotient topology.
\end{proposition}

\begin{proof}
By definition, every object in $\mathcal I^0$ satisfies $c_2(F)=c_2(E)$ and $\deg M=\deg L$. \Cref{lem:pairclassification} therefore identifies its underlying pair type with $(E,L)$.
Since the equivalence relations agree, the natural inclusion
is a bijection onto $\mathcal I^0$.

For $\mathfrak a\in\mathcal I^0$, its defect support is empty. Condition (i) in \Cref{def:idealbase} forces every object in the neighborhood to also have zero defects. Hence $\mathcal I^0$ is open in $\mathcal I$.

It remains to compare the induced topology with the weighted gauge quotient topology. Let $\mathcal S_\beta:=\{B\in \mathcal A^{1,p}_\beta:F_B^+=0\}$ be the weighted ASD solution space and let $\pi:\mathcal S_\beta\to \mathcal M_\beta^{\mathrm{KM}}$ be the quotient map. For $A\in\mathcal S_\beta$, the sets
\begin{align*}
V_\varepsilon([A])
:=
\pi\bigl(\{B\in\mathcal S_\beta:
\|B-A\|_{W_\beta^{1,p}}<\varepsilon\}\bigr)
\end{align*}
form a neighborhood basis at $[A]$ in the quotient topology. On the other hand, let $\mathfrak a=([A],\varnothing,\varnothing)$, then $U_\varepsilon=U(\mathfrak a;\varnothing,\varnothing,\varepsilon)$ form a neighborhood basis in the induced topology.
We claim that, for every $\varepsilon>0$, there exists $0<\delta<\varepsilon$ such that 
\begin{align*}
    V_\delta([A])\subset U_\varepsilon\subset V_\varepsilon([A]),
\end{align*}
after the identification.

Continuity of the curvature map
$W_\beta^{1,p}\to L^p$, together with $p>2$, implies $|F_B^\perp|^2\to|F_A^\perp|^2$ in $L^1$ as $B\to A$ in $W^{1,p}_\beta$.
Hence $d_w(\nu_B,\nu_A)\to0$. We may therefore choose
$0<\delta<\varepsilon$ such that
$\|B-A\|_{W_\beta^{1,p}}<\delta$ implies both
$\|B-A\|_{W_\beta^{1,p}}<\varepsilon$ and
$d_w(\nu_B,\nu_A)<\varepsilon$. This proves the first inclusion.

For the second inclusion, let $\mathfrak b=([B],\varnothing,\varnothing)\in U_\varepsilon$.
By \Cref{def:idealbase}, there is a global unitary pair isomorphism $\Phi$ such that $\|\Phi^*B-A\|_{W_\beta^{1,p}}<\varepsilon$.
Since both $B$ and $\Phi^*B$ are weighted connections, the regularity result of \Cref{prop:gauge} shows that
$\Phi$ is a weighted $W^{2,p}$ pair isomorphism.
Thus $[B]=[\Phi^*B]\in V_\varepsilon([A])$.

The two inclusions show that the induced topology on
$\mathcal I^0$ coincides with the weighted gauge quotient topology.
\end{proof}

\begin{definition}\label{def:compactification}
The \textbf{analytic compactification} $\overline{\mathcal M^{\mathrm{par,HE}}_h}(E,L;\mathcal P)$ of $\mathcal M^{\mathrm{par,HE}}_h(E,L;\mathcal P)$ is defined by the closure in $\mathfrak I^{\mathrm{par,HE}}_\alpha(E,L;\mathcal P)$.
\end{definition}

\begin{corollary}\label{cor:compactification}
The analytic compactification $\overline{\mathcal M^{\mathrm{par,HE}}_h}(E,L;\mathcal P)$ is compact and metrizable.
\end{corollary}
\begin{proof}
It is a closed subset of a compact Hausdorff metrizable space, hence also has these properties.
\end{proof}

\bibliographystyle{amsalpha}
\bibliography{ref}

\end{document}